\documentclass{amsproc}

\usepackage{amsmath}

\newtheorem{theorem}{Theorem}[section]
\newtheorem{lemma}[theorem]{Lemma}
\newtheorem{conjecture}[theorem]{Conjecture}
\newtheorem{proposition}[theorem]{Proposition}

\theoremstyle{definition}

\newtheorem{example}[theorem]{Example}

\theoremstyle{remark}
\newtheorem{remark}[theorem]{Remark}

\numberwithin{equation}{section}

\makeatletter
\def\eqnarray{%
  \stepcounter{equation}%
  \let\@currentlabel=\theequation
  \global\@eqnswtrue
  \global\@eqcnt\z@
  \tabskip\@centering
  \let\\=\@eqncr
  $$\halign to \displaywidth\bgroup\@eqnsel\hskip\@centering
  $\displaystyle\tabskip\z@{##}$&\global\@eqcnt\@ne
  \hfil$\displaystyle{{}##{}}$\hfil
  &\global\@eqcnt\tw@$\displaystyle\tabskip\z@{##}$\hfil
  \tabskip\@centering&\llap{##}\tabskip\z@\cr}
\makeatother

\def\La{\Lambda}
\def\la{\lambda}
\def\ol#1{\overline{#1}}
\def\Uq{U_q(X_n^{(1)})}
\def\Uqcl{U_q(X_n)}
\def\Uqp{U'_q(X_n^{(1)})}
\def\ch{\mbox{\sl ch}\,}
\def\lev{\mbox{\sl lev}\,}
\def\wt{\mbox{\sl wt}\,}
\def\cd{\cdots}
\def\ot{\otimes}
\def\bb{{\bf b}}
\def\bp{{\bf p}}
\def\P{{\mathcal P}}
\def\veps{\varepsilon}
\def\vphi{\varphi}
\def\C{{\mathbb C}}

\begin{document}

\title{Remarks on Fermionic Formula}

% author one information
\author{G. Hatayama}
\address{Institute of Physics, University of Tokyo, Tokyo 153-8902, Japan}
%\curraddr{}
\email{hata@gokutan.c.u-tokyo.ac.jp}
%\thanks{}

% author two information
\author{A. Kuniba}
\address{Institute of Physics, University of Tokyo, Tokyo 153-8902, Japan}
%\curraddr{}
\email{atsuo@hep1.c.u-tokyo.ac.jp}
%\thanks{}

% author three information
\author{M. Okado}
\address{Department of Informatics and Mathematical Science,
    Graduate School of Engineering Science,
Osaka University, Osaka 560-8531, Japan}
%\curraddr{}
\email{okado@sigmath.es.osaka-u.ac.jp}
%\thanks{}

% author four information
\author{T. Takagi}
\address{Department of Mathematics and Physics, National Defense Academy,
Kanagawa 239-8686, Japan}
%\curraddr{}
\email{takagi@cc.nda.ac.jp}
%\thanks{}

% author five information
\author{Y. Yamada}
\address{Department of Mathematics, Faculty of Science,
Kobe University, Hyogo 657-8501, Japan}
%\curraddr{}
\email{yamaday@math.kobe-u.ac.jp}
%\thanks{}

%\date{20/Oct/98}

\begin{abstract}
Fermionic formulae originate in the Bethe ansatz in 
solvable lattice models.
They are specific expressions of some 
$q$-polynomials as sums of products of $q$-binomial 
coefficients.
We consider the fermionic formulae associated with 
general non-twisted quantum affine algebra 
$U_q(X^{(1)}_n)$ and discuss several aspects 
related to representation theories and combinatorics.
They include crystal base theory, one dimensional sums,
spinon character formulae, $Q$-system and combinatorial 
completeness of the string hypothesis for arbitrary $X_n$.
\end{abstract}

\maketitle

\section{Introduction}

Many important $q$-polynomials and $q$-series arising in
representation theory and combinatorics can be expressed
as sums of products of $q$-binomial or multinomial
coefficients. Among them there are a class of formulae 
connected with the Bethe ansatz.

Consider the tensor product of the $L$-copies of the 
vector representation $W^{(1)}_1 \simeq
{\mathbb C}^2$ of $sl(2)$:
$$
W = W^{(1)}_1 \otimes \cdots \otimes W^{(1)}_1.
$$
Decomposing $W$ into irreducible components one gets
$$
W = \bigoplus_\la M(W,\la)V(\la),
$$
where the multiplicity $M(W,\la)$ of the 
$\la$-highest weight representation $V(\la)$ is the well known 
Kostka number $K_{\la, (1^L)}$.

In physics one regards $W$ as a Hilbert space of a 
statistical mechanical model, spin $\frac{1}{2}$ Heisenberg chain,
and tries to diagonalize an $sl(2)$-linear operator
on $W$ called the row transfer matrix (RTM).
The Bethe ansatz is a method to construct 
the eigenvectors of RTM that are $sl(2)$-highest
%%%%%%%%%%%%%%%%%%%%%%%%%%
\footnote{This is so for systems with Yangian symmetry
as in the XXX Heisenberg chain. Those systems associated with 
trigonometric solutions to the Yang-Baxter equation can have
a RTM that does not commute with
$U_q(sl(2))$ nor $sl(2)$. In such a case, the Bethe ansatz 
produces not only the highest weight vectors but
also the other weight vectors.}.
%%%%%%%%%%%%%%%%%%%%%%%%%%
%
Thus $M$ is counted as the number of the solutions to the 
Bethe equation. 
Under the string hypothesis \cite{Be} it leads to 
{\bf Fermionic formula}: $M(W,\la) = M(W,\la,q=1)$ where
$$
M(W,\la,q)=\sum_{\{m\}} q^{c(\{m\})} \prod_{j \geq 1}
\left[\begin{array}{c} p_i + m_i \\ m_i \end{array} \right]_q
%\left[p_i+m_i \atop m_i \right]_q,
$$
\begin{eqnarray*}
c(\{m\})&=&\sum_{i,j \geq 1} \min(i,j) m_i m_j-
L \sum_{i \geq 1} m_i + \frac{L(L-1)}{2}, \cr
p_i&=&L -2 \sum_{j \geq 1} \min(i,j) m_j.
\end{eqnarray*}
Here the $\{m\}$ sum is taken over
$m_i \geq 0$ ($i \geq 1$), such that
$2 \sum_{i \geq 1} i m_i=L-k$ for $\la$
corresponding to the $(k+1)$-dimensional 
representation.
The $q$-analogue $M(W,\la,q)$%
%%%%%%%%%%%%%%%%%%%%%%%%%
\footnote{$M(W,\la,q)$ here is 
$q^{L(L-1)/2}$ times that in (\ref{eq:mvq1}).}
%%%%%%%%%%%%%%%%%%%%%%%%%
is counting the 
spectrum of the momentum (cyclic lattice shift operator)
which is a specialization of the RTM.
The formula involves ($q$-)binomial without signs reflecting the 
``fermionic" nature of the counting in the Bethe ansatz.
The idea to apply the Bethe ansatz
to representation theory and combinatorics has been
initiated in \cite{KR1,KR2} and has led to fruitful results.
For example, $M(W,\la,q)$ in the above turns out to be 
the Kostka-Foulkes polynomial $K_{\la,(1^L)}(q)$ \cite{Ma}.
See \cite{Ki3} and reference therein.

On the other hand, there is another important 
idea {}from Baxter's
corner transfer matrix (CTM) method \cite{ABF,B}. 
It indicates that in the limit $L \rightarrow \infty$, 
the space $W$ has a 
structure of an integrable highest weight module 
over the quantum affine algebra $U_q(\widehat{sl}(2))$
%%%%%%%%%%%%%%%%%%%%%%%%%%%%%%%%%%%%
\footnote{Here we have the $q$-deformed XXZ chain in mind rather than
the XXX chain.}.
%%%%%%%%%%%%%%%%%%%%%%%%%%%%%%%%%%%%
%
To be more precise we invoke 
the crystal base theory of Kashiwara \cite{Kas} and regard 
$W^{(1)}_1$ as a $U_{q}'(\widehat{sl}(2))$-module,
which has a {\em perfect crystal} $B$.
Then in the limit $L \rightarrow \infty$, a subset of 
the semi-infinite tensor product $\cdots \otimes B \otimes B$
turns out to be isomorphic to the crystal base of the level 1
integrable  highest weight module over $U_{q}(\widehat{sl}(2))$
\cite{KMN1}. Moreover  $B^{\otimes L}$
with finite $L$ is isomorphic to the crystal base 
of a Demazure module \cite{KMOU}. 
By using the notion of path and energy function, one
can define {\bf One dimensional sum (1dsum)}:
$$
X=\sum_{\{b\}} q^{\sum_{j=1}^{L-1}
(L-j)H(b_j \otimes b_{j+1})},
$$
where 
$$H(1\otimes 1)=H(2 \otimes 1)=H(2 \otimes 2)=1, \quad 
H(1 \otimes 2)=0,
$$
and the $\{b\}$ sum is taken over
$b_j \in \{1,2\}$  ($1 \leq j \leq L$), such that
$\sum_{j=1}^m(\delta_{b_j,1}-\delta_{b_j,2}) \ge 0$
for $1 \leq m \leq L-1$ and $=k$ for $m=L$.
This is the branching coefficient of 
the $(k+1)$-dimensional $U_q(sl(2))$-module in the Demazure
module, which is graded by the null root $q = e^{-\delta}$. 
See $X(B^{\otimes L},\lambda,q)$ in 
(\ref{eq:def-X}) for the general definition.
When $q=1$, both $M$ and $X$ give the multiplicity
of the same representation in $W$.
In $q$ case, they are related to the spectrum of the 
RTM and CTM, respectively.
Their spectra are different in general although
they are expected to be the same 
in the ``CFT" limit, which involves $L \rightarrow \infty$ 
at least.
It is therefore remarkable that 
$$
X=M
$$
is known to hold for {\em finite} $L$
%%%%%%%%%%%%%%%%%%%%%%%%%%%%%%%%%%%%%%%
\footnote{For $q=1$ this was proved by Bethe \cite{Be} 
already in 1931.}.
%%%%%%%%%%%%%%%%%%%%%%%%%%%%%%%%%%%%%%
%
Note here that the Bethe ansatz 
play a complementary role with the CTM and the crystal base theory.
The latter provides a conceptual definition of the 1dsum $X$, 
while the former offers a specific formula $M$
which is fermionic.
This kind of phenomena have been conjectured or proved
for a variety of setting and a wide class of representations.
There are immense literatures on this, see for example 
\cite{Ber, BMS, BMSW, BLS, DF, DKKMM, FLW, FOW, FS, Ge, HKKOTY, KM, Ki3, KR1, 
KSS, NRT, OPW, SW, Wa} and the reference therein.
Especially there are extensive results on
fractional level case of $A^{(1)}_1$ cf. \cite{BMS,FLW}.
However, beyond the $A^{(1)}_1$ or $A^{(1)}_n$ case to which these 
literatures are
mostly devoted,
there are relatively few works concerning general $X^{(1)}_n$ especially at 
higher level, cf. \cite{KR2, KNS, Y}.

The aim of this paper is to treat all the  
non-twisted quantum affine algebra  $U_q(X^{(1)}_n)$ on an equal footing.
We formulate a general $X = M$ conjecture and discuss its
several consequences and applications that generalize or unify
many earlier results.
Let us give an overview of them aligned in the 
sections 2 -- 8 in the main text.

In section 2, we introduce a conjectural family of crystals 
$B^{r,s}$ of the finite dimensional representation $W^{(r)}_s$ of 
$U_q'(X^{(1)}_n)$ $(1 \le r \le n, s \in {\mathbb Z}_{\ge 1})$. 
Its existence is suggested {}from the 
Bethe ansatz argument in \cite{KR2} for the Yangian $Y(X_n)$ case.
We propose a criterion telling whether $B^{r,s}$ is perfect or not 
according to $r, s$  and the root data of $X_n$.
When 
$B^{r,s}$ is perfect, the theory of \cite{KMN1} applies and 
the semi-infinite tensor product 
$\cdots \otimes B^{r,s} \otimes B^{r,s}$ 
becomes isomorphic to the crystal base of an integrable
highest weight $U_q(X^{(1)}_n)$-module.
When $B^{r,s}$ is non-perfect, we conjecture that 
$\cdots \otimes B^{r,s} \otimes B^{r,s}$ 
is isomorphic to the crystal base of a certain tensor 
product module.
This conjecture, indicated again {}from the Bethe ansatz \cite{Ku},
has been proved in \cite{HKKOT} in some cases.

In section 3,
the 1dsums, either classically restricted one
$X(B,\la,q)$ or 
level restricted one $X_l(B,\la,q)$, are defined 
for arbitrary inhomogeneous 
$B = B^{r_1,s_1} \otimes \cdots \otimes B^{r_L,s_L}$.
By definition, the classically restricted one $X(W,\la,q)$ at $q=1$ 
coincides with the multiplicity 
$[W:\la]$ of $V(\la)$ in 
$W = W^{(r_1)}_{s_1} \otimes \cdots \otimes W^{(r_L)}_{s_L}$ regarded as
a $U_q(X_n)$-module.
Relations of  $L \rightarrow \infty$ limit of the 1sums 
(homogeneous case)
to the branching functions are also stated 
including the non-perfect case.
We shall also formulate the conjecture $X = M$, 
where the definition of the
fermionic form $M$ is postponed to section 4.

In section 4, we define the 
fermionic formula $M(W,\la,q)$ and its restricted version
$M_l(W,q)$ for general inhomogeneous 
$W = W^{(r_1)}_{s_1} \otimes \cdots \otimes W^{(r_L)}_{s_L}$.
We shall also introduce a modified one 
$N_l(W,\la,q)$, for which we allow 
negative ``vacancy numbers" $p^{(a)}_i$ in the Bethe ansatz terminology.
This is a noteworthy difference {}from $M$'s.
Nevertheless we conjecture that they are the same for $\la$
in the dominant chamber, 
see {\sc Conjecture} \ref{con:MN}.
The $N_{\infty}(W,\la,1)$ has a nice behaviour 
under the Weyl group  (cf. {\sc Remark} \ref{rem:weyl}, 
{\sc Conjecture} \ref{con:nmweyl}) and plays an essential role in section 8.

One of the important advantages of the
fermionic formulae is that
we can derive the spinon character formulae for the branching functions 
by taking the limit $L \rightarrow \infty$ 
in $M(W^{(r) \otimes L}_s,\la,q)$.
We do this in section 5, generalizing
the calculations in \cite{HKKOTY}.
In general, the spinon character formulae are
related with the particle structure of the model,
or the domain wall description of the paths \cite{NY1,NY2,ANOT,BPS,BS,BLS}.
For simply laced $X_n$, the spinon character formulae derived here
are consistent
with the conjectured particle structure of the model,
i.e. for any level $k$, the particles are labeled  by 
the fundamental representations \cite{FT} and their
$S$-matrices are the tensor products of those for
$U_q(X^{(1)}_n)$ and its RSOS analogues \cite{R}. 

In section 6, a recursion relation satisfied by
the fermionic formulae is proved.

In section 7, we discuss closely related  difference equations among 
the characters $Q^{(r)}_s = ch W^{(r)}_s$ \cite{Ki1, KR2}
($Q$-system).
We verify that the $Q$-system is indeed valid for $A_n, B_n, C_n$ and $D_n$  
as announced in \cite{KR2} and include a proof for $B_n$ case. 
In addition we establish  
some asymptotic property of the $Q^{(r)}_s$-functions needed in section 8.

Finally, in section 8, we formulate and 
prove  {\sc Theorem} \ref{th:completeness} related to the combinatorial
completeness of string hypothesis of the Bethe ansatz for arbitrary $X_n$
and $W = W^{(r_1)}_{s_1} \otimes \cdots \otimes W^{(r_L)}_{s_L}$.
By combinatorial completeness it is usually meant the equality 
$[W:\la] = M(W,\la,q=1)$, which is an indication (not a proof) 
that there are as many solutions to 
the Bethe equation as the $X_n$-highest weight vectors in $W$.
Instead of this, what our {\sc Theorem} \ref{th:completeness} 
tells is that 
$[W:\la] = N_\infty(W,\la,q=1)$ is reduced to checking 
the $Q$-system and the 
asymptotic property of the $Q^{(r)}_s$-functions
%%%%%%%%%%%%%%%%%%%%%%%%%%%%%
\footnote{In section 7 we have been able to establish these conditions 
for $A_n, B_n, C_n$ and $D_n$ only. 
But {\sc Theorem} \ref{th:completeness} itself 
works also for $E_{6,7,8}, F_4$ and $G_2$ once these conditions
are guaranteed. }.
%%%%%%%%%%%%%%%%%%%%%%%%%%%%
%
This is a support of our $X=M$ conjecture in a weak sense
in that we restrict to $q=1$ and need {\sc Conjecture} \ref{con:MN} 
to replace $M$ by $N_\infty$. 
($X = [W:\la]$ holds by definition, cf. (\ref{eq:count-mult}).)
The heart of our proof of the {\sc Theorem} is to derive 
an integral representation
of the fermionic formulae by means of the $Q$-system.
It is a generalization of the $A_n$ case in \cite{Ki2},
where {\em not} $[W:\la] = M$ but $[W:\la] = N_\infty$ was shown 
similarly.

String hypothesis is an origin of the fermionic formulae.
However we emphasize that our {\sc Theorem} \ref{th:completeness}
stands totally independent of it and elucidates the 
following points which have not necessarily been recognized 
so evidently:
(i) The combinatorial completeness is reduced to the 
$Q$-system and asymptotic property of the $Q^{(r)}_s$-functions 
for general $X_n$.
(ii) The relevant fermionic form is not $M(W,\la,1)$ but 
$N_\infty(W,\la,1)$ which contains contributions 
{}from negative vacancy numbers.
(iii) Despite the significant difference in their definitions,
$M = N_\infty$ is valid ({\sc Conjecture} \ref{con:MN}).

At present no explanation is available for the remarkable 
coincidence of the two fermionic forms 
$M$ and $N_\infty$.
Although it is not ``physical" to allow 
the negative vacancy numbers $p^{(a)}_i$
in the Bethe ansatz context, $N_\infty$ enjoys a nice symmetry 
under the Weyl group. Compare (\ref{eq:mzero}) and 
{\sc Conjecture} \ref{con:nmweyl}.

In Appendix \ref{app:hata}, we list
explicit examples of the fermionic formulae $M(W^{(r)}_s,\la,q^{-1})$,
which is a $q$-version of the list in \cite{Kl}.
Some of those data require a very long CPU time
to be evaluated.

In Appendix \ref{app:taka} we illustrate the calculation of 
the 1dsum and the fermionic form with an example {}from 
$C^{(1)}_2$.

In Appendix \ref{app:qsys} we give the explicit form of the $Q$-system for 
non-simply laced $X_n$.

\vspace{0.4cm}
\noindent
{\bf Acknowledgements} \hspace{0.1cm}
The authors thank K. Aomoto, M. Kashiwara, A.N. Kirillov, Y. Koga,
A. Lascoux, B. Leclerc, J. Lepowsky, T. Miwa, T. Nakanishi, 
M. Shimozono and J.-Y. Thibon for valuable 
discussions and interest in this 
work.
A. K. and M.O. thank N. Jing and K.C. Misra, 
organizers of the
``Conference on Affine and Quantum Affine Algebras and Related Topics''
held at North Carolina State University, Raleigh during May 21-24, 1998,
for the invitation and warm hospitality.

\section{A conjectural family of crystals}

\subsection{Preliminaries}
Let $X_n$ denote one of the classical simple Lie algebra 
$A_n (n \ge 1), B_n (n \ge 2), C_n (n \ge 2), D_n (n \ge 4), 
E_{6,7,8}, F_4$ and $G_2$, and let $X^{(1)}_n$ be the associated 
non-twisted affine Lie algebra.
We recapitulate necessary facts and notations concerning 
the quantum affine algebra $U_q(X^{(1)}_n)$. Let $\alpha_i,h_i,
\La_i$ ($i \in I=\{0,1, \ldots, n\}$) be the simple roots, simple coroots,
fundamental weights of $X^{(1)}_n$. 
We enumerate the vertices of the Dynkin diagram 
as in {\sc Table} \ref{tab:Dynkin}, which is the same with
TABLE Fin and Aff1 in \cite{Kac}.

Let $(\cdot|\cdot)$ be the standard bilinear
form normalized by $(\alpha_i|\alpha_i)=2$ with $\alpha_i$ 
a long root. We shall write $\vert x \vert^2$ to mean $(x \vert x)$.
Let $\delta$ and $c$ denote the null root and
the canonical central element, respectively. Let
$P=\bigoplus_{i \in I}{\mathbb Z}\La_i\bigoplus{\mathbb Z}\delta$ 
be the weight lattice.
We define the following subsets of $P$: $P^+=\sum_{i \in I}
{\mathbb Z}_{\ge 0}\La_i$, $P^+_l=\{\la\in P^+\mid \langle\la,c\rangle=l\}$,
$\ol{P}=\sum_{i = 1}^n{\mathbb Z}\ol{\La}_i$, $\ol{P}^+=\sum_{i=1}^{n}
{\mathbb Z}_{\ge 0}\ol{\La}_i$. 
Here $\ol{\La}_i=\La_i-\langle \La_i,c \rangle 
\La_0$ is the classical
part of $\La_i$. This map $\ol{\phantom{\La}}$ is extended to
a map on $P$ so that it is ${\mathbb Z}$-linear and 
$\ol{\delta} = 0$. 
To consider finite dimensional
$U'_q(X^{(1)}_n)$-modules, the classical weight lattice 
$P_{cl}=P/{\mathbb Z}\delta$
is also needed. 
In this paper we canonically identify $P_{cl}$ with 
$\bigoplus_{i\in I}{\mathbb Z}\La_i \subset P$.
For a precise treatment, see section 3.1 of \cite{KMN1}.
We let 
$\ol{Q} = \bigoplus_{i=1}^n {\mathbb Z} \alpha_i$ denote
the classical root lattice. 

For later convenience we introduce a few more notations 
concerning the classical Lie algebra $X_n$.
\begin{eqnarray}
t_a  & =  & \frac{2}{(\alpha_a \vert \alpha_a)} \in \{1, 2, 3\}, \quad
t = \mbox{max}_{\begin{subarray}{c}
1\le a \le n \end{subarray}} t_a = 
\begin{cases}1 & A_n, D_n, E_{6,7,8} \\
2 & B_n, C_n, F_4 \\
3 & G_2\end{cases},\label{eq:cartan}\\
C_{a b} & = & t_a(\alpha_a \vert \alpha_b), \quad 
C^{-1}_{a b} = t_a ( \ol{\La}_a \vert \ol{\La}_b) 
\quad 1 \le a, b \le n, \label{eq:cartan1}\\
\alpha_a & = & \sum_{b=1}^n C_{b a} \ol{\La}_b.
\end{eqnarray}
$C$ and $C^{-1}$ are the Cartan and the inverse Cartan matrices.
We shall write $a \sim b$ to mean that 
$C_{a b} < 0$.
The following explicit form will be of later use.
\begin{eqnarray}\label{eq:inversecartan}
( \ol{\La}_a \vert \ol{\La}_b) & = & K^{(n+1)}_{a, b} \ (A_n), 
\quad \frac{\mbox{min}(a,b)}{t_a t_b} \ (B_n), \quad
\frac{\mbox{min}(a,b)}{2} \ (C_n),\\
& & 
\left(\begin{array}{cccc}
2 & 3 & 2 & 1 \\
3 & 6 & 4 & 2 \\
2 & 4 & 3 & \frac{3}{2} \\
1 & 2 & \frac{3}{2} & 1 \end{array} \right) \ (F_4), \quad 
 \left(\begin{array}{cc}
2 & 1 \\
1 & \frac{2}{3} \end{array} \right) \ (G_2), \nonumber
\end{eqnarray}
where $K^{(n+1)}_{a, b}$ is defined by
\begin{equation}
K^{(l)}_{i, j}  =  K^{(l)}_{l-i,l-j} =
\mbox{min}(i,j) - \frac{i j}{l}.\label{eq:GKdef3}
\end{equation}

To each non-simply laced algebra $X_n = B_n, C_n, F_4$ and $G_2$,
we associate the two algebras $Z_n$ and $Y_n$ such that 
$Z_n \subset X_n \subset Y_n$.
%

%%%%%%%%%%%%%%%%%%%%%%%%%%%%%%%%%%%%%%%%%%%%%%%%%%%%%
\unitlength=.95pt
\begin{table}
\caption{Dynkin diagrams}
\label{tab:Dynkin}
\begin{tabular}[t]{rl}
$A_n$:&
\begin{picture}(106,20)(-5,-5)
\multiput( 0,0)(20,0){2}{\circle{6}}
\multiput(80,0)(20,0){2}{\circle{6}}
\multiput( 3,0)(20,0){2}{\line(1,0){14}}
\multiput(63,0)(20,0){2}{\line(1,0){14}}
\multiput(39,0)(4,0){6}{\line(1,0){2}}
\put(0,-5){\makebox(0,0)[t]{$1$}}
\put(20,-5){\makebox(0,0)[t]{$2$}}
\put(80,-5){\makebox(0,0)[t]{$n\!\! -\!\! 1$}}
\put(100,-5){\makebox(0,0)[t]{$n$}}
\end{picture}
\\
&
\\
$B_n$:&
\begin{picture}(106,20)(-5,-5)
\multiput( 0,0)(20,0){2}{\circle{6}}
\multiput(80,0)(20,0){2}{\circle{6}}
\multiput( 3,0)(20,0){2}{\line(1,0){14}}
\multiput(63,0)(20,0){1}{\line(1,0){14}}
\multiput(82.85,-1)(0,2){2}{\line(1,0){14.3}}
\multiput(39,0)(4,0){6}{\line(1,0){2}}
\put(90,0){\makebox(0,0){$>$}}
\put(0,-5){\makebox(0,0)[t]{$1$}}
\put(20,-5){\makebox(0,0)[t]{$2$}}
\put(80,-5){\makebox(0,0)[t]{$n\!\! -\!\! 1$}}
\put(100,-5){\makebox(0,0)[t]{$n$}}
\end{picture}
\\
&
\\
$C_n$:&
\begin{picture}(106,20)(-5,-5)
\multiput( 0,0)(20,0){2}{\circle{6}}
\multiput(80,0)(20,0){2}{\circle{6}}
\multiput( 3,0)(20,0){2}{\line(1,0){14}}
\multiput(63,0)(20,0){1}{\line(1,0){14}}
\multiput(82.85,-1)(0,2){2}{\line(1,0){14.3}}
\multiput(39,0)(4,0){6}{\line(1,0){2}}
\put(90,0){\makebox(0,0){$<$}}
\put(0,-5){\makebox(0,0)[t]{$1$}}
\put(20,-5){\makebox(0,0)[t]{$2$}}
\put(80,-5){\makebox(0,0)[t]{$n\!\! -\!\! 1$}}
\put(100,-5){\makebox(0,0)[t]{$n$}}
\end{picture}
\\
&
\\
$D_n$:&
\begin{picture}(106,40)(-5,-5)
\multiput( 0,0)(20,0){2}{\circle{6}}
\multiput(80,0)(20,0){2}{\circle{6}}
\put(80,20){\circle{6}}
\multiput( 3,0)(20,0){2}{\line(1,0){14}}
\multiput(63,0)(20,0){2}{\line(1,0){14}}
\multiput(39,0)(4,0){6}{\line(1,0){2}}
\put(80,3){\line(0,1){14}}
\put(0,-5){\makebox(0,0)[t]{$1$}}
\put(20,-5){\makebox(0,0)[t]{$2$}}
\put(80,-5){\makebox(0,0)[t]{$n\!\! - \!\! 2$}}
\put(103,-5){\makebox(0,0)[t]{$n\!\! -\!\! 1$}}
\put(85,20){\makebox(0,0)[l]{$n$}}
\end{picture}
\\
&
\\
$E_6$:&
\begin{picture}(86,40)(-5,-5)
\multiput(0,0)(20,0){5}{\circle{6}}
\put(40,20){\circle{6}}
\multiput(3,0)(20,0){4}{\line(1,0){14}}
\put(40, 3){\line(0,1){14}}
\put( 0,-5){\makebox(0,0)[t]{$1$}}
\put(20,-5){\makebox(0,0)[t]{$2$}}
\put(40,-5){\makebox(0,0)[t]{$3$}}
\put(60,-5){\makebox(0,0)[t]{$4$}}
\put(80,-5){\makebox(0,0)[t]{$5$}}
\put(45,20){\makebox(0,0)[l]{$6$}}
\end{picture}
\\
&
\\
$E_7$:&
\begin{picture}(106,40)(-5,-5)
\multiput(0,0)(20,0){6}{\circle{6}}
\put(40,20){\circle{6}}
\multiput(3,0)(20,0){5}{\line(1,0){14}}
\put(40, 3){\line(0,1){14}}
\put( 0,-5){\makebox(0,0)[t]{$1$}}
\put(20,-5){\makebox(0,0)[t]{$2$}}
\put(40,-5){\makebox(0,0)[t]{$3$}}
\put(60,-5){\makebox(0,0)[t]{$4$}}
\put(80,-5){\makebox(0,0)[t]{$5$}}
\put(100,-5){\makebox(0,0)[t]{$6$}}
\put(45,20){\makebox(0,0)[l]{$7$}}
\end{picture}
\\
&
\\
$E_8$:&
\begin{picture}(126,40)(-5,-5)
\multiput(0,0)(20,0){7}{\circle{6}}
\put(80,20){\circle{6}}
\multiput(3,0)(20,0){6}{\line(1,0){14}}
\put(80, 3){\line(0,1){14}}
\put( 0,-5){\makebox(0,0)[t]{$1$}}
\put(20,-5){\makebox(0,0)[t]{$2$}}
\put(40,-5){\makebox(0,0)[t]{$3$}}
\put(60,-5){\makebox(0,0)[t]{$4$}}
\put(80,-5){\makebox(0,0)[t]{$5$}}
\put(100,-5){\makebox(0,0)[t]{$6$}}
\put(120,-5){\makebox(0,0)[t]{$7$}}
\put(85,20){\makebox(0,0)[l]{$8$}}
\end{picture}
\\
&
\\
$F_4$:&
\begin{picture}(66,20)(-5,-5)
\multiput( 0,0)(20,0){4}{\circle{6}}
\multiput( 3,0)(40,0){2}{\line(1,0){14}}
\multiput(22.85,-1)(0,2){2}{\line(1,0){14.3}}
\put(30,0){\makebox(0,0){$>$}}
\put(0,-5){\makebox(0,0)[t]{$1$}}
\put(20,-5){\makebox(0,0)[t]{$2$}}
\put(40,-5){\makebox(0,0)[t]{$3$}}
\put(60,-5){\makebox(0,0)[t]{$4$}}
\end{picture}
\\
&
\\
$G_2$:&
\begin{picture}(26,20)(-5,-5)
\multiput( 0, 0)(20,0){2}{\circle{6}}
\multiput(2.68,-1.5)(0,3){2}{\line(1,0){14.68}}
\put( 3, 0){\line(1,0){14}}
\put( 0,-5){\makebox(0,0)[t]{$1$}}
\put(20,-5){\makebox(0,0)[t]{$2$}}
\put(10, 0){\makebox(0,0){$>$}}
\end{picture}
\\
&
\\
\end{tabular}
\begin{tabular}[t]{rl}
$A_1^{(1)}$:&
\begin{picture}(26,20)(-5,-5)
\multiput( 0,0)(20,0){2}{\circle{6}}
\multiput(2.85,-1)(0,2){2}{\line(1,0){14.3}}
\put(0,-5){\makebox(0,0)[t]{$0$}}
\put(20,-5){\makebox(0,0)[t]{$1$}}
\put( 6, 0){\makebox(0,0){$<$}}
\put(14, 0){\makebox(0,0){$>$}}
\end{picture}
\\
&
\\
\begin{minipage}[b]{4em}
\begin{flushright}
$A_n^{(1)}$:\\$(n \ge 2)$
\end{flushright}
\end{minipage}&
\begin{picture}(106,40)(-5,-5)
\multiput( 0,0)(20,0){2}{\circle{6}}
\multiput(80,0)(20,0){2}{\circle{6}}
\put(50,20){\circle{6}}
\multiput( 3,0)(20,0){2}{\line(1,0){14}}
\multiput(63,0)(20,0){2}{\line(1,0){14}}
\multiput(39,0)(4,0){6}{\line(1,0){2}}
\put(2.78543,1.1142){\line(5,2){44.429}}
\put(52.78543,18.8858){\line(5,-2){44.429}}
\put(0,-5){\makebox(0,0)[t]{$1$}}
\put(20,-5){\makebox(0,0)[t]{$2$}}
\put(80,-5){\makebox(0,0)[t]{$n\!\! -\!\! 1$}}
\put(100,-5){\makebox(0,0)[t]{$n$}}
\put(55,20){\makebox(0,0)[lb]{$0$}}
\end{picture}
\\
&
\\
\begin{minipage}[b]{4em}
\begin{flushright}
$B_n^{(1)}$:\\$(n \ge 3)$
\end{flushright}
\end{minipage}&
\begin{picture}(126,40)(-5,-5)
\multiput( 0,0)(20,0){3}{\circle{6}}
\multiput(100,0)(20,0){2}{\circle{6}}
\put(20,20){\circle{6}}
\multiput( 3,0)(20,0){3}{\line(1,0){14}}
\multiput(83,0)(20,0){1}{\line(1,0){14}}
\put(20,3){\line(0,1){14}}
\multiput(102.85,-1)(0,2){2}{\line(1,0){14.3}} %double line
\multiput(59,0)(4,0){6}{\line(1,0){2}} %dash line
\put(110,0){\makebox(0,0){$>$}}
\put(0,-5){\makebox(0,0)[t]{$1$}}
\put(20,-5){\makebox(0,0)[t]{$2$}}
\put(40,-5){\makebox(0,0)[t]{$3$}}
\put(100,-5){\makebox(0,0)[t]{$n\!\! -\!\! 1$}}
\put(120,-5){\makebox(0,0)[t]{$n$}}
\put(25,20){\makebox(0,0)[l]{$0$}}
\end{picture}
\\
&
\\
\begin{minipage}[b]{4em}
\begin{flushright}
$C_n^{(1)}$:\\$(n \ge 2)$
\end{flushright}
\end{minipage}&
\begin{picture}(126,20)(-5,-5)
\multiput( 0,0)(20,0){3}{\circle{6}}
\multiput(100,0)(20,0){2}{\circle{6}}
\multiput(23,0)(20,0){2}{\line(1,0){14}}
\put(83,0){\line(1,0){14}}
\multiput( 2.85,-1)(0,2){2}{\line(1,0){14.3}} %double line
\multiput(102.85,-1)(0,2){2}{\line(1,0){14.3}} %double line
\multiput(59,0)(4,0){6}{\line(1,0){2}} %dash line
\put(10,0){\makebox(0,0){$>$}}
\put(110,0){\makebox(0,0){$<$}}
\put(0,-5){\makebox(0,0)[t]{$0$}}
\put(20,-5){\makebox(0,0)[t]{$1$}}
\put(40,-5){\makebox(0,0)[t]{$2$}}
\put(100,-5){\makebox(0,0)[t]{$n\!\! -\!\! 1$}}
\put(120,-5){\makebox(0,0)[t]{$n$}}
\end{picture}
\\
&
\\
\begin{minipage}[b]{4em}
\begin{flushright}
$D_n^{(1)}$:\\$(n \ge 4)$
\end{flushright}
\end{minipage}&
\begin{picture}(106,40)(-5,-5)
\multiput( 0,0)(20,0){2}{\circle{6}}
\multiput(80,0)(20,0){2}{\circle{6}}
\multiput(20,20)(60,0){2}{\circle{6}}
\multiput( 3,0)(20,0){2}{\line(1,0){14}}
\multiput(63,0)(20,0){2}{\line(1,0){14}}
\multiput(39,0)(4,0){6}{\line(1,0){2}}
\multiput(20,3)(60,0){2}{\line(0,1){14}}
\put(0,-5){\makebox(0,0)[t]{$1$}}
\put(20,-5){\makebox(0,0)[t]{$2$}}
\put(80,-5){\makebox(0,0)[t]{$n\!\! - \!\! 2$}}
\put(103,-5){\makebox(0,0)[t]{$n\!\! -\!\! 1$}}
\put(25,20){\makebox(0,0)[l]{$0$}}
\put(85,20){\makebox(0,0)[l]{$n$}}
\end{picture}
\\
&
\\
$E_6^{(1)}$:&
\begin{picture}(86,60)(-5,-5)
\multiput(0,0)(20,0){5}{\circle{6}}
\multiput(40,20)(0,20){2}{\circle{6}}
\multiput(3,0)(20,0){4}{\line(1,0){14}}
\multiput(40, 3)(0,20){2}{\line(0,1){14}}
\put( 0,-5){\makebox(0,0)[t]{$1$}}
\put(20,-5){\makebox(0,0)[t]{$2$}}
\put(40,-5){\makebox(0,0)[t]{$3$}}
\put(60,-5){\makebox(0,0)[t]{$4$}}
\put(80,-5){\makebox(0,0)[t]{$5$}}
\put(45,20){\makebox(0,0)[l]{$6$}}
\put(45,40){\makebox(0,0)[l]{$0$}}
\end{picture}
\\
&
\\
$E_7^{(1)}$:&
\begin{picture}(126,40)(-5,-5)
\multiput(0,0)(20,0){7}{\circle{6}}
\put(60,20){\circle{6}}
\multiput(3,0)(20,0){6}{\line(1,0){14}}
\put(60, 3){\line(0,1){14}}
\put( 0,-5){\makebox(0,0)[t]{$0$}}
\put(20,-5){\makebox(0,0)[t]{$1$}}
\put(40,-5){\makebox(0,0)[t]{$2$}}
\put(60,-5){\makebox(0,0)[t]{$3$}}
\put(80,-5){\makebox(0,0)[t]{$4$}}
\put(100,-5){\makebox(0,0)[t]{$5$}}
\put(120,-5){\makebox(0,0)[t]{$6$}}
\put(65,20){\makebox(0,0)[l]{$7$}}
\end{picture}
\\
&
\\
$E_8^{(1)}$:&
\begin{picture}(146,40)(-5,-5)
\multiput(0,0)(20,0){8}{\circle{6}}
\put(100,20){\circle{6}}
\multiput(3,0)(20,0){7}{\line(1,0){14}}
\put(100, 3){\line(0,1){14}}
\put( 0,-5){\makebox(0,0)[t]{$0$}}
\put(20,-5){\makebox(0,0)[t]{$1$}}
\put(40,-5){\makebox(0,0)[t]{$2$}}
\put(60,-5){\makebox(0,0)[t]{$3$}}
\put(80,-5){\makebox(0,0)[t]{$4$}}
\put(100,-5){\makebox(0,0)[t]{$5$}}
\put(120,-5){\makebox(0,0)[t]{$6$}}
\put(140,-5){\makebox(0,0)[t]{$7$}}
\put(105,20){\makebox(0,0)[l]{$8$}}
\end{picture}
\\
&
\\
$F_4^{(1)}$:&
\begin{picture}(86,20)(-5,-5)
\multiput( 0,0)(20,0){5}{\circle{6}}
\multiput( 3,0)(20,0){2}{\line(1,0){14}}
\multiput(42.85,-1)(0,2){2}{\line(1,0){14.3}} %double line
\put(63,0){\line(1,0){14}}
\put(50,0){\makebox(0,0){$>$}}
\put( 0,-5){\makebox(0,0)[t]{$0$}}
\put(20,-5){\makebox(0,0)[t]{$1$}}
\put(40,-5){\makebox(0,0)[t]{$2$}}
\put(60,-5){\makebox(0,0)[t]{$3$}}
\put(80,-5){\makebox(0,0)[t]{$4$}}
\end{picture}
\\
&
\\
$G_2^{(1)}$:&
\begin{picture}(46,20)(-5,-5)
\multiput( 0,0)(20,0){3}{\circle{6}}
\multiput( 3,0)(20,0){2}{\line(1,0){14}}
\multiput(22.68,-1.5)(0,3){2}{\line(1,0){14.68}}
\put( 0,-5){\makebox(0,0)[t]{$0$}}
\put(20,-5){\makebox(0,0)[t]{$1$}}
\put(40,-5){\makebox(0,0)[t]{$2$}}
\put(30,0){\makebox(0,0){$>$}}
\end{picture}
\\
&
\\
\end{tabular}
\end{table}
\unitlength=1pt
%%%%%%%%%%%%%%%%%%%%%%%%%%%%%

\begin{table}[h]
\caption{}\label{tab:YZ}
\begin{center}
\begin{tabular}{c|cccc}
%\noalign{\hrule height0.8pt}
$X_n$ & $B_n$ & $C_n$ & $F_4$ & $G_2$ \\
\hline
$Y_n$ & $D_{n+1}$ & $A_{2n-1}$ & $E_6$ & $B_3$ \\
$Z_n$ & $A_1$ & $A_{n-1}$ & $A_2$ & $A_1$
\end{tabular}
\end{center}
\end{table}

Here $Z_n$ is obtained {}from $X_n$ by removing those vertices $a$
with $t_a = 1$ in the Dynkin diagram.
The embedding $X_n \hookrightarrow Y_n$ is well known.
Let $D = (D_{a b})_{1 \le a, b \le n, t_a = t_b = t}$ 
be the Cartan matrix of the subalgebra $Z_n$.
Then its inverse has the matrix elements
\begin{equation}\label{eq:dmat}
D^{-1}_{a b} = \begin{cases}
\frac{1}{2}\delta_{a n}\delta_{b n} & B_n, \\
K^{(n)}_{a, b} & C_n, \\
K^{(3)}_{a-2, b-2} & F_4,\\
\frac{1}{2}\delta_{a 2}\delta_{b 2} & G_2. \end{cases}
\end{equation}

Given a positive integer $l$ we consider the following
subsets of ${\mathbb Z} \times {\mathbb Z}$:
\begin{eqnarray}
H_l & = & \{(a,j) \mid 1 \le a \le n, 1 \le j \le t_al\}, 
\label{eq:GKdef1}\\
{\overline H}_l & = &\{(a,j) \mid 1 \le a \le n, 1 \le j \le t_al-1\},
\label{eq:GKdef2}\\
H_l[i] & = & \{(a,j) \mid 1 \le a \le n, \frac{t_a}{t}(i-1) < j \le t_al \}
\quad 1 \le i \le tl+1.\label{eq:Hdef1} 
\end{eqnarray}
Thus 
$H_l[tl+1]= \emptyset$ and $H_l[1] = H_l$.  One 
can check
\begin{equation}\label{eq:Hdef2}
H_l[i] = H_l[i+1] \sqcup \{(a,\frac{t_a}{t}i) \mid
1 \le a \le n, \frac{t_a}{t}i \in {\mathbb Z} \}.
\end{equation}

\subsection{Terminology for crystals}

A crystal base $B$ of a $\Uq$ ($\Uqp,\Uqcl$)-module can be regarded as
a set of basis vectors of the module at $q=0$. On $B$ one still has Chevalley-like
generators 
\[
{\tilde e}_i,{\tilde f}_i:\;
B\longrightarrow B\sqcup\{0\},
\]
which are sometimes called Kashiwara operators. For $b\in B$ we set
$\veps_i(b)=\max\{k\ge0\mid{\tilde e}_i^kb\ne0\},
\vphi_i(b)=\max\{k\ge0\mid{\tilde f}_i^kb\ne0\}$.
If $B_1$ and $B_2$ are crystals, the crystal structure on the tensor product
$B_1\ot B_2$ is given by
\begin{eqnarray}
{\tilde e}_i(b_1\ot b_2)&=&
\begin{cases}
{\tilde e}_ib_1\ot b_2&\mbox{if }\vphi_i(b_1)\ge\veps_i(b_2),\\
b_1\ot {\tilde e}_ib_2&\mbox{otherwise},
\end{cases}\label{eq:tensor_e}\\
{\tilde f}_i(b_1\ot b_2)&=&
\begin{cases}
{\tilde f}_ib_1\ot b_2&\mbox{if }\vphi_i(b_1)>\veps_i(b_2),\\
b_1\ot {\tilde f}_ib_2&\mbox{otherwise}.
\end{cases}\label{eq:tensor_f}
\end{eqnarray}

We mainly use two categories of crystals. The first one contains the crystal
base $B(\la)$ of the irreducible integrable $\Uq$-module 
$\mathcal{L}(\la)$ with highest weight
$\la\in P^+$. $B(\la)$ is a $P$-weighted crystal.
The other one contains a crystal base $B$ of a finite-dimensional
$\Uqp$-module. As opposed to $B(\la)$, $B$ is a finite set. It is 
$P_{cl}$-weighted. We shall call it a {\em finite crystal}. 
For a finite crystal $B$, we set
$\veps(b)=\sum_i\veps_i(b)\La_i$, $\vphi(b)=\sum_i\vphi_i(b)\La_i$,
$\wt b = \vphi(b) - \veps(b)$, 
and introduce the {\em level} of $B$ by
\[
\lev B=\min\{\langle c,\veps(b)\rangle\mid b\in B\}.
\]
We further set $B_{\min}=\{b\in B\mid\langle c,\veps(b)\rangle=\lev B\}$,
and call an element of $B_{\min}$ {\em minimal}.
A $\Uqp$-module can be viewed as a $\Uqcl$-module.
For the latter the irreducible representation with 
highest weight $\la \in \ol{P}^+$ will be denoted by $V(\la)$.
\subsection{Family $W^{(r)}_s$}

We shall present a conjectural family of finite-dimensional
$\Uqp$-modules having crystal bases. This is a revelation of 
the Bethe Ansatz.

\begin{conjecture}
For each $(r,s)$ ($1\le r\le n,s\ge1$), there exists an irreducible
finite-dimensional $\Uqp$-module $W^{(r)}_s$ having the following features:
\begin{itemize}
\item[(1)] $W^{(r)}_s$ has a crystal base $B^{r,s}$. $B^{r,s}$ is a 
           finite crystal of level 
$\left\lceil\frac{s}{t_r}\right\rceil$
%%%%%%%%%%%%%%%%%%%%%%%%%%%%%%%%%%%%%%%%%%%%%%%%%%%%%
\footnote{The symbol $\lceil x \rceil$ denotes the smallest integer
not less than $x$.}.
%%%%%%%%%%%%%%%%%%%%%%%%%%%%%%%%%%%%%%%%%%%%%%%%%%%%%
           Moreover, it is perfect if $\frac{s}{t_r}$
           is an integer, and not perfect if not an integer.
\item[(2)] As a $\Uqcl$-module, $W^{(r)}_s$ decomposes itself into
\[
W^{(r)}_s=\bigoplus_{\la\in \ol{P}^+}M(W^{(r)}_s,\la,q=1)\,V(\la),
\]
where $M(W^{(r)}_s,\la,q)$ is defined in (\ref{eq:mvq1}).
\item[(3)] Set $Q^{(r)}_s=\ch W^{(r)}_s$, then $Q^{(r)}_s$ satisfies
           the $Q$-system (\ref{eq:qsys}).
\end{itemize}
\end{conjecture}

{}From (\ref{eq:mvq1}) the decomposition (2) of $W^{(r)}_s$ has the form
\[
W^{(r)}_s=V(s\ol{\La}_r)\oplus\cdots,
\]
where $\cdots$ contains irreducible modules with highest weights
strictly lower than $s\ol{\La}_r$ only. 
For $X_n=A_n$, $W^{(r)}_s=V(s\ol{\La}_r)$,
i.e. representation corresponding to a rectangular Young diagram. 
In the other cases, we need the $\cdots$ part  in
general. (See \cite{CP1,KR2} for the Yangian case.) 
A list of such decompositions is available in Appendix A (by setting
$q=1$).

The notion of perfect crystals is introduced in \cite{KMN1}. {}From a perfect 
crystal $B$ we can construct a set of paths, which is isomorphic to the 
crystal of an irreducible integrable $\Uq$-module. In \cite{HKKOT} the
definition of a set of paths is generalized to non-perfect crystals. 
Also in this case, the set of paths is isomorphic to the crystal of an
integrable module, but not necessarily irreducible as we will see next.

\begin{remark}
For $v \in {\mathbb C}^{\times}$, let $W^{(r)}_s(v)$ 
denote the pull back of $W^{(r)}_s$ by
the Hopf algebra automorphism $\tau_v$ 
defined on the Drinfel'd generators \cite{CP2}.
Then we expect that 
$W^{(r)}_s(v)$ is the irreducible $U'_q(X^{(1)}_n)$-module  
characterized by the Drinfel'd polynomial 
$P_a(u) = \left((1-q^{\frac{s-1}{t_a}}uv)(1-q^{\frac{s-3}{t_a}}uv)
\cdots (1-q^{\frac{-s+1}{t_a}}uv)\right)^{\delta_{a, r}}$ 
up to a normalization of $v$.
\end{remark}

\begin{remark}
To our knowledge, the following crystals have 
been shown to be perfect for 
the non-twisted case so far.
\[
\begin{array}{ll}
\mbox{\cite{KMN2}} &A_n:B^{r,s},\,B_n:B^{1,s},\,C_n:B^{n,s},\,
                    D_n:B^{1,s}B^{n-1,s}B^{n,s},\\
\mbox{\cite{KKM}} &C_n:B^{1,2s},\\
\mbox{\cite{Ko}} &B_n:B^{r,t_r},\,C_n:B^{r,2}(r\ne n),\,
                  D_n:B^{r,1}(r\ne n-1,n),\\
\mbox{\cite{Ya}} &G_2:B^{1,s}.
\end{array}
\]
\end{remark}

We give some examples of $B^{r,s}$.

\begin{example}\label{ex:crystalstructure-A}
For $X_n=A_n$, $B^{1,s}$ is isomorphic to $B(s\ol{\La}_1)$
as a crystal for $U_q(A_n)$. As a set,
\begin{equation*}
B^{1,s}=\{ (x_1,\dots ,x_{n+1}) \in {\mathbb Z}^{n+1} |\,
x_i \ge 0,\, \sum_{i=1}^{n+1} x_i=s \}.
\end{equation*}
The actions of $\tilde{e}_i$, $\tilde{f}_i$ and 
$\varepsilon,\varphi$
are defined as follows (cf. e.g. \cite{HKKOTY}):\newline
For $b=(x_1,\dots ,x_{n+1}) \in B^{1,s}$,
{\allowdisplaybreaks %
\begin{align*}
\tilde{e}_0 b &= (x_1-1,x_2,\dots,x_{n+1}+1), \\
\tilde{f}_0 b &= (x_1+1,x_2,\dots,x_{n+1}-1), \\
\tilde{e}_i b &= (x_1,\dots,x_i+1,x_{i+1}-1,\dots,x_{n+1}) 
\text{ for } i=1,\dots,n, \\
\tilde{f}_i b &= (x_1,\dots,x_i-1,x_{i+1}+1,\dots,x_{n+1})
\text{ for } i=1,\dots,n,
\end{align*}
}
where the RHS is regarded as 0 if it is not an element of $B^{1,s}$, and
\[
\varepsilon(b) = \sum_{i=0}^{n} x_{i+1} \Lambda_i,\,\,
\varphi(b) = x_{n+1} \Lambda_0 + \sum_{i=1}^{n} x_i \Lambda_i.
\]
In this case, we have $(B^{1,s})_{\min}=B^{1,s}$.
$B^{1,s}$ is perfect of level $s$.
\end{example}

\begin{example}\label{ex:crystalstructure-C-row}
For $X_n=C_n$, $B^{1,s}$ is isomorphic to 
$B(s\ol{\Lambda}_1) \oplus B((s-2)\ol{\Lambda}_1) \oplus \cdots 
\oplus (B(\ol{\Lambda}_1)\mbox{ or }B(0))$
as a crystal for $U_q(C_n)$. As a set,
\begin{equation*}
\begin{split}
B^{1,s}&=\{ (x_1,\dots ,x_{n},\bar{x}_{n},\dots,\bar{x}_1) \in 
{\mathbb Z}^{2n} |\,
x_i, \bar{x}_i \ge 0,\, \\ &\quad\quad
\sum_{i=1}^{n} (x_i+\bar{x}_i) \le s,\,
\sum_{i=1}^{n} (x_i+\bar{x}_i) \equiv s \pmod{2} \}.
\end{split}
\end{equation*}
The actions of $\tilde{e}_i$, $\tilde{f}_i$ and
$\varepsilon,\varphi$ are defined 
as follows (cf. \cite{KKM,HKKOT}):\newline
For $b=(x_1,\dots ,x_{n},\bar{x}_{n},\dots,\bar{x}_1) \in B^{1,s}$,
{\allowdisplaybreaks %
\begin{align*}
\tilde{e}_0 b &= \begin{cases}
	(x_1-2,x_2,\dots,\bar{x}_{2},\bar{x}_{1}) 
        &\text{ if } x_1 \ge \bar{x}_1+2,\\
	(x_1-1,x_2,\dots,\bar{x}_2,\bar{x}_1+1) 
        &\text{ if } x_1=\bar{x}_1+1, \\
	(x_1,x_2,\dots,\bar{x}_2,\bar{x}_1+2) &\text{ if } x_1 \le \bar{x}_1,
\end{cases} \\
\tilde{f}_0 b &= \begin{cases}
	(x_1+2,x_2,\dots,\bar{x}_{2},\bar{x}_{1}) 
        &\text{ if } x_1 \ge \bar{x}_1,\\
	(x_1+1,x_2,\dots,\bar{x}_2,\bar{x}_1-1) 
        &\text{ if } x_1=\bar{x}_1-1, \\
	(x_1,x_2,\dots,\bar{x}_2,\bar{x}_1-2) &\text{ if } x_1 \le \bar{x}_1-2,
\end{cases} \\
\tilde{e}_i b &= \begin{cases}
	(x_1,\dots,x_i+1,x_{i+1}-1,\dots,\bar{x}_{1})
	\text{ if } x_{i+1} > \bar{x}_{i+1}\\
	(x_1,\dots,\bar{x}_{i+1}+1,\bar{x}_{i}-1,\dots,\bar{x}_{1})
	\text{ if } x_{i+1} \le \bar{x}_{i+1}
	\end{cases}
\text{ for } i=1,\dots,n-1,\\
\tilde{f}_i b &= \begin{cases}
	(x_1,\dots,x_i-1,x_{i+1}+1,\dots,\bar{x}_{1})
	\text{ if } x_{i+1} \ge \bar{x}_{i+1}\\
	(x_1,\dots,\bar{x}_{i+1}-1,\bar{x}_{i}+1,\dots,\bar{x}_{1})
	\text{ if } x_{i+1} < \bar{x}_{i+1}
	\end{cases}
\text{ for } i=1,\dots,n-1,\\
\tilde{e}_n b &= (x_1,\dots,x_n+1,\bar{x}_n-1,\dots,\bar{x}_1), \\
\tilde{f}_n b &= (x_1,\dots,x_n-1,\bar{x}_n+1,\dots,\bar{x}_1), \\
\end{align*}
where the RHS is regarded as 0 if it is not an element of $B^{1,s}$, and
\begin{align*}
\varepsilon(b) &= 
\left( \frac{s-\sum_{i=1}^{n}(x_i+\bar{x}_i)}{2}
	+ (x_1-\bar{x}_1)_{+} \right) \Lambda_0 +
\sum_{i=1}^{n-1} \left( \bar{x}_i + (x_{i+1}-\bar{x}_{i+1})_{+}
	\right) \Lambda_i +
\bar{x}_n \Lambda_n,\\
\varphi(b) &= 
\left( \frac{s-\sum_{i=1}^{n}(x_i+\bar{x}_i)}{2}
	+ (\bar{x}_1-x_1)_{+} \right) \Lambda_0 +
\sum_{i=1}^{n-1} \left( x_i + (\bar{x}_{i+1}-x_{i+1})_{+}
	\right) \Lambda_i +
x_n \Lambda_n,\\
\end{align*}
}
where $(x)_{+}:=\max(x,0)$.
In this case, we have
\[
\left( B^{1,s} \right)_{\min} = \begin{cases}
\{ (x_1,\dots,x_n,x_n,\dots,x_1) \in B^{1,s} \}
\text{ if $s$ is even,}\\[5pt]
\{ (x_1,\dots,x_k,\dots,x_n,x_n,\dots,x_k \pm 1,\dots,x_1 )
	\in B^{1,s}\,|\,k=1,\dots,n \}  \\
	\quad\quad \text{ if $s$ is odd.}
\end{cases}
\]
$B^{1,s}$ is perfect (level $\frac{s}2$) if $s$ is even, and 
non-perfect (level $\frac{s+1}2$) if $s$ is odd.
\end{example}

\begin{example}\label{ex:crystalstructure-C-col}
For $X_n=C_n$, $B^{r,1}$ is isomorphic to 
$B(\ol{\Lambda}_r)$ as a crystal for $U_q(C_n)$. As a set,
\begin{equation*}
\begin{split}
B^{r,1}&=\{ (x_1,\dots ,x_{n},\bar{x}_{n},\dots,\bar{x}_1) \in 
{\mathbb Z}^{2n} |\,
x_i, \bar{x}_i =0 \text{ or }1,\\
&\quad\quad
\sum_{i=1}^{n} (x_i+\bar{x}_i) =r,\, 
\text{ if $x_k=\bar{x}_k=1$ then $\sum_{i=1}^k(x_i+\bar{x}_i)\le k$}
\}.
\end{split}
\end{equation*}
The actions of $\tilde{e}_i$ and $\tilde{f}_i$
are defined as follows (cf. \cite{AK}):\newline
For $b=(x_1,\dots ,x_{n},\bar{x}_{n},\dots,\bar{x}_1) \in B^{r,1}$,
{\allowdisplaybreaks %
\begin{align*}
\tilde{e}_0 b &=
	(x_1-1,x_2,\dots,\bar{x}_2,\bar{x}_1+1),\\
\tilde{f}_0 b &=
	(x_1+1,x_2,\dots,\bar{x}_2,\bar{x}_1-1), \\
\tilde{e}_i b &= \begin{cases}
	(x_1,\dots,\bar{x}_{i+1}+1,\bar{x}_{i}-1,\dots,\bar{x}_{1})\\
	\quad\quad\text{ if $(\bar{x}_{i+1},\bar{x}_{i})=(0,1)$ and
	$(x_i,x_{i+1}) \ne (1,0)$}\\[5pt]
	(x_1,\dots,x_{i}+1,x_{i+1}-1,\dots,\bar{x}_{1})\\
	\quad\quad\text{ if $(\bar{x}_{i+1},\bar{x}_{i}) \ne (0,1)$ and
	$(x_i,x_{i+1}) = (0,1)$}\\[5pt]
	0 \text{ otherwise }
	\end{cases}
\text{ for } i=1,\dots,n-1,\\
\tilde{f}_i b &= \begin{cases}
	(x_1,\dots,x_i-1,x_{i+1}+1,\dots,\bar{x}_{1})\\
	\quad\quad\text{ if $(x_{i},x_{i+1})=(1,0)$ and
	$(\bar{x}_{i+1},\bar{x}_{i}) \ne (0,1)$}\\[5pt]
	(x_1,\dots,\bar{x}_{i+1}-1,\bar{x}_{i}+1,\dots,\bar{x}_{1})\\
	\quad\quad\text{ if $(x_{i},x_{i+1}) \ne (1,0)$ and
	$(\bar{x}_{i+1},\bar{x}_{i}) = (1,0)$}\\[5pt]
	0 \text{ otherwise }
	\end{cases}
\text{ for } i=1,\dots,n-1,\\
\tilde{e}_n b &= (x_1,\dots,x_n+1,\bar{x}_n-1,\dots,\bar{x}_1), \\
\tilde{f}_n b &= (x_1,\dots,x_n-1,\bar{x}_n+1,\dots,\bar{x}_1),
\end{align*}}%
where the RHS is regarded as 0 if it is not an element of $B^{r,1}$.
$B^{r,1}$ is not perfect except when $r=n$. The level of $B^{r,1}$ is 1.
\end{example}

\subsection{Set of paths} \label{subsec:paths}

Let $B$ be a finite crystal of level $k$. {}From $B$ we construct a 
subset of $\cd\ot B\ot\cd\ot B$ called a {\em set of paths}. 
See also \cite{HKKOT}. First we fix a reference path 
$\bp=\cd\ot\bb_j\ot\cd\ot\bb_2\ot\bb_1$. For any $j$, $\veps(\bb_j)$ should
have level $k$, and satisfy 
\begin{equation} \label{eq:phi=eps}
\vphi(\bb_{j+1})=\veps(\bb_j).
\end{equation}
Set 
\begin{equation} \label{eq:set-of-paths}
\P(\bp,B)=\{p=\cd\ot b_j\ot\cd\ot b_2\ot b_1\mid 
b_j\in B,b_J=\bb_J\mbox{ for }J\gg1\}.
\end{equation}
An element of $\P(\bp,B)$ is called a {\em path}. To a path we can associate
a weight in $P$. For this purpose we consider the energy function
$H:\,B\times B\rightarrow {\mathbb Z}$, 
which is determined up to constant difference 
by the following rule:
\begin{eqnarray*}
H({\tilde e}_i(b_1\ot b_2))&=H(b_1\ot b_2)+1
&\mbox{ if }i=0,\vphi_0(b_1)\geq\veps_0(b_2),\\
&=H(b_1\ot b_2)-1
&\mbox{ if }i=0,\vphi_0(b_1)<\veps_0(b_2),\\
&\hspace{-6mm}=H(b_1\ot b_2)&\mbox{ if }i\neq0.
\end{eqnarray*}
Using this function we define the energy $E(p)$ and weight $\wt p$ of
$p\in\P(\bp,B)$ by 
\begin{eqnarray*}
E(p)&=&\sum_{j=1}^\infty j(H(b_{j+1}\ot b_j)-H(\bb_{j+1}\ot\bb_j)),\\
\wt p&=&\vphi(\bb_1)+\sum_{j=1}^\infty(\wt b_j-\wt\bb_j)-E(p)\delta.
\end{eqnarray*}
Notice that they have finite values. To calculate the weight of an element
of $B$, we identify $P_{cl}$ with $\bigoplus_{i\in I}{\mathbb Z}\La_i$.

Let us start with examining $\P(\bp,B)$ when $B$ is perfect of level $k$.
Take $\la\in P^+_k$ and choose a unique $\bb_1$ such that $\vphi(\bb_1)=\la$.
Then all $\bb_j$ is uniquely fixed {}from (\ref{eq:phi=eps}). We denote this
reference path by $\bp^{(\la)}$. It is known in \cite{KMN1} that we have 
the following isomorphism of $P$-weighted crystals.
\[
\P(\bp^{(\la)},B)\simeq B(\la).
\]
Of course, the actions of ${\tilde e}_i,{\tilde f}_i$ on the left hand side
obey the tensor-product rule of crystals (\ref{eq:tensor_e},\ref{eq:tensor_f}).
``Signature rule'' is a good way to calculate the actions on multiple tensor
products. See \cite{KMOU} for that.

\begin{remark}
The set of paths $\P(\bp,B)$ admits a natural filtration
\begin{eqnarray*}
&&\P_0(\bp,B)\subset\P_1(\bp,B)\subset\cd\subset\P_L(\bp,B)
\subset\cd\subset\P(\bp,B),\\
&&\P_L(\bp,B)=\{p\in\P(\bp,B)\mid b_j=\bb_j\mbox{ for }j>L\}.
\end{eqnarray*}
If $B$ is perfect and $\la=k\La_0$, the finite subset $\P_L(\bp^{(\la)},B)$
can be identified with the crystal base corresponding to a Demazure module
$\mathcal{L}_w(\la)$ for a special affine Weyl group element $w$.
See \cite{KMOU,KMOTU1} for details.
\end{remark}

\begin{example}
We give examples of $\bp^{(\la)}=\cd\ot\bb^{(\la)}_j\ot\cd\ot\bb^{(\la)}_1$.
Set $\la=\sum_{i=0}^n\la_i\La_i$. 

For $X_n=A_n$ and $B=B^{1,s}$ 
(cf. {\sc Example} \ref{ex:crystalstructure-A}),
\[
\bb^{(\la)}_j=(\la_{(j)},\la_{(j+1)},\cd,\la_{(j+n)})\in B^{1,s}.
\]
Here $(a)$ denotes the integer such that $(a)\equiv a$ (mod $n+1$),
$0\le(a)\le n$.

For $X_n=C_n$ and $B=B^{1,s}$ ($s$: even) (cf. {\sc Example} 
\ref{ex:crystalstructure-C-row}),
\[
\bb^{(\la)}_j=(\la_1,\cd,\la_n,\la_n,\cd,\la_1)\in B^{1,s}.
\]
\end{example}

Next we consider the non-perfect case. Assume $B$ 
is non-perfect and its level is $k$. Then even if we fix $\bb_1$, the
reference path is not unique. There are actually 
infinitely many reference
paths. Among them, we conjecture the existence of $\bp^{(\la,\mu)}$ which is
friendly to representation theory in the following sense.
For $B=B^{r,s}$ set $B^\dagger=B^{r,1}$. For our $X_n$ recall $Y_n$ in 
{\sc Table} \ref{tab:YZ}.
We take $\la$ {}from $P^+_{k-1}$ of $X^{(1)}_n$, but $\mu$ {}from $P^+_1$ of 
$Y^{(1)}_n$. Then there
exist reference paths $\bp^{(\la,\mu)}$ for $B$ and $\bp^{\dagger(\mu)}$
for $B^\dagger$ and the following isomorphism.
\[
\P(\bp^{(\la,\mu)},B)\simeq B(\la)\ot\P(\bp^{\dagger(\mu)},B^\dagger).
\]
Moreover, at the character level, $\P(\bp^{\dagger(\mu)},B^\dagger)$ 
can be identified with the $Y^{(1)}_n$-module $\mathcal{L}^Y(\mu)$ considered as 
an $X^{(1)}_n$-module through the natural embedding 
$X^{(1)}_n\hookrightarrow Y^{(1)}_n$. This is proven for 
$U'_q(C^{(1)}_n)$-crystal $B^{1,s}$ ($s$: odd) in \cite{HKKOT}.

\begin{example}
In this $U'_q(C^{(1)}_n)$-crystal $B^{1,s}$ ($s$: odd) case,
$\bp^{(\la,\mu)}=\cd\ot\bb^{(\la,\mu)}_j\ot\cd\ot\bb^{(\la,\mu)}_1$,
$\bp^{\dagger(\mu)}=\cd\ot\bb^{\dagger(\mu)}_j\ot\cd\ot\bb^{\dagger(\mu)}_1$
are given as follows. Set $\la=\sum_{i=0}^n\la_i\La_i\in P^+_{\frac{s-1}2}$
of $C^{(1)}_n$ and take $\mu=\La_0$ of $A^{(1)}_{2n-1}$.
\begin{eqnarray*}
\bb^{(\la,\La_0)}_j&=&\begin{cases}
(\la_1,\cd,\la_n,\la_n,\cd,\la_i+1,\cd,\la_1)
& \mbox{if }j\equiv i\,(\mbox{mod }2n)\\
& \mbox{ for some }i\,(1\le i\le n)\\
(\la_1,\cd,\la_i+1,\cd,\la_n,\la_n,\cd,\la_1)
& \mbox{if }j\equiv 1-i\,(\mbox{mod }2n)\\
& \mbox{ for some }i\,(1\le i\le n),
\end{cases}\\
\bb^{\dagger(\La_0)}_j&=&\bb^{(0,\La_0)}_j
\quad(\mbox{i.e. set formally $\la=0$ in the above}).
\end{eqnarray*}
Thus we have the isomorphism 
\[
\P(\bp^{(\la,\La_0)},B^{1,s})\simeq B(\la)\ot\P(\bp^{\dagger(\La_0)},B^{1,1}),
\]
and the character of $\P(\bp^{\dagger(\La_0)},B^{1,1})$ agrees with that of
the $A^{(1)}_{2n-1}$-module $\mathcal{L}^Y(\La_0)$ 
considered as a $C^{(1)}_n$-module.
\end{example}

\section{One dimensional sums}

Fix any $B=B^{r,s}$ and let $k$ be the level of $B$. Take $b_0\in B$ 
such that $\vphi(b_0)=k\La_0$. In all known cases, such $b_0$ is 
unique for $B$.
For $\la\in \ol{P}^+$, we define
\begin{equation} \label{eq:def-X}
X(B^{\ot L},\la,q)=\mathop{{\sum}^*}q^{\sum_{j=0}^{L-1}(L-j)H(b_j\ot b_{j+1})},
\end{equation}
where the sum $\sum^*$ is taken over all 
$b_1\ot\cd\ot b_L\in B^{\ot L}$ satisfying $\sum_{j=1}^L\wt b_j=\la$
and
\begin{equation} \label{eq:hwc}
{\tilde e}_i(b_1\ot\cd\ot b_L)=0\quad\forall i\neq0.
\end{equation}
Note that $b_0$ is fixed. (\ref{eq:hwc}) can be rewritten 
as local conditions.
\[
{\tilde e}_i^{\langle h_i,\wt b_1+\cd+\wt b_{j-1}\rangle+1}b_j=0
\quad\forall i\neq0,1\le j\le L.
\]
(\ref{eq:def-X}) is introduced in \cite{KMOTU2} in more general setting
($b_0$ not limited to the above one) 
and called ``classically restricted 1dsum''. 
It is easy to see {}from the definition that
\begin{equation} \label{eq:count-mult}
X(B^{\ot L},\la,1)=[W^{\ot L}:\la],
\end{equation}
where $W$ is the corresponding finite-dimensional $\Uqp$-module of $B$
and 
\begin{equation}\label{eq:cl-mult}
[M:\mu]=\dim\C\langle v\in M\mid\wt v=\mu,e_iv=0\,\forall i\neq0\rangle.
\end{equation}
If $B$ is perfect of level $k$, it is known \cite{KMOTU2} 
that for $\la\in\ol{P}^+$,
\begin{equation} \label{eq:bran-func}
\lim_{\begin{subarray}{c}
L\to\infty\\ L\equiv0\,\scriptstyle{mod}\,\gamma
\end{subarray}}q^{-c_L}
X(B^{\ot L},\la,q)=\sum_i[\mathcal{L}(k\La_0):\la+k\La_0-i\delta]\,q^i.
\end{equation}
Here $c_L=\sum_{j=0}^{L-1}(L-j)H(\bb_j\ot\bb_{j+1})$, $\bb_0,\cd,\bb_L$
are the last $L+1$ components of $\bp^{(k\La_0)}$ (read {}from left to right), 
and $\gamma$ is the period of the sequence $\bb_1,\bb_2,\cd$. Note that
when $L\equiv0$ mod $\gamma$, $\bb_0$ agrees with the fixed $b_0$. Note also 
that in $[\mathcal{L}(k\La_0):\mu]$, $\mu$ is a weight in $P$ as opposed to $\la$ in
(\ref{eq:count-mult}) in $\ol{P}$.
The right hand side of (\ref{eq:bran-func}) is called a branching function.
If $B$ is non-perfect, 
the conjecture in section \ref{subsec:paths} implies
\begin{equation}\label{eq:bran-func2}
\lim_{\begin{subarray}{c}
L\to\infty\\ L\equiv0\,\scriptstyle{mod}\,\gamma
\end{subarray}}q^{-c_L}
X(B^{\ot L},\la,q)
=\sum_i[\mathcal{L}((k-1)\La_0)\ot \mathcal{L}^Y(\La_0):\la+k\La_0-i\delta]\,q^i.
\end{equation}
Here $\mathcal{L}^Y(\La_0)$ is the integrable $Y^{(1)}_n$-module with highest weight 
$\La_0$ regarded as an $X^{(1)}_n$-module.

We are to state a conjecture for the expression of $X$.

\begin{conjecture}\label{con:XM}
Let $B=B^{r,s}$, set $k=\left\lceil\frac{s}{t_r}\right\rceil$ and take
$b_0,\, b^\natural \in B$ such that 
$\vphi(b_0)=k\La_0,\, \wt b^\natural = s\ol{\La}_r$. Then we have 
\[
X(B^{\ot L},\la,q)=q^c M(W^{(r)\ot L}_s,\la,q)
\]
with  $c = \frac{L(L-1)}{2}H(b^\natural \otimes b^\natural) + 
L H(b_0 \otimes b^\natural)$, where $M$ is defined in 
(\ref{eq:mvq1}) -- (\ref{eq:mvq4}).
\end{conjecture}

Surprisingly, this conjecture admits ``inhomogeneous'' version.
Consider the tensor product of crystals
\[
B_0\ot B^{r_1,s_1}\ot\cd\ot B^{r_L,s_L}.
\]
Here $B_0=B^{r_{i_1},s_{i_1}}\ot\cd\ot B^{r_{i_t},s_{i_t}}$
and $(r_{i_1},s_{i_1}),\cd,(r_{i_t},s_{i_t})$ are mutually distinct
elements of $\{(r_1,s_1),\cd,(r_L,s_L)\}$. To explain the corresponding
classically restricted 1dsum, one has to redefine the energy function
$H$ due to inhomogeneity \cite{NY3,HKKOTY}.

Suppose $B_1$ and $B_2$ are finite crystals, and $b_1\ot b_2\in
B_1\ot B_2$ mapped to $b'_2\ot b'_1$ under the isomorphism $B_1\ot B_2
\simeq B_2\ot B_1$. The energy function $H$ on $B_1\ot B_2$ is defined 
up to constant difference by the following rule:
\begin{eqnarray*}
H(\tilde{e}_i(b_1\ot b_2))&=H(b_1\ot b_2)+1
&\mbox{ if }i=0,\vphi_0(b_1)\geq\veps_0(b_2),\vphi_0(b'_2)\geq\veps_0(b'_1),
\\
&=H(b_1\ot b_2)-1
&\mbox{ if }i=0,\vphi_0(b_1)<\veps_0(b_2),\vphi_0(b'_2)<\veps_0(b'_1),
\\
&\hspace{-6mm}=H(b_1\ot b_2)&\mbox{ otherwise}.
\end{eqnarray*}
Next consider the tensor product of finite crystals $B_1\ot\cd\ot B_L$.
We define $b^{(i)}_j$ ($i<j$) by
\begin{eqnarray*}
&&
\begin{array}{ccccc}\hspace{-5mm}
B_i\ot\cd\ot B_{j-1}\ot B_j&\simeq&
B_i\ot\cd\ot B_j\ot B_{j-1}&\simeq&\cd\\
b_i\ot\cd\ot b_{j-1}\ot b_j&\mapsto&
b_i\ot\cd\ot b^{(j-1)}_j\ot b'_{j-1}&\mapsto&\cd
\end{array}\\
&&\hspace{5cm}
\begin{array}{ccc}
\cd&\simeq&B_j\ot B_i\ot\cd\ot B_{j-1}\\
\cd&\mapsto&b^{(i)}_j\ot b'_i\ot\cd\ot b'_{j-1},
\end{array}
\end{eqnarray*}
and set $b^{(i)}_i=b_i$.

Take some $b_0\in B_0$. The classically restricted 1dsum is defined similarly.
\[
X(B^{r_1,s_1}\ot\cd\ot B^{r_L,s_L},\la,q)
=\mathop{{\sum}^*}q^{\sum_{0\le i<j\le L}H(b_i\ot b^{(i+1)}_j)}.
\]
Note that $b_0\in B_0$, $b_j\in B^{r_j,s_j}$ ($1\le j\le L$) and $X$ depends
also on $b_0$.
Of course, we have a similar formula to (\ref{eq:count-mult}).

Let  $b^\natural_0 = b_0 \in B_0$ and 
$b^\natural_j \in B^{r_j, s_j}$ be the highest weight element 
for $1 \le j \le L$.
\begin{conjecture} \label{conj:X=M}
Set $k=\lev B_0$. There exists an element $b_0\in B_0$ such that
$\vphi(b_0)=k\La_0$, and we have 
\[
X(B^{r_1,s_1}\ot\cd\ot B^{r_L,s_L},\la,q)
=q^c M(W^{(r_1)}_{s_1}\ot\cd\ot W^{(r_L)}_{s_L},\la,q)
\]
with $c = \sum_{0 \le i < j \le L}
H(b^\natural_i\otimes (b^\natural_j)^{(i+1)})$. Here $M$ is defined in 
(\ref{eq:mvq1}) -- (\ref{eq:mvq4}).
\end{conjecture}

The authors are informed \cite{KSS} that this conjecture is proved for 
$X_n=A_n$. In this case, if all $r_j=1$, the conjecture has already been
proved combining the results of \cite{KR1} and \cite{NY3}.

%This has been proposed for $C^{(1)}_2$, $\forall s_i = 1$ case 
%in \cite{Y}.

Let us give the ``level restricted'' version of {\sc Conjecture} \ref{conj:X=M}.
For $\la\in\ol{P}^+$ and a positive integer $l$, define
\[
X_l(B^{r_1,s_1}\ot\cd\ot B^{r_L,s_L},\la,q)
=\mathop{{\sum}^{\dagger}}q^{\sum_{0\le i<j\le L}H(b_i\ot b^{(i+1)}_j)},
\]
where the sum $\sum^{\dagger}$ is taken over all
$b_1\ot\cd\ot b_L\in B^{r_1,s_1}\ot\cd\ot B^{r_L,s_L}$ satisfying
$\sum_{j=1}^L\wt b_j=\la$ and 
\[
{\tilde e}_i^{\langle h_i,l\La_0\rangle+1}(b_1\ot\cd\ot b_L)=0
\quad\forall i.
\]
Note that the above condition is the same as (\ref{eq:hwc}) except $i=0$.
Thus it is clear that $X_l=X$ for sufficiently large $l$.

\begin{conjecture}
With the same assumptions and $c$ as in {\sc Conjecture} \ref{conj:X=M}, we have
\[
X_l(B^{r_1,s_1}\ot\cd\ot B^{r_L,s_L},0,q)
=q^c M_l(W^{(r_1)}_{s_1}\ot\cd\ot W^{(r_L)}_{s_L},q),
\] 
where $M_l$ is defined in 
(\ref{eq:mlq1}) -- (\ref{eq:mlq5}).
\end{conjecture}

For $X_n=A_n$ there are several proven cases. See e.g. 
\cite{Ber,FOW,Wa}.

\begin{remark}
In the homogeneous case, there is a similar result to (\ref{eq:bran-func}).
If $B$ is perfect of level $k \,(\le l)$, it is known that for $\la\in\ol{P}^+$,
\begin{equation}\label{eq:coset-bran}
\lim_{\begin{subarray}{c}
L\to\infty\\ L\equiv0\,\scriptstyle{mod}\,\gamma
\end{subarray}}q^{-c_L}
X_l(B^{\ot L},\la,q)
=\sum_i[[\mathcal{L}((l-k)\La_0)\ot 
\mathcal{L}(k\La_0):\la+l\La_0-i\delta]]\,q^i,
\end{equation}
where 
\begin{equation}\label{eq:aff-mult}
[[M:\mu]]=\dim\C\langle v\in M\mid\wt v=\mu,e_iv=0\,\forall i\rangle.
\end{equation}
Note the difference between (\ref{eq:cl-mult}) and (\ref{eq:aff-mult}).
The RHS of (\ref{eq:coset-bran}) is called a coset branching 
function.
If $B$ is non-perfect and level $k$, we conjecture 
(cf. \cite{Ku}, p230)
\begin{equation*}
\lim_{\begin{subarray}{c}
L\to\infty\\ L\equiv0\,\scriptstyle{mod}\,\gamma
\end{subarray}}q^{-c_L}
X_l(B^{\ot L},\la,q)
=\sum_i[[\mathcal{L}((l-k)\La_0)\ot \mathcal{L}((k-1)\La_0)
\ot \mathcal{L}^Y(\La_0):\la+l\La_0-i\delta]]\,q^i.
\end{equation*}
\end{remark}

\section{Fermionic forms}

For $m \in {\mathbb Z}_{\ge 0}$ and $p \in {\mathbb Z}$,
we define 
\begin{align}
\left\{\begin{array}{c} p + m \\ m \end{array} \right\}_q & = 
\frac{(q^{p+1})_\infty(q^{m+1})_\infty}
{(q)_\infty(q^{p+m+1})_\infty}, \label{eq:qbinomial1} \\
\left[\begin{array}{c} p + m \\ m \end{array} \right]_q & = 
\begin{cases} 
\left\{\begin{array}{c} p + m \\ m \end{array} \right\}_q 
&\mbox{ for } p \ge 0, \\
0 &\mbox{ for } p < 0, \end{cases}\label{eq:qbinomial2}
\end{align}
where $(x)_\infty = \prod_{j=0}^\infty(1-xq^j)$.
(\ref{eq:qbinomial2}) is the usual 
$q$-binomial coefficient. (\ref{eq:qbinomial1}) is an 
extended version, which vanishes only for $-m \le p \le -1$ 
and has non-zero value 
$(-q^{p+(m+1)/2})^m\left[\begin{array}{c} -p-1 \\ 
m \end{array} \right]_q$
for $p \le -m-1$.
In the $q \rightarrow 1$ limit they become 
\begin{equation*}
\left\{\begin{array}{c} p + m \\ m \end{array} \right\}_1 = 
\frac{\Gamma(p+m+1)}{\Gamma(p+1) \Gamma(m+1)}, \quad
\left[\begin{array}{c} p + m \\ m \end{array} \right]_1 = 
\left(\begin{array}{c} p+m \\ m \end{array} \right).
\end{equation*}
We shall also use the notation 
$(q)_k = (q)_\infty/(q^{k+1})_\infty$ for $k \in {\mathbb Z}_{\ge0}$.

Given any $\{ \nu^{(a)}_j \in {\mathbb Z}_{\ge 0} \mid 
j \ge 1, 1 \le a \le n \}$, we set 
$W = \bigotimes_{a=1}^n \bigotimes_{j \ge 1}
\bigl(W^{(a)}_j\bigr)^{\otimes \nu^{(a)}_j}$.
For each $1 \le a \le n$ 
we shall always assume that $\nu^{(a)}_j$'s are non-zero 
only for finitely many $j$'s.
Let $\la \in \ol{P}$. 
We define the fermionic form $M(W,\la,q)$ by
\begin{eqnarray}
M(W,\lambda,q) & = & \sum_{\{m \}} q^{c(\{m\})}
\prod_{\begin{subarray}{c} 1 \le a \le n  \\
i \ge 1 \end{subarray}}
\left[ \begin{array}{c} p^{(a)}_i +  m^{(a)}_i
 \\   m^{(a)}_i \end{array} \right]_q,\label{eq:mvq1}\\
c(\{m\}) & = & \frac{1}{2}\sum_{1 \le a, b \le n} (\alpha_a \vert \alpha_b)
\sum_{j, k \ge 1} \mbox{min}(t_bj, t_ak) m^{(a)}_j m^{(b)}_k \label{eq:mvq2}\\
&& \qquad\quad -\sum_{a=1}^n\sum_{j,k \ge 1} \nu^{(a)}_j\mbox{min}(j, k)
m^{(a)}_k,\nonumber\\
p^{(a)}_i& = & \sum_{j \ge 1}\nu^{(a)}_j\mbox{min}(i,j)
- \sum _{b=1}^{n} (\alpha_a \vert \alpha_b) \sum_{k \ge 1} 
\mbox{min}(t_bi,t_ak) m^{(b)}_k,
\label{eq:mvq3}
\end{eqnarray}
where the sum $\sum_{\{ m \}}$ is taken over
$\{ m^{(a)}_i \in {\mathbb Z}_{\ge 0} \mid 1 \le a \le n, \, i \ge 1 \}$
satisfying $p^{(a)}_i \ge 0$ for $1 \le a \le n, i \ge 1$, and
\begin{equation}\label{eq:mvq4}
\sum_{a=1}^n \sum_{i \ge 1} i m^{(a)}_i \alpha_a = 
\sum_{a=1}^n\sum_{i \ge 1} i\nu^{(a)}_i \ol{\La}_a - \lambda
\qquad \mbox{ for } \ 1 \le a \le n.
\end{equation}
By definition $M(W,\la,q) = 0$ if the RHS of (\ref{eq:mvq4})
does not belong to 
$\sum_{a=1}^n {\mathbb Z}_{\ge 0} \alpha_a$.
More strictly we have 
\begin{equation}\label{eq:mzero}
M(W,\la,q) = 0 \quad \mbox{ unless } \la \in 
\Bigl(\sum_{a=1}^n \sum_{i \ge 1} i \nu^{(a)}_i \ol{\La}_a
 - \sum_{a=1}^n {\mathbb Z}_{\ge 0} \alpha_a\Bigr) \cap 
\ol{P}^+.
\end{equation}
To see this, compare (\ref{eq:mvq3}) and (\ref{eq:mvq4}), 
which leads to
$p^{(a)}_\infty = (t_a \alpha_a \vert \la)$.
Thus the above property (\ref{eq:mzero}) holds 
because of (\ref{eq:qbinomial2}).
It is easy to see 
\begin{equation*}
M(W,\la,q) \in \mathbb{Z}_{\ge 0}[q^{-1}].
\end{equation*}
Note that
$
p^{(a)}_i = -\frac{\partial c(\{m\})}{\partial m^{(a)}_i}
$
{}from (\ref{eq:mvq2}) and (\ref{eq:mvq3}), regarding 
all $m^{(a)}_i$ as independent continuous variables.
Similar relations are valid also between $c_l(\{m\})$ and
$p^{(a)}_i$'s in the later 
equations (\ref{eq:mlq2})-(\ref{eq:mlq3}) and 
(\ref{eq:mlq22})-(\ref{eq:mlq23}).

To introduce a ``restricted version'' of $M(W,\la,q)$, 
we fix a positive integer $l$.
For $W = \bigotimes_{(a,j) \in H_l} 
\bigl(W^{(a)}_j\bigr)^{\otimes \nu^{(a)}_j}$ 
we define
\begin{eqnarray}
M_l(W,q) & = & \sum_{\{m \}} q^{c_l(\{m\})}
\prod_{(a,i) \in H_l}
\left[ \begin{array}{c} p^{(a)}_i +  m^{(a)}_i
 \\   m^{(a)}_i \end{array} \right]_q,\label{eq:mlq1}\\
c_l(\{m\}) & = & \frac{1}{2}\sum_{(a,j), (b,k) \in H_l} 
(\alpha_a \vert \alpha_b)
 \mbox{min}(t_bj, t_ak) m^{(a)}_j m^{(b)}_k \label{eq:mlq2}\\
&& \qquad\quad -\sum_{(a,j) \in H_l}
\sum_{k=1}^{t_al} \nu^{(a)}_j\mbox{min}(j, k)
m^{(a)}_k,\nonumber\\
p^{(a)}_i& = & \sum_{j=1}^{t_al}\nu^{(a)}_j\mbox{min}(j,i)
- \sum _{(b,k) \in H_l} (\alpha_a \vert \alpha_b) 
\mbox{min}(t_bi,t_ak) m^{(b)}_k,
\label{eq:mlq3}
\end{eqnarray}
where the sum $\sum_{\{ m \}}$ is taken over
$\{ m^{(a)}_i \in {\mathbb Z}_{\ge 0} \mid (a,i) \in H_l \}$
satisfying $p^{(a)}_i \ge 0$ for $(a,i) \in H_l$, and
\begin{equation}\label{eq:mlq4}
\sum_{(a,i) \in H_l} i\, m^{(a)}_i \alpha_a = 
\sum_{(a,i) \in H_l} i\nu^{(a)}_i \ol{\La}_a
\qquad \mbox{ for } 1 \le a \le n,
\end{equation}
or equivalently,
\begin{equation}\label{eq:mlq5}
\sum_{k=1}^{t_al}k m^{(a)}_k = 
\sum_{(b,i) \in H_l} C^{-1}_{a b}i \nu^{(b)}_i
\qquad \mbox{ for } 1 \le a \le n.
\end{equation} 

Notice that (\ref{eq:mlq4}) is equivalent to $p^{(a)}_{t_al} = 0$.
Eliminating $m^{(a)}_{t_al}$ ($1 \le a \le n$) {}from 
(\ref{eq:mlq5}) one can rewrite (\ref{eq:mlq1}) -- (\ref{eq:mlq4}) 
as follows.
\begin{eqnarray}
M_l(W,q) & = & \sum_{\{m \}} q^{c_l(\{m\})}
\prod_{(a,i) \in {\overline H}_l}
\left[ \begin{array}{c} p^{(a)}_i +  m^{(a)}_i
 \\   m^{(a)}_i \end{array} \right]_q,\label{eq:mlq21}\\
c_l(\{m\}) & = & \frac{1}{2}\sum_{(a,j), (b,k) \in {\overline H}_l} 
(\alpha_a \vert \alpha_b)
 K^{(t_at_bl)}_{t_bj, t_ak} m^{(a)}_j m^{(b)}_k \label{eq:mlq22}\\
&&  -\sum_{(a,j) \in {\overline H}_l}
\sum_{i=1}^{t_al-1}K^{(t_al)}_{i, j}
 \nu^{(a)}_i m^{(a)}_j
- \frac{\vert \sum_{(a,j) \in H_l} j \nu^{(a)}_j \ol{\La}_a \vert ^2}{2l},
\nonumber\\
p^{(a)}_i& = & \sum_{j=1}^{t_al-1}K^{(t_al)}_{i, j}\nu^{(a)}_j
- \sum _{(b,k) \in {\overline H}_l} (\alpha_a \vert \alpha_b) 
K^{(t_at_bl)}_{t_bi,t_ak} m^{(b)}_k.
\label{eq:mlq23}
\end{eqnarray}
Here the sum $\sum_{\{ m \}}$ is taken over
$\{ m^{(a)}_i \in {\mathbb Z}_{\ge 0} \mid (a,i) \in {\overline H}_l \}$
satisfying $p^{(a)}_i \ge 0$ for $(a,i) \in {\overline H}_l$, and
$m^{(a)}_{t_al}$ defined by (\ref{eq:mlq5}) is a non-negative
integer for $1 \le a \le n$.

\begin{remark}
The two fermionic forms $M$ and $M_l$ are related by
$$M(W,0,q) = M_\infty(W,q).$$
In the RHS, the limit $l \rightarrow \infty$ poses 
no subtlety.
\end{remark}

\begin{remark}
For $W = \bigotimes_{(a,i) \in H_l} \left( W^{(a)}_j \right)^{
\otimes \nu^{(a)}_j}$, set
$\ol{M}_l(W,q) = q^{\vert \La \vert^2/2l} M_l(W,q)$, 
$\La = \sum_{(a,j) \in H_l} j \nu^{(a)}_j \ol{\La}_a$,
which corresponds to dropping the last term in (\ref{eq:mlq22}).
For simply laced $X_n$, $\ol{H}_{l=1}$ (\ref{eq:GKdef2}) is empty 
and therefore
{}from (\ref{eq:mlq21}), 
$$\ol{M}_1(W,q) = 1 \quad \mbox{ for } \ A_n, D_n, E_{6,7,8}.
$$
On the other hand for non-simply laced $X_n$, 
$\ol{H}_1$ is bijective to $\ol{H}_t$ of $Z_n$ given in {\sc Table} 
\ref{tab:YZ}. ($t > 1$ here is the one for $X_n$ and not for $Z_n$.)
Moreover a simple manipulation tells that
$$\ol{M}_1(W,q) \ \mbox{ for } \ X_n = 
\ol{M}_t(W',q) \ \mbox{ for } \ Z_n,
$$
where $W' = \bigotimes_{(a,i) \in H_t} \left( W^{(a)'}_j \right)^{
\otimes \nu^{(a)'}_j}$ is specified by
($H_t$ and $W^{(a)'}_j$ here denote those for $Z_n$.)
\begin{eqnarray*}
B_n &:& \nu^{(1)'}_1 = \nu^{(n)}_1, \ \
\nu^{(1)'}_2 = \sum_{a=1}^{n-1} a \nu^{(a)}_1 + 
\frac{n-1}{2} \nu^{(n)}_1 + n \nu^{(n)}_2, \\
C_n &:& \nu^{(a)'}_1 = \nu^{(a)}_1, \ \
\nu^{(a)'}_2 = \nu^{(a)}_2 + \frac{\delta_{a, n-1}}{2}
\left( \sum_{a=1}^{n-1}a(\nu^{(a)}_1 + 2 \nu^{(a)}_2) + n\nu^{(n)}_1
\right) \ 1 \le a \le n-1,\\
F_4 &:& \nu^{(a)'}_j = \nu^{(a+2)}_j\ 1 \le a, j \le 2 \ \mbox{except }
\ \nu^{(1)'}_2 = 3\nu^{(1)}_1 + 6\nu^{(2)}_1 + 
4\nu^{(3)}_1 + 
9\nu^{(3)}_2 + 
2\nu^{(4)}_1 + 
4\nu^{(4)}_2,\\
G_2 &:& \nu^{(1)'}_j = \nu^{(2)}_j \ 1 \le j \le 2, \ \
\nu^{(1)'}_3 = 2\nu^{(1)}_1 + 
\nu^{(2)}_1 + 
2\nu^{(2)}_2.
\end{eqnarray*}
In the above the condition (\ref{eq:mlq5}) $\in {\mathbb Z}$ ensures 
$\forall \nu^{(a)'}_j \in {\mathbb Z}$. 
\end{remark}

The above fermionic forms are defined in terms of the usual
$q$-binomial 
(\ref{eq:qbinomial2}).  However we shall also concern those
involving (\ref{eq:qbinomial1}).
Given $l \in {\mathbb Z}_{\ge 1}$, 
$\la = \la_1 \ol{\La}_1 + \cdots + 
\la_n \ol{\La}_n \in \ol{P}$ and 
$W = \bigotimes_{(a,j) \in H_l} 
\bigl(W^{(a)}_j\bigr)^{\otimes \nu^{(a)}_j}$, 
we introduce the third fermionic form
\begin{equation}
N_l(W,\la,q) = \sum_{\{m \}} q^{c_l(\{m\})}
\prod_{(a,i) \in H_l}
\left\{ \begin{array}{c} p^{(a)}_i +  m^{(a)}_i
 \\   m^{(a)}_i \end{array} \right\}_q,\label{eq:nlq1}
\end{equation}
where $c_l(\{m\})$ and $p^{(a)}_i$ are defined by the same formulae
as (\ref{eq:mlq2}) and (\ref{eq:mlq3}), respectively.
The sum $\sum_{\{ m \}}$ is taken over
$\{ m^{(a)}_i \in {\mathbb Z}_{\ge 0} \mid (a,i) \in H_l \}$
{\em not} bounded by the restriction 
$p^{(a)}_i \ge 0$ for $(a,i) \in H_l$, but only subject to
the condition
\begin{equation}\label{eq:nlq2}
\sum_{(a,i) \in H_l} i\, m^{(a)}_i \alpha_a = 
\sum_{(a,i) \in H_l} i\nu^{(a)}_i \ol{\La}_a -\la
\qquad \mbox{ for } 1 \le a \le n.
\end{equation}
By introducing 
\begin{equation}
\gamma^{(a)}_j = \sum_{k=1}^{t_al} \nu^{(a)}_k \mbox{min }
(j,k), \quad 
\mu_a = \gamma^{(a)}_\infty - \la_a,
\label{eq:nlq3}
\end{equation}
one can rewrite (\ref{eq:nlq2}) and $p^{(a)}_i$ (\ref{eq:mlq3}) as
\begin{align}
\mu_a &= 
\sum_{(b,k)\in H_l} C_{a b}k m^{(b)}_k,
\label{eq:nlq4}\\
p^{(a)}_i &= \gamma^{(a)}_i - \mu_a + \sum_{b=1}^n
C_{a b} \sum_{\frac{t_b}{t_a}i < k \le t_bl}
\left(k-\frac{t_b}{t_a}i\right)m^{(b)}_k.
\label{eq:nlq5}
\end{align}
Although $N_l(W,\la,q)$ is defined 
for all $\la \in \ol{P}$, it is vanishing unless 
$\nu^{(a)}_i$'s and $\la$ are chosen so that 
$\sum_{b=1}^n C^{-1}_{a b}\mu_b \in {\mathbb Z}_{\ge 0}$
for all $1 \le a \le n$.

In contrast with $M$ and $M_l$, the above fermionic form 
contains, in general, non-zero 
contributions with possibly negative signs 
{}from $p^{(a)}_i < -m^{(a)}_i$ cases.
See the remarks after (\ref{eq:qbinomial2}).
Nevertheless our computer experiments suggest
\begin{conjecture}\label{con:MN}
\begin{align}
M(W,\la,q) &= N_\infty(W,\la,q) \quad
\mbox{ for } \la \in \ol{P}^+, \label{con:MN1} \\
M_l(W,q) &= N_l(W,0,q).\label{con:MN2}
\end{align}
\end{conjecture}
The conjecture means that those ``unphysical" contributions in the 
Bethe ansatz context cancel out totally.
So far we have checked it for most of the 
fundamental representations $W = W^{(a)}_1$. 
(The case $X_n = A_1, \, \nu^{(1)}_j=L\delta_{j, 1}$
is exceptional in that the 
``unphysical contributions" are absent also in $N_l(W,\la,q)$  
because  $p^{(1)}_1 \ge p^{(1)}_2 \ge \cdots \ge p^{(1)}_l 
= (\alpha_1 \vert \la) \ge 0$ for $\la \in \ol{P}^+$.)
On the contrary, for $\la  \not \in  \ol{P}^+$ 
the two fermionic forms 
$M$ and $N_\infty$ are significantly different.
Compare (\ref{eq:mzero}) and {\sc Conjecture} \ref{con:nmweyl}.
One can define $M_l(W,\la,q)$ by restricting the summands in  
$N_l(W,\la,q)$ to those satisfying 
$\forall p^{(a)}_i \ge 0$. However 
$M_l(W,\la,q) \neq N_l(W,\la,q)$ in general for 
$\la \neq 0$ as opposed to (\ref{con:MN2}).

%%%%%%%%%%%%%%%%%%%%%%%%%%%%%%%%%%%%%%%%%%%%
\section{Generalized spinon character formula}
In this section we concentrate on the 
fermionic form $M(W^{(r)\otimes L}_s,\lambda,q)$
for the homogeneous quantum space $W = W^{(r)\otimes L}_s$
and its limit $L \rightarrow \infty$ such that 
$Ls\ol{\La}_r \in \ol{Q}$.
Due to (\ref{eq:bran-func})--(\ref{eq:bran-func2})
 and {\sc Conjecture} \ref{con:XM}, 
it yields explicit formulae for the 
branching functions. 
The result turns out to be a generalization of the spinon character formulae 
\cite{NY1,NY2,ANOT,BPS,BS,BLS}.
Such a connection between the spinon character formula and the 
limit of the fermionic form was first shown in \cite{HKKOTY} for 
$A^{(1)}_n$.

First we seek the minimum point of the quadratic form 
$c(\{m\})$ (\ref{eq:mvq2}) for $\nu^{(a)}_j = L\delta_{a, r}\delta_{j, s}$.
The solution of the linear equation 
$\frac{\partial c(\{m\})}{\partial m^{(a)}_j} = 0$ is 
given by
\begin{lemma}\label{lem:critical} The simultaneous equation
\begin{equation}\label{eq:critical}
\sum_{b=1}^n \sum_{k \ge 1} (\alpha_a \vert \alpha_b) 
\mbox{min}(t_bj, t_ak) m^{(b)}_{k,0} = L \delta_{a, r}
\mbox{min}(s,j)
\end{equation}
for $1 \le a \le n, j \ge 1$  
is satisfied by
\begin{equation}
m^{(a)}_{j,0}  =  LC^{-1}_{r a} \delta_{j, \frac{t_a}{t_r}s}
\nonumber
\end{equation}
for any $X_n$  provided that $\frac{s}{t_r} \in {\mathbb Z}$.
In case $\frac{s}{t_r} \not\in {\mathbb Z}$, (\ref{eq:critical}) holds for 
\begin{equation}
B_n:\quad
m^{(a)}_{j,0}  =  \begin{cases}
\frac{La}{2}(\delta_{j,\frac{s-1}{2}} + \delta_{j,\frac{s+1}{2}})
& \quad 1 \le a \le n-1,\\
\frac{L}{4}\bigl((n-1)(\delta_{j,s-1}+\delta_{j,s+1})+
2\delta_{j,s}\bigr) & \quad a = n,
\end{cases}\nonumber
\end{equation}
\begin{equation}
C_n:\quad
m^{(a)}_{j,0} = \frac{Lar}{2n}(\delta_{j,\frac{t_a}{2}(s-1)} +
\delta_{j,\frac{t_a}{2}(s+1)}) + L(\mbox{min}(a,r) - \frac{ar}{n})
\delta_{j, \frac{t_a}{2}s},\nonumber
\end{equation}
\begin{eqnarray}
F_4:\quad
m^{(a)}_{j,0} & = & \frac{L}{3}(2\delta_{r,3}+ \delta_{r,4})
\bigl(C^{-1}_{a 2}(\delta_{j,\frac{t_a}{2}(s-1)} +
\delta_{j,\frac{t_a}{2}(s+1)}) +
\delta_{a, 3}\delta_{j,s}\bigr) \nonumber\\
& & + \frac{L}{3}
(\delta_{r,3}+ 2\delta_{r,4})\delta_{a,4}\delta_{j,s},\nonumber
\end{eqnarray}
where the inverse Cartan matrix of $F_4$ is given 
by (\ref{eq:cartan1}) and (\ref{eq:inversecartan}).
\begin{eqnarray}
G_2: \quad m^{(1)}_{j,0} & = & \begin{cases}
L(2\delta_{j, \frac{s-1}{3}} + \delta_{j,\frac{s+2}{3}}) 
& s \in 3{\mathbb Z}+1,\\
L(\delta_{j, \frac{s-2}{3}} + 2\delta_{j,\frac{s+1}{3}}) 
& s \in 3{\mathbb Z}+2,\end{cases} \nonumber \\
m^{(2)}_{j,0}& = & \begin{cases}
L(\delta_{j, s-1} + \frac{1}{2}\delta_{j,s} 
+ \frac{1}{2}\delta_{j,s+2}) & s \in 3{\mathbb Z}+1,\\
L(\frac{1}{2}\delta_{j, s-2} + \frac{1}{2}\delta_{j,s} 
+ \delta_{j,s+1}) & s \in 3{\mathbb Z}+2.\end{cases}\nonumber
\end{eqnarray}
\end{lemma}
\begin{remark}\label{rem:critical1}
In all the cases in {\sc Lemma} \ref{lem:critical}, one has
$$
\sum_{a=1}^n \sum_{j \ge 1} j m^{(a)}_{j,0} \alpha_a = Ls\ol{\La}_r.
$$
When $s/t_r \in {\mathbb Z}$, 
$m^{(a)}_{i,0}$ is an integer iff $LsC^{-1}_{r a} \in {\mathbb Z}$.
When $s/t_r \not \in {\mathbb Z}$, it is so iff
\begin{eqnarray*}
B_n: && L \in 2{\mathbb Z} \,\, (n:\mbox{ odd}), 
\quad L \in 4{\mathbb Z} \,\, (n:\mbox{ even}),\\
C_n: && L \in 2n{\mathbb Z},\\
F_4: && L \in 3{\mathbb Z},\\
G_2: && L \in 2{\mathbb Z}.
\end{eqnarray*}
For such $L$ we have $Ls\ol{\La}_r \in \ol{Q}$.
\end{remark}
\begin{remark}\label{rem:critical2}
$ m^{(a)}_{j, 0} \neq 0$ iff
\begin{equation*}
(a,j)  \in \begin{cases}
\{ (a,\frac{t_a}{t_r}) \mid 1 \le a \le n \}\quad
\mbox{ for } \frac{s}{t_r} \in {\mathbb Z},\\[5pt]
\{ (a,\frac{t_a(s-s_0)}{t_r}),
 (a,\frac{t_a(s-s_0+t_r)}{t_r}) \mid 1 \le a \le n \}\\
\cup \{(a,s) \mid 1 \le a \le n, t_a = t_r(>1)\} \quad 
\mbox{ for } \frac{s}{t_r} \not\in {\mathbb Z},
\end{cases}
\end{equation*}
where $s_0 \equiv s$ mod $t_r{\mathbb Z}$ and $1 \le s_0 \le t_r-1$.
This is essentially the same with (3.5b) in \cite{Ku}, 
which appeared in the 
Bethe ansatz analysis.
It corresponds to the physical picture 
that the ground state in ``regime III-like" region 
is a Dirac sea of color $a$ length $j$-strings for those
$(a,j)$ such that $m^{(a)}_{j,0} \rightarrow \infty$ as 
$L \rightarrow \infty$. The value $m^{(a)}_{i,0}$ should coincide
with the $0$-th Fourier component of the ground state 
density function of the $(a,j)$-string.
\end{remark}

As in {\sc Lemma} \ref{lem:critical} the limit of 
$M(W^{(r)\otimes L}_s,\lambda,q)$ is considerably 
complicated if $s/t_r \not\in {\mathbb Z}$.
We therefore treat the case $s/t_r \in {\mathbb Z}$ first.
\begin{theorem} \label{th:spinon1}
Assume that $s/t_r \in {\mathbb Z}$. Then
\begin{eqnarray}
\lim
q^{-c(\{m^{(a)}_{i,0}\})} M(W^{(r)\otimes L}_s,\lambda,q)
& = & \sum_{\zeta}\frac{M(W(\zeta),\lambda,q^{-1}) 
M_{\frac{s}{t_r}}(W(\zeta),q^{-1})}
{(q)_{\zeta_1} \cdots (q)_{\zeta_n}},\\
W(\zeta) & = & \bigotimes_{a=1}^n W^{(a) \otimes \zeta_a}_1.
\end{eqnarray}
Here the limit 
$L \rightarrow \infty$ is taken so that 
$LC^{-1}_{r a} \in {\mathbb Z}$ for $1 \le a \le n$.
$c(\{m^{(a)}_{i,0}\}) = -L^2sC^{-1}_{r r}/2 \in {\mathbb Z}$.
The sum is taken over $\zeta_1, \ldots, \zeta_n \in {\mathbb Z}_{\ge 0}$
such that $\zeta = \sum_{a=1}^n \zeta_a \ol{\La}_a \in \ol{Q}$.
\end{theorem}
\begin{proof}
We start with the expression (\ref{eq:mvq1}) -- (\ref{eq:mvq4}).
The limit is to be expanded {}from the minimum 
$m^{(a)}_i = m^{(a)}_{i,0}$ of $c(\{m\})$ 
determined in {\sc Lemma} \ref{lem:critical}. 
The non-zero $m^{(a)}_{i,0}$'s are proportional to $L$.
Thus under the identification $\zeta_a = p^{(a)}_{\frac{st_a}{t_r}}$,
the factor
$\prod_{a=1}^{n} \left[ \begin{array}{c} 
p^{(a)}_{st_a/t_r} +  m^{(a)}_{st_a/t_r}
 \\   m^{(a)}_{st_a/t_r} \end{array} \right]_q$
in (\ref{eq:mvq1}) gives rise to
$\left((q)_{\zeta_1} \cdots (q)_{\zeta_{n}}\right)^{-1}$
as $L \rightarrow \infty$.
After the shift
$m^{(a)}_i \rightarrow m^{(a)}_i + m^{(a)}_{i, 0}$, the relation
(\ref{eq:mvq3}) and its $i=\frac{st_a}{t_r}$ case become
\begin{eqnarray}
p^{(a)}_i & = & - \sum_{b=1}^{n} (\alpha_a \vert \alpha_b)
\sum_{k \ge 1} \mbox{min}(t_bi,t_ak)m^{(b)}_k, \label{pr:spinon12}\\
\sum_{b=1}^{n}C^{-1}_{a b} \zeta_b & = &
- \sum_{k \ge 1} \mbox{min}(\frac{t_as}{t_r},k) m^{(a)}_k.
\label{pr:spinon13}
\end{eqnarray}
Eliminating $m^{(a)}_{\frac{t_as}{t_r}}$ with (\ref{pr:spinon13}),
one rewrites (\ref{pr:spinon12}) as
$$
p^{(a)}_i =
\left\{
\begin{array}{ll}
\zeta_a - \sum_{b=1}^n(\alpha_a \vert \alpha_b)
 \sum_{k > \frac{t_bs}{t_r}}
\mbox{min}(t_b(i-\frac{t_as}{t_r}),t_a(k-\frac{t_bs}{t_r}))
m^{(b)}_k & \quad i > \frac{t_as}{t_r}\\
K^{(\frac{t_as}{t_r})}_{1, \frac{t_as}{t_r}-i}\zeta_a -
\sum_{b=1}^{n}(\alpha_a \vert \alpha_b)\sum_{k=1}^{\frac{t_bs}{t_r}-1}
K^{(\frac{t_at_bs}{t_r})}_{t_bi,t_ak}m^{(b)}_k& 
\quad 1 \le i < \frac{t_as}{t_r}\end{array} \right. .
$$
By setting
$p^{(a)}_i = \hat{p}^{(a)}_{i-\frac{t_as}{t_r}}$ and
$m^{(a)}_i = \hat{m}^{(a)}_{i-\frac{t_as}{t_r}}$, 
$i > \frac{t_as}{t_r}$ case
in the above becomes the relation (\ref{eq:mvq3}) 
between $\hat{p}^{(a)}_i$ and $\hat{m}^{(a)}_i$ 
for $M(W(\zeta),\lambda, q^{-1})$.
Similarly, by setting 
$p^{(a)}_i = \tilde{p}^{(a)}_{\frac{t_as}{t_r}-i}$ and
$m^{(a)}_i = \tilde{m}^{(a)}_{\frac{t_as}{t_r}-i}$, 
$i < \frac{t_as}{t_r}$ case
becomes the relation (\ref{eq:mlq23}) 
between $\tilde{p}^{(a)}_i$ and $\tilde{m}^{(a)}_i$ 
for $M_{\frac{s}{t_r}}(W(\zeta), q^{-1})$.
Thus in the product 
$\prod_{1 \le a \le n,i \ge 1}$ of (\ref{eq:mvq1}), 
we are to extract $M(W(\zeta),\lambda,q^{-1})$ 
{}from $\prod_{1 \le a \le n} \prod_{i > \frac{st_a}{t_r}}$ part and 
$M_{\frac{s}{t_r}}(W(\zeta),q^{-1})$ {}from 
$\prod_{(a,i) \in \overline{H}_{\frac{s}{t_r}}}$ part.
See (\ref{eq:GKdef2}) for the definition of 
$\overline{H}_{\frac{s}{t_r}}$.
Actually, the decomposition
\begin{eqnarray*}
&&-c(\{m^{(a)}_{i,0}\}) 
+ c(\{m^{(a)}_i \rightarrow m^{(a)}_i+ m^{(a)}_{i,0}\}) \\
& = & -\sum_{1 \le a \le n, i \ge 1}\hat{m}^{(a)}_i\hat{p}^{(a)}_i
-c(\{\hat{m}\}) 
-\sum_{(a,i) \in \overline{H}_l}\tilde{m}^{(a)}_i\tilde{p}^{(a)}_i
-c_{\frac{s}{t_r}}(\{\tilde{m}\}) 
\end{eqnarray*}
is valid among the quadratic forms (\ref{eq:mvq2}) and 
(\ref{eq:mlq22}). 
The first (resp. last) two terms on the RHS 
yield the quadratic form in the variables
$\{ \hat{m}^{(a)}_i \}$ (resp. $\{ \tilde{m}^{(a)}_i \}$)
for $M(W(\zeta),\lambda,q^{-1})$ 
(resp. $M_{\frac{s}{t_r}}(W(\zeta),q^{-1})$).
It remains to establish that the sum over
$\{m^{(a)}_j \mid 1 \le a \le n, j \ge 1 \}$ is 
properly translated into the sums over
$\{\hat{m}^{(a)}_j \mid 1 \le a \le n, j \ge 1 \}$, 
$\{\zeta_a \mid 1 \le a \le n \}$ and 
$\{\tilde{m}^{(a)}_j \mid (a,j) \in \overline{H}_l \}$.
For this it suffices to check 
(i) $\sum_{a=1}^n \zeta_a \ol{\La}_a \in \ol{Q}$, 
(ii) $\sum_{a=1}^n \zeta_a \ol{\La}_a - \sum_{a=1}^n\sum_{i \ge 1}
i \hat{m}^{(a)}_i \alpha_a = \lambda$ (see (\ref{eq:mvq4})).
(iii) $\tilde{m}^{(a)}_{\frac{t_as}{t_r}}$ determined {}from
$\sum_{k=1}^{\frac{t_as}{t_r}}k \tilde{m}^{(a)}_k = 
\sum_{b=1}^n C^{-1}_{a b}\zeta_b$ (see (\ref{eq:mlq5}))
is a non-negative integer.
(i) Due to (\ref{pr:spinon13}), 
$\sum_{a=1}^n \zeta_a \ol{\La}_a = 
- \sum_{a=1}^n \sum_{k \ge 1} \mbox{min}(\frac{t_as}{t_r},k) m^{(a)}_k
\alpha_a \in \ol{Q}$ holds.
(ii) Making the shift $m^{(a)}_i \rightarrow m^{(a)}_i + m^{(a)}_{i,0}$
in (\ref{eq:mvq4}) and applying {\sc Remark} \ref{rem:critical1}, we have
$\lambda = -\sum_{a=1}^n \sum_{k \ge 1} k m^{(a)}_k \alpha_a$.
Combining this with the above we have 
$\sum_{a=1}^n \zeta_a \ol{\La}_a - \lambda = 
- \sum_{a=1}^n\sum_{k > \frac{t_as}{t_r}}(k-\frac{t_as}{t_r})
m^{(a)}_k \alpha_a$.
(iii) Replace the RHS of 
$\sum_{k=1}^{\frac{t_as}{t_r}}k \tilde{m}^{(a)}_k = 
\sum_{b=1}^n C^{-1}_{a b}\zeta_b$ with that of (\ref{pr:spinon13}),  
leading to
$\tilde{m}^{(a)}_{\frac{t_as}{t_r}} = 
\frac{t_r}{t_as}\left(
\sum_{b=1}^nC^{-1}_{a b} \zeta_b - \sum_{i=1}^{\frac{t_as}{t_r}-1}
(\frac{t_as}{t_r} - i) m^{(a)}_i \right)
= - \sum_{i \ge 1} m^{(a)}_i \in {\mathbb Z}$.
Thus we are left to verify $\tilde{m}^{(a)}_{\frac{t_as}{t_r}} \ge 0$.
Setting $i = t_a$ in (\ref{pr:spinon12}) we have
$- \sum_{k \ge 1}\mbox{min}(t_a, k) m^{(a)}_k = 
\sum_{b=1}^n C^{-1}_{a b} p^{(b)}_{t_b} \ge 0$.
This tells that $\tilde{m}^{(a)}_{\frac{t_as}{t_r}} \ge 0$ indeed holds
if $t_a = 1$.
Henceforth we assume that $t_a > 1$ and use (\ref{pr:spinon12}) 
again for $i=1$,
which reads 
$p^{(a)}_1 = - \sum_{b=1}^n C_{b a} \sum_{k \ge 1}m^{(b)}_k$.
{}From the above expression of $\tilde{m}^{(a)}_{\frac{t_as}{t_r}}$,
we rewrite this as 
\begin{equation}\label{pr:spinon14}
\sum_{1 \le b \le n, t_b > 1} C_{b a} \tilde{m}^{(b)}_{\frac{t_bs}{t_r}}
= p^{(a)}_1 - 
\sum_{1 \le b \le n, t_b = 1} C_{b a} \tilde{m}^{(b)}_{\frac{t_bs}{t_r}},
\end{equation}
which is valid for each $1 \le a \le n$ such that $t_a > 1$
For $t_b = 1$ 
we have already shown that $\tilde{m}^{(b)}_{\frac{t_bs}{t_r}} \ge 0$, 
and moreover 
$C_{b a} \le 0$ because $t_a > 1 = t_b$. 
Therefore the RHS of (\ref{pr:spinon14}) is non-negative.
The coefficients $(C_{b a})_{t_a, t_b > 1}$ form a 
symmetric submatrix, which is the Cartan matrix
of the algebra $Z_n$ in {\sc Table} \ref{tab:YZ}.
(Simply laced $X_n$ is out of question here.)
Thus all the elements $D^{-1}_{a b}$ (\ref{eq:dmat})
of their inverse are positive,
which completes the proof.
\end{proof}

The above proof is a straightforward generalization of the one for 
$X_n = A_n$ given in \cite{HKKOTY}.
Next we proceed to the case $s/t_r \in {\mathbb Z}+\frac{1}{2}$.
\begin{theorem}\label{th:spinon2}
Let $X_n = B_n, C_n, F_4$ and assume that 
$t_r = 2$, $s \in 2{\mathbb Z} + 1$. Then we have
\begin{eqnarray}
&&\lim q^{-c(\{m^{(a)}_{i,0}\})} M(W^{(r)\otimes L}_s,\lambda,q)\\
&& = \sum_{\xi,\eta,\zeta}\frac{q^{\phi}M(W(\eta),\lambda,q^{-1}) 
M_{\frac{s-1}{2}}(W(\xi),q^{-1})}
{\left(\prod_{1 \le a \le n}(q)_{\eta_a}(q)_{\xi_a}\right)
\prod_{1 \le a \le n,t_a = 2}(q)_{\zeta_a}},\nonumber\\
&&W(\eta) =  \bigotimes_{a=1}^n W^{(a) \otimes \eta_a}_1, \quad
W(\xi)  =  \bigotimes_{a=1}^n W^{(a) \otimes \xi_a}_1, \\
&&\phi = \frac{1}{2}\vert \sum_{a=1}^n(\eta_a - \xi_a)\ol{\La}_a \vert^2
+ \sum_{\begin{subarray}{c}
1\le a, b \le n \\ t_a = t_b = 2 \end{subarray}}
D^{-1}_{a b} (\zeta_a - \frac{\xi_a+\eta_a}{2})
(\zeta_b - \frac{\xi_b+\eta_b}{2}),
\end{eqnarray}
where $D^{-1}$ has been defined in (\ref{eq:dmat}). 
The limit 
$L \rightarrow \infty$ is taken under 
the condition in {\sc Remark} \ref{rem:critical1}.
The sum runs over non-negative integers
$\eta_a, \xi_a \ (1 \le a \le n)$ and 
$\zeta_a \ (1 \le a \le n, t_a = 2)$ satisfying the
constraints
\begin{eqnarray*}
B_n:&&\quad \xi_n, \ \ \eta_n, \ \ 
\zeta_n + \sum_{a=1}^{n-1}a(\xi_a + \eta_a) +
\frac{n-1}{2}(\xi_n + \eta_n) \in 2{\mathbb Z},\\
C_n:&& \quad \sum_{a=1}^n a \xi_a,\ 
 \sum_{a=1}^n a \eta_a \in 2{\mathbb Z},\ \ 
\sum_{a=1}^{n-1}a\zeta_a - 
\frac{1}{2}\sum_{a=1}^na(\xi_a + \eta_a) \in n{\mathbb Z},\\
F_4:&&\quad  \xi_3 + 2\xi_4 + \eta_3 + 2 \eta_4 + \zeta_3 + 2\zeta_4 
\in 3{\mathbb Z}.
\end{eqnarray*}
\end{theorem}

Finally we consider $s/t_r \in {\mathbb Z} \pm \frac{1}{3}$ case.
\begin{theorem}\label{th:spinon3}
Let $X_n = G_2$ and $r=2$ hence $t_r=3$. Assume that 
$s \in 3{\mathbb Z}+1$ or $3{\mathbb Z}+2$, for which  
$s_0= 1$ or $2$, respectively as in {\sc Remark} \ref{rem:critical2}. 
Then we have 
\begin{eqnarray}
&&\lim q^{-c(\{m^{(a)}_{i,0}\})} M(W^{(2)\otimes L}_s,\lambda,q)\\
&& = \sum_{\xi,\eta,\zeta,m}\frac{q^{\phi}M(W(\eta),\lambda,q^{-1}) 
M_{\frac{s-s_0}{3}}(W(\xi),q^{-1})\Phi}
{\left(\prod_{1 \le a \le 2}(q)_{\eta_a}(q)_{\xi_a}\right)
(q)_{\zeta}},\nonumber\\
&&W(\eta) =  \bigotimes_{a=1}^2 W^{(a) \otimes \eta_a}_1, \quad
W(\xi)  =  \bigotimes_{a=1}^2 W^{(a) \otimes \xi_a}_1, \\
&&\phi = \frac{1}{2}\vert \sum_{a=1}^2(\eta_a - \xi_a)\ol{\La}_a \vert^2
+ \frac{m^2}{2} + 
\frac{3}{8}\left(\zeta - \frac{s_0\eta_2 + (3-s_0)\xi_2}{3}\right)^2,\\
&&\Phi = \left[ \begin{array}{c} 
\frac{1}{2}\left(\zeta+(2-s_0)\eta_2+(s_0-1)\xi_2 \right)
 \\  m \end{array} \right]_q.
\end{eqnarray}
Here the limit 
$L \rightarrow \infty$ is taken so that $L \in 2{\mathbb Z}$  
as in {\sc Remark} \ref{rem:critical1}.
The sum runs over non-negative integers
$\eta_a, \xi_a \ (1 \le a \le 2), \  
\zeta$ and $m$ satisfying 
$$
\zeta + (2-s_0)\eta_2 + (s_0-1)\xi_2, \ \ 
m + \frac{(3-s_0)\xi_2 + s_0\eta_2 + \zeta}{2} \in 2{\mathbb Z}.
$$
\end{theorem}

Our proofs of {\sc Theorem} \ref{th:spinon2} and \ref{th:spinon3} are
essentially the same with that for {\sc Theorem} \ref{th:spinon1}.
Here we shall only mention a few key points.
{}From {\sc Lemma} \ref{lem:critical},  the variables $m^{(a)}_j$ 
with those $(a,j)$ in {\sc Remark} \ref{rem:critical2} are
tending to infinity.
Correspondingly we have set 
$\xi_a = p^{(a)}_{\frac{t_a}{t_r}(s-s_0)}, \
\eta_a = p^{(a)}_{\frac{t_a}{t_r}(s-s_0+t_r)}$
and $\zeta_a = p^{(a)}_s \ (t_a = t_r)$.
For $G_2$, $\zeta_2$ is simply denoted by $\zeta$ and 
$\Phi$ is in fact equal to 
$\left[\begin{array}{c}
p^{(2)}_{s+3-2s_0}+ m^{(2)}_{s+3-2s_0}\\
m^{(2)}_{s+3-2s_0} \end{array} \right]_q$.
The congruence conditions among the variables 
$\xi, \eta, \zeta$ (and $m$ for $G_2$) are equivalent to 
the condition that the following quantities are integers 
for all $1 \le a \le n$ ;
$\sum_{b=1}C^{-1}_{a b} \xi_b, \ 
\sum_{b=1}C^{-1}_{a b} \eta_b, \ 
m^{(a)}_{\frac{t_a}{t_r}(s-s_0)}, \
m^{(a)}_{\frac{t_a}{t_r}(s-s_0+t_r)}, \
m^{(a)}_s \, (t_a = t_r)$.
{}From the condition, $\phi \in {\mathbb Z}$ follows easily. 

\section{Recursion relation of fermionic forms}

The fermionic form obeys a recursion relation
compatible with the $Q$-system discussed in the next section.

\begin{theorem}\label{th:rec}
Fix $1 \le a_0 \le n$ and $j_0 \in {\mathbb Z}_{\ge 1}$
arbitrarily.
Given any $\la \in \ol{P}^+$ and 
$W = \bigotimes_{1 \le a \le n}\bigotimes_{j \ge 1}
(W^{(a)}_j)^{\otimes \nu^{(a)}_j}$, set
\begin{eqnarray*}
W_1 & = & W^{(a_0)}_{j_0} \otimes W^{(a_0)}_{j_0} \otimes W, \\
W_2 & = & W^{(a_0)}_{j_0 + 1} 
\otimes W^{(a_0)}_{j_0 - 1} \otimes W, \\
W_3 & = & \bigotimes_{b \sim a_0} \bigotimes_{k=0}^{-C_{a_0 b}-1}
W^{(b)}_{\left[ (C_{b a_0}j_0 - k)/C_{a_0 b}\right]} \otimes W,
\end{eqnarray*}
where the symbol $[x]$ denotes the largest integer
not exceeding $x$.
Then we have
\begin{eqnarray}
M(W_1,\lambda,q) & = & M(W_2,\lambda, q) 
+ q^{-\theta} M(W_3,\lambda,q), \label{th:rec1}\\
M_l(W_1,q) & = & M_l(W_2,q) + q^{-\theta}M_l(W_3,q), 
\label{th:rec2}\\
N_l(W_1,\la,q) & = & N_l(W_2,\la,q) + q^{-\theta}N_l(W_3,\la,q), 
\label{th:rec5}\\
\theta & = & j_0 + \sum_{k \ge 1} \nu^{(a_0)}_k \mbox{ min }
(j_0, k),\label{th:rec3}
\end{eqnarray}
where in (\ref{th:rec2}) and (\ref{th:rec5}) we assume 
$0 \le j \le t_al$ for all $W^{(a)}_j$'s contained in 
$W_1, W_2$ and $W_3$.
\end{theorem}
When $j_0 = 1$ we understand that 
$W_2 = W^{(a_0)}_2 \otimes W$  and similarly for $W_3$.
\begin{proof}
We prove (\ref{th:rec1}). (\ref{th:rec2}) and (\ref{th:rec5})
are similar.
By means of the properties
$\left[ \begin{array}{c} p + m
 \\   m \end{array} \right]_{q^{-1}} = 
q^{-mp}\left[ \begin{array}{c} p + m
 \\   m \end{array} \right]_q$ and 
$\left[ \begin{array}{c} p + m
 \\   m \end{array} \right]_q = 
\left[ \begin{array}{c} p + m - 1
 \\   m \end{array} \right]_q + q^p 
\left[ \begin{array}{c} p + m - 1
 \\   m - 1\end{array} \right]_q$,
we rewrite (\ref{eq:mvq1}) as
\begin{equation}\label{th:rec4}
\begin{split}
M(W_1, \la, q^{-1}) & = 
\sum_{\{ m \}}q^{\ol{c}(\{m\})}
\prod_{(a,i) \neq (a_0,j_0)} 
\left[ \begin{array}{c} p^{(a)}_i + m^{(a)}_i 
 \\   m^{(a)}_i \end{array} \right]_q \\
& \times  \left(
\left[ \begin{array}{c} p^{(a_0)}_{j_0} + m^{(a_0)}_{j_0} - 1
 \\   m^{(a_0)}_{j_0} \end{array} \right]_q + 
q^{p^{(a_0)}_{j_0}}
\left[ \begin{array}{c} p^{(a_0)}_{j_0} + m^{(a_0)}_{j_0} - 1
 \\   m^{(a_0)}_{j_0} - 1\end{array} \right]_q
\right).
\end{split}
\end{equation}
Here $\ol{c}(\{m\}) = \frac{1}{2}\sum_{1 \le a, b \le n}
(\alpha_a \vert \alpha_b)
\sum_{j, k \ge 1} \mbox{min}(t_bj, t_ak) m^{(a)}_j m^{(b)}_k$,
which does not contain $\nu^{(a)}_j$'s explicitly as opposed to 
(\ref{eq:mvq2}).
Let us write (\ref{eq:mvq3}) as 
$p^{(a)}_i(W,\{ m \})$ to exhibit 
the dependence on 
$W = \bigotimes_{a,j}
(W^{(a)}_j)^{\otimes \nu^{(a)}_j}$
and $\{m \} = \{m^{(b)}_k\}$.
Similarly, (\ref{eq:mvq4}) is denoted by (\ref{eq:mvq4})${}_W$.
Then $p^{(a)}_i = p^{(a)}_i(W_1,\{m\})$ for all $(a,i)$ in 
(\ref{th:rec4}).
Suppose $\{m\}$ obey 
(\ref{eq:mvq4})${}_{W_1}$, and set 
$m^{(b)'}_k = m^{(b)}_k - \delta_{b, a_0}\delta_{k, j_0}$.
Then it is easy to see
\begin{eqnarray*}
p^{(a)}_i(W_2,\{m\}) & = & p^{(a)}_i(W_1,\{m\}) 
- \delta_{a, a_0}\delta_{i, j_0},\quad \{m\} \mbox{ obeys } 
(\ref{eq:mvq4})_{W_2},\\
p^{(a)}_i(W_3,\{m'\}) & = & p^{(a)}_i(W_1,\{m\}),
\quad \qquad\qquad\quad\{m'\} \mbox{ obeys } 
(\ref{eq:mvq4})_{W_3},\\
\ol{c}(\{m\}) & = & \ol{c}(\{m'\}) - 
p^{(a_0)}_{j_0}(W_1,\{m\}) + \theta.
\end{eqnarray*}
Thus in the expansion (\ref{th:rec4}), 
the first and the second terms on the RHS 
yield $M(W_2,\la,q^{-1})$ and 
$q^\theta M(W_3,\la,q^{-1})$, respectively 
up to boundary effects.
By boundary effects we mean that the above
expansion should not be applied if 
$p^{(a_0)}_{j_0}(W_1,\{m\}) = 
m^{(a_0)}_{j_0} = 0$, because it amounts to allowing 
$p^{(a_0)}_{j_0}(W_2,\{m\}) = -1$ or 
$m^{(a_0)'}_{j_0} = -1$ in either interpretation
$\left[ \begin{array}{c} -1
 \\   0\end{array} \right]_q = 1$ or 
$\left[ \begin{array}{c}  - 1
 \\   - 1\end{array} \right]_q = 1$.
In fact $p^{(a_0)}_{j_0}(W_1,\{m\}) = 
m^{(a_0)}_{j_0} = 0$ does not hold for any summand in
$M(W_1,\la,q^{-1})$. This can be guaranteed by combining  
(\ref{eq:mvq3}) for $p^{(a)}_j = p^{(a)}_j(W_1,\{m\})$ as
\begin{equation*}
\begin{split}
p^{(a)}_{j+1} - 2p^{(a)}_{j} + p^{(a)}_{j-1} & =
- \nu^{(a)}_{j} - 2 \delta_{a,a_0}\delta_{j,j_0} + 2m^{(a)}_j \\
& + \sum_{b \sim a} C_{a b}
\left( \frac{t_b}{t_a}m^{(b)}_{\frac{t_b}{t_a}j} + 
\sum_{i=1}^{t_b-1}(t_b-i)(m^{(b)}_{t_bj-i} + m^{(b)}_{t_bj+i})\right),
\end{split}
\end{equation*}
where we interpret $m^{(b)}_k = 0$ whenever $k \not \in {\mathbb Z}$.
Thus $p^{(a_0)}_{j_0} = 
m^{(a_0)}_{j_0} = 0$  contradicts 
$(a,j) = (a_0,j_0)$ case of the above 
since $\forall m^{(a)}_j, p^{(a)}_j \ge 0$ and 
$C_{a b} < 0$ for $b \sim a$.
\end{proof}

For $X_n = A_n$ a similar observation has been made also in \cite{SW}.

%%%%%%%%%%%%%%%%%%%%%%%%%%%%%%%%%%%%%%%%%%%%%%%%
\section{$Q$-system}\label{sec:qsys}

Let $x_a = e^{\ol{\La}_a}$ be a complex variable.
We shall use the notations 
$e^{\alpha_a}$ and $e^{\ol{\rho}}$, etc  to stand for 
$\prod_{b=1}^n x_b^{C_{b a}}$ and $\prod_{b=1}^nx_b$, etc.
For any $\la \in \ol{P}^+$ the character
${\mbox ch } V(\la)$ of the associated irreducible finite dimensional
$X_n$-module belongs to 
${\mathbb Z}[x^{\pm 1}_1, \ldots, x^{\pm 1}_n]$.

Let $\{Q^{(a)}_j \mid 1 \le a \le n, j \in {\mathbb Z}_{\ge 0} \}$
be the functions satisfying $Q^{(a)}_0 = 1$ and 
\begin{equation}
Q^{(a)^2}_j = Q^{(a)}_{j+1}Q^{(a)}_{j-1} + 
\prod_{b \sim a} \prod_{k=0}^{-C_{a b}-1}
Q^{(b)}_{\left[\frac{C_{b a}j - k}{C_{a b}}\right]}\quad j \ge 1.
\label{eq:qsys}
\end{equation}
We shall call (\ref{eq:qsys}) $Q$-system of type $X_n$.
It first appeared in \cite{KR2} for 
$X_n = A_n, B_n, C_n$ and $D_n$ and in 
\cite{Ki1} for the exceptional case.
For the simply laced case $X_n = A_n, D_n, E_{6,7,8}$, 
it has the simple form
\begin{equation*}
Q^{(a)^2}_j
 = Q_{j+1}^{(a)}Q_{j-1}^{(a)} + \prod_{b \sim a} Q^{(b)}_j,
\end{equation*}
where the symbol $b \sim a$ is defined before 
(\ref{eq:inversecartan}).
We have listed the explicit forms for the non-simply laced 
case in Appendix \ref{app:qsys}.

The rest of this section is devoted to a 
solution of (\ref{eq:qsys}) for the classical 
algebras $X_n = A_n, B_n, C_n$ and $D_n$.
Let us introduce the functions
\begin{eqnarray}
A_n:& \chi^{(a)}_j &= \mbox{ch}V(j\ol{\La}_a), 
\label{eq:dominoa}\\
C_n:& \chi^{(a)}_j &= \begin{cases}
\sum \mbox{ch}V(k_1\ol{\La}_1+ k_2\ol{\La}_2 + \cdots + k_a\ol{\La}_a) 
&\quad 1 \le a \le n-1,\\
\mbox{ch}V(j\ol{\La}_n) &\quad a = n,
\end{cases}
\label{eq:dominoc}\\
B_n, D_n:&
\chi^{(a)}_j & = \sum 
\mbox{ch}V(k_{a_0}\ol{\La}_{a_0} + k_{a_0 + 2}\ol{\La}_{a_0 + 2} 
+ \cdots +
k_a\ol{\La}_a) \,\, 1 \le a \le n^\prime,\label{eq:dominobd}\\
&n^\prime &=  n \mbox{ for }  B_n, \, \,\, n-2  \mbox{ for } D_n,\quad
a_0 \equiv a \hbox{ mod } 2, \quad a_0 = 0 \hbox{ or } 1,
\nonumber\\
&\chi^{(a)}_j &=  \mbox{ch}V(j\ol{\La}_a) \quad a = n-1, n
\quad \hbox{ only for } D_n.\nonumber
\end{eqnarray}
Here $\ol{\La}_0 = 0$ and 
$\mbox{ch}V(\lambda)$ denotes the irreducible $X_n$ character  with 
highest weight $\lambda$.
The sum in (\ref{eq:dominoc}) is taken over non-negative integers
$k_1, \ldots, k_a$ that satisfy 
$k_1 + \cdots + k_a \le j, k_b \equiv j\delta_{a, b}$
mod 2 for all $ 1 \le b \le a$.
The sum in (\ref{eq:dominobd}) extends over non-negative integers
$k_{a_0}, k_{a_0 + 2}, \ldots, k_a$ obeying the constraint
$t_a(k_{a_0} + k_{a_0+2} + \cdots + k_{a-2}) + k_a = j$.
If one depicts the highest weights in 
the sum (\ref{eq:dominoc}) (resp. (\ref{eq:dominobd}))
by the Young diagrams as usual, they correspond to those
obtained {}from the $a \times j$ rectangular one by successively
removing $1\times 2$ (resp. $2 \times 1$) dominoes. 
(For $B_n$ with $a = n$, one needs to say a 
bit more precisely, see Appendix \ref{app:hata}.)
Now we have 
\begin{theorem}\label{th:domino}
Let $X_n = A_n, B_n, C_n$ or $D_n$. Then,
\begin{enumerate}
\renewcommand{\theenumi}{\Alph{enumi}}
\renewcommand{\labelenumi}{(\theenumi)}
\item 
$\chi^{(a)}_j = \sum_{\la \in \ol{P}^+} d_\la \mbox{ ch } V(\la),
d_\la \in {\mathbb C},\ d_{j\ol{\La}_a} = 1$,\ 
$d_\la = 0$ unless $\la \in j\ol{\La}_a - 
\sum_{b=1}^n {\mathbb Z}_{\ge 0} \alpha_b,
\quad 1 \le a \le n, j \in {\mathbb Z}_{\ge 0}$,
\item
$Q^{(a)}_j = \chi^{(a)}_j$ solves the $Q$-system (\ref{eq:qsys}),
\item
$\lim_{j \rightarrow \infty}
\left(\frac{\chi^{(a)}_j}{\chi^{(a)}_{j+1}}\right)
= x^{-1}_a$ 
in the domain $\vert e^{\alpha_1} \vert, \ldots, 
\vert e^{\alpha_n} \vert > 1$.
\end{enumerate}
\end{theorem}
(A) is trivially true. 
The functions and (\ref{eq:dominoa})-(\ref{eq:dominobd}) 
were first introduced in \cite{KR2}, where 
the statement that they satisfy the $Q$-system was
also announced.
However as they did not present a proof, we prove 
the property (B) for the 
typical case $X_n = B_n$. 
($A_n$ case is too much straightforward.)
The property (C) is needed in section \ref{sec:complete},
for which we also sketch a proof.

\begin{proof} ({\sc Theorem} \ref{th:domino} (B) for $B_n$.)
Give a partition 
$\lambda = (\la_1, \la_2, \ldots)$ with any length $l(\la)$,
define
\begin{equation}\label{eq:qladef}
Q(\la) = \det_{1 \le i,j \le l(\la)}(Q^{(1)}_{\la_i-i+j}),
\end{equation}
where we understand that $Q^{(1)}_i = 0$ for $i < 0$.
\begin{lemma}\label{lem:qbsol}
For $j \in {\mathbb Z}_{\ge 0}$ the following solves the $Q$-system 
(see Appendix C for $X_n=B_n$) 
in the polynomial ring of infinitely many variables 
$Q^{(n)}_1, Q^{(1)}_1, Q^{(1)}_2, \ldots$.
\begin{align}
Q^{(a)}_{t_aj} &= Q((j^a)) \quad 1 \le a \le n, \label{lem:qbsol1} \\
Q^{(n)}_{2j+1} &= Q^{(n)}_1\sum_{k=0}^j(-1)^k Q((j^{n-1},k)).
\label{lem:qbsol2}
\end{align}
\end{lemma}
This follows as a corollary of Theorem 5.1 in \cite{KOS} by
dropping the spectral parameter dependence.
(To get (\ref{lem:qbsol2}) one needs to expand the determinant therein.)
Thus it suffices to show
\begin{equation}\label{eq:aim}
\begin{split}
&\mbox{RHS of (\ref{lem:qbsol1}), (\ref{lem:qbsol2}) under the
substitution } 
Q^{(n)}_1 \rightarrow \chi^{(n)}_1,
Q^{(1)}_j \rightarrow \chi^{(1)}_j (\forall j)\\
&= \mbox{RHS of (\ref{eq:dominobd})}.
\end{split}
\end{equation}
In the remainder of  this section we always assume 
this substitution for $Q(\la)$ in (\ref{eq:qladef}).
A partition $\la = (\la_1, \la_2,\ldots)$ is identified with the 
Young diagram in the usual way.
We denote by 
$\la'$ its transpose, 
$\vert \la \vert = \la_1 + \la_2 + \cdots$ and 
$2\la = (2\la_1, 2\la_2, \ldots)$.
Let $s_\la$ be the Schur function of $gl_n$ associated with the 
partition $\la$. 
We denote by $LR^\la_{\mu, \nu}$ the Littlewood-Richardson
coefficient, i.e., 
$s_\mu s_\nu = \sum_\la LR^\la_{\mu, \nu}s_\la$. 
It is well known that $LR^\la_{\mu, \nu}=0$ unless
$\vert \la \vert = \vert \mu \vert + \vert \nu \vert$.
(cf. \cite{Ma})
\begin{lemma}\label{pr:qdecomp}
For a partition $\la$ with $l(\la) \le n$, 
$Q(\la)$ decomposes into the irreducible $B_n$ characters as 
\begin{equation}\label{pr:qdecomp1}
Q(\la) = \sum_{\kappa, \mu} LR^{\la}_{(2\kappa)', \mu} ch V(\mu).
\end{equation}
Here the $\kappa, \mu$ sums run over all partitions 
such that $(2\kappa)', \mu \subset \la$.
The partition $\mu = (\mu_1, \mu_2,\ldots, \mu_n)$ is identified 
with the $B_n$ highest weight $(\mu_1-\mu_2)\ol{\La}_1 + \cdots +
(\mu_{n-1}-\mu_n)\ol{\La}_{n-1} + 2\mu_n\ol{\La}_n$.
\end{lemma}
(\ref{eq:aim}) follows {}from {\sc Lemma} \ref{pr:qdecomp}.
To see this for (\ref{lem:qbsol1}), one only has to notice that 
any non-zero $LR^{(j^a)}_{(2\kappa)',\mu}$ is 1 exactly 
when $\mu$ is obtained {}from $(j^a)$ by removing $(2\kappa)'$
made of $2 \times 1$ dominoes.
A similar argument tells that 
\begin{equation*}
\sum_{k=0}^j(-1)^k Q((j^{n-1},k)) = 
\sum_{\mu \subset (j^n)}(-1)^{nj-\vert \mu \vert} ch V(\mu).
\end{equation*}
Combining this with $\chi^{(n)}_1 = ch V(\ol{\La}_n)$ 
and 
\begin{equation*}
\begin{split}
& ch V(\ol{\La}_n) ch V(\mu) \\
& = \sum_{\epsilon_1, \ldots, \epsilon_n = \pm 1, \nu \in \ol{P}^+}
ch V(\nu = \mu + \frac{\epsilon_1-\epsilon_n}{2}\ol{\Lambda}_1 + 
\ldots + \frac{\epsilon_{n-1}-\epsilon_n}{2}\ol{\Lambda}_{n-1}
+\epsilon_n \ol{\La}_n)
\end{split}
\end{equation*}
(cf. (5.1.15) in \cite{N}),
one can verify (\ref{eq:aim}) for (\ref{lem:qbsol2}).
It remains to show {\sc Lemma} \ref{pr:qdecomp}.
For this purpose we introduce 
generating functions of characters:
\begin{align*}
\frac{1}{\phi(z)} & := \frac{1+z}{\prod_{i=1}^n(1-y_iz)(1-y^{-1}_iz)}
= \sum_{j \ge 0} \chi^{(1)}_j z^j,\\
\frac{1}{\psi(z)} & := \frac{1}{(1-z^2)\phi(z)} = \sum_{j \ge 0} p_j z^j,
\end{align*}
where $y_i = x_i^{t_i}/x_{i-1}\, (1 \le i \le n, x_0 = 1)$.
Here $\chi^{(1)}_j = ch V(j\ol{\La}_1)$ is an irreducible 
$B_n$ character while $p_j$ is a reducible one in general.
\begin{lemma}\label{lem:phipsi}
\begin{align}
\frac{1}{\phi(z_1) \cdots \phi(z_n)} &= 
\sum_\la Q(\la) s_\la(z),\label{lem:phipsi1}\\
\frac{1}{\phi(z_1) \cdots \phi(z_n)} &= 
\frac{\sum_\la ch V(\la) s_\la(z)}
{\prod_{1 \le i<j \le n}(1-z_iz_j)},\label{lem:phipsi2}\\
\frac{1}
{\prod_{1 \le i<j \le n}(1-z_iz_j)} &=
\sum_\kappa s_{(2\kappa)'}(z),\label{lem:phipsi3}
\end{align}
where $\la$ (resp. $\kappa$) sum extends over all the partitions 
with $l(\la) \le n$ (resp. $l((2\kappa)') \le n$).
\end{lemma}
$Q(\la)$ and $ch V(\la)$ are characters of $B_n$ hence
the functions of $y_1,\ldots, y_n$ only.
Thus {\sc Lemma} \ref{pr:qdecomp} immediately follows
{}from {\sc Lemma} \ref{lem:phipsi}.
(\ref{lem:phipsi3}) is eq.(11.9;2) in \cite{L} or
Example 5(b) on p77 of \cite{Ma}.
Let us verify (\ref{lem:phipsi2}).
Below we write any $n \times n$ determinant 
$\det_{1 \le i,j\le n}(A_{i,j})$ simply as 
$\vert A_{i,j} \vert$.
Due to factorization of the Vandermonde type
determinant 
$\vert z^{n-j}_i-z^{n+j}_i \vert = 
\prod_{1 \le i \le j \le n}(1-z_iz_j) 
\prod_{1 \le i < j \le n}(z_i - z_j)$,
we have
\begin{equation*}
\begin{split}
&\frac{\prod_{i<j}(1-z_iz_j)(z_i-z_j)}
{\phi(z_1) \cdots \phi(z_n)} = 
\frac{\vert z^{n-j}_i-z^{n+j}_i \vert}
{\psi(z_1) \cdots \psi(z_n)}\\
&= \sum_{\sigma \in \frak{S}_n, 
\epsilon_1, \ldots, \epsilon_n = \pm 1} 
\epsilon_1 \cdots \epsilon_n \, sgn \, \sigma 
\sum_{j_1, \ldots, j_n \ge 0} z_1^{n-\epsilon_1\sigma_1+j_1} \cdots 
z_n^{n-\epsilon_n\sigma_n+j_n}p_{j_1} \cdots p_{j_n}.
\end{split}
\end{equation*}
Setting $l_i = n-\epsilon_i\sigma_i + j_i$ this equals 
$\sum_{l_1, \ldots, l_n} \vert p_{l_i-n+j} - p_{l_i-n-j} \vert
z^{l_1}_1 \cdots z^{l_n}_n$.
In view of the anti-symmetry under the interchange 
of $l_1, \ldots, l_n$, we set $l_i = \la_i+n-i$ and rewrite this as
a sum over partitions $\la = (\la_1, \ldots, \la_n)$:
$$\sum_{\la, l(\la) \le n} \vert p_{\la_i-i+j} - p_{\la_i-i-j} \vert
\vert z^{\la_i+n-i}_j \vert.$$
Thus (\ref{lem:phipsi2}) follows {}from the well known character formulae
$\vert z^{\la_i+n-i}_j \vert = s_\la(z) \prod_{i<j}(z_i-z_j)$
and $\vert p_{\la_i-i+j} - p_{\la_i-i-j} \vert = ch V(\la)$
(cf. eq.(9.13) in \cite{W}).
(\ref{lem:phipsi1}) can be proved similarly by starting with 
$\vert z^{n-j}_i \vert/(\phi(z_1) \cdots \phi(z_n))$.
\end{proof}
\begin{proof} ({\sc Theorem} \ref{th:domino}(C) for $B_n$.)
Let us first assume $a<n$.
If $a$ is even (resp. odd) (\ref{eq:dominobd}) reads
\begin{equation*}
\chi^{(a)}_m = \sum_k {\rm ch} 
V(k_a \ol{\Lambda}_a + k_{a-2} \ol{\Lambda}_{a-2} + 
\cdots + k_2 \ol{\Lambda}_2 (\mbox{resp. $k_1 \ol{\Lambda}_1$} ) ),
\end{equation*}
where the summation is taken over the set 
of non-negative integers 
$k_a,k_{a-2},\ldots$ such that $k_a + k_{a-2} + \cdots + k_2 \leq m$ 
(resp. $k_a + k_{a-2} + \cdots + k_1 = m$).
The $B_n$ character with highest weight 
$\la = (l_1 - l_2 -1) \ol{\La}_1 + \cdots +
(l_{n-1} - l_n -1) \ol{\La}_{n-1} + (2 l_n -1) \ol{\La}_n$ is
\begin{equation*}
D \cdot {\rm ch} V(\la) = \sum_{\sigma \in {\frak S}_n,
\epsilon_1, \ldots, \epsilon_n = \pm 1} {\mbox sgn}(\sigma)
\epsilon_1 \cdots \epsilon_n y^{\epsilon_1 l_1}_{\sigma_1} 
\cdots y^{\epsilon_n l_n}_{\sigma_n},
\end{equation*}
where $D$ is the Weyl denominator independent of $\la$.
First let us prove the case $a$ even, for which we have
\begin{eqnarray*}
D \cdot \chi^{(a)}_m &=& \sum_k e^{\ol{\rho}}
 (y_1 y_2)^{k_2} (y_1 y_2 y_3 y_4)^{k_4} \cdots (y_1 \cdots y_a)^{k_a}
 \nonumber\\
&\pm& (\mbox{other terms given by $y_i \rightarrow y^{\pm 1}_{\sigma_i}$}).
\end{eqnarray*}
Owing to the formula (\ref{eq:formula}) explained later, 
this expression can be recast into 
\begin{eqnarray*}
D \cdot \chi^{(a)}_m &=& e^{\ol{\rho}} 
\left\{ C_2(y_1 y_2)^m +C_4(y_1 y_2 y_3 y_4)^m+ \cdots 
+ C_a (y_1 \cdots y_a)^m 
\vphantom{\frac{1}{(1- z_1 z_2) \cdots (1-z_1 \cdots z_a)}} \right. 
\nonumber\\
& & \left. \qquad + 
\frac{1}{(1- y_1 y_2) \cdots (1-y_1 \cdots y_a)} \right\} \nonumber\\
&\pm& (\mbox{other terms given by $y_i \rightarrow y^{\pm 1}_{\sigma_i}$}),
\end{eqnarray*}
where the coefficients $C_i$'s have no dependence on $m$.
Recall that the assumption of the property (C),
$\vert e^{\alpha_1} \vert, \ldots, 
\vert e^{\alpha_n} \vert > 1$, now reads as
$\vert y_1 \vert > \cdots > \vert y_n \vert > 1$.
Therefore in the $m \rightarrow \infty$ limit we have
\begin{equation*}
D \cdot \chi^{(a)}_m = (y_1 \cdots y_a)^m (\sum (\pm ) 
\tau ( C_a e^{\ol{\rho}}) + {\mathcal O}(\gamma^m)), \quad 
\vert \gamma \vert < 1.
\end{equation*}
Here the sum is taken over those $\tau \in {\frak S}_n$ such that 
$y_{\tau_1} \cdots y_{\tau_a} = y_1 \cdots y_a$.
(Action of $\tau$ is the natural permutation of the $y$-variables.) 
As the leading term does not vanish we get
$\chi^{(a)}_m/\chi^{(a)}_{m+1} \rightarrow (y_1 \cdots y_a)^{-1}=
e^{-\ol{\La}_a}$ in this limit.

In the case $a$ odd the same 
estimation is obtained by taking the $k$-sum via
\begin{eqnarray*}
&&\sum_k e^{\ol{\rho}}(y_1)^{k_1} (y_1 y_2 y_3 )^{k_3} 
\cdots (y_1 \cdots y_a)^{k_a} \nonumber\\
&& = e^{\ol{\rho}} \left\{ \tilde{C}_1(y_1 )^m 
+\tilde{C}_3(y_1 y_2 y_3 )^m+ \cdots 
+ \tilde{C}_a (y_1 \cdots y_a)^m \right\},
\end{eqnarray*}
where the coefficients $\tilde{C}_i$'s again have no dependence on $m$.
The key here is an elementary identity
\begin{equation}\label{eq:formula}
\sum_j \xi_1^{j_1} \cdots \xi_s^{j_s} = 
\sum_{i=1}^s \frac{\xi_i^m}{\prod_{k=1,(k \ne i)}^s (1-\xi_k/\xi_i)},
\end{equation}
where the sum in the left is taken over the non-negative
integers $j_1,\ldots,j_s$ such that $j_1+\cdots+j_s=m$.
(For even $a$  we have used this with one of the variables equal to 1.)

Finally when  $a=n$, 
the factor $(y_1 \cdots y_a)^{k_a}$ in the above is replaced 
by $(y_1 \cdots y_n)^{k_n/2}$. 
The $k$-sum is taken under the condition 
$k_n + 2(k_{n-2} + \cdots + k_2) \leq m$ 
(resp. $k_n + 2(k_{n-2} + \cdots + k_1) = m$) if $n$ is even (resp. odd).
Following the similar argument to the above, one gets the
leading $m$-dependence proportional to
$(y_1 \cdots y_n)^{m/2} = e^{m\ol{\La}_n}$,
hence the desired result.
\end{proof}

%%%%%%%%%%%%%%%%%%%%%%%%%%%%%%%%%%%%%%%%%%%%%%%%
\section{Combinatorial completeness for $X_n$}\label{sec:complete}

\subsection{Main theorem}

Our goal in this section is to prove
\begin{theorem}\label{th:completeness}
Suppose that a linear combination of characters
\begin{equation*}
Q^{(a)}_j = \sum_{\la \in \ol{P}^+} d_\la \mbox{ ch } V(\la) 
\quad 1 \le a \le n, j \in {\mathbb Z}_{\ge 0}
\end{equation*}
possesses the properties:
\begin{enumerate}
\renewcommand{\theenumi}{\Alph{enumi}}
\renewcommand{\labelenumi}{(\theenumi)}
\item 
$d_\la \in {\mathbb C},\ d_{j\ol{\La}_a} = 1$,\ 
$d_\la = 0$ unless $\la \in j\ol{\La}_a - 
\sum_{b=1}^n {\mathbb Z}_{\ge 0} \alpha_b$,
\item
$\{Q^{(a)}_j\}$ satisfies the $Q$-system (\ref{eq:qsys}),
\item
$\lim_{j \rightarrow \infty}
\left(\frac{Q^{(a)}_j}{Q^{(a)}_{j+1}}\right)
= x^{-1}_a$ 
in the domain $\vert e^{\alpha_1} \vert, \ldots, 
\vert e^{\alpha_n} \vert > 1$.
\end{enumerate}
Then for arbitrary
$\nu^{(a)}_j \in {\mathbb Z}_{\ge 0}$'s we have 
\begin{equation}\label{th:completeness1}
\prod_{a=1}^n\prod_{j \ge 1}\left(Q^{(a)}_j\right)^{\nu^{(a)}_j} = 
\sum_{\la} N_\infty(W,\la,1) \mbox{ch }V(\la)
\end{equation}
for $W = \bigotimes_{a=1}^n\bigotimes_{j \ge 1}
(W^{(a)}_j)^{\otimes \nu^{(a)}_j}$.
The sum $\sum_{\la}$ runs over 
$\la \in 
(\sum_{a=1}^n \gamma^{(a)}_\infty \ol{\La}_a
 - \sum_{b=1}^n {\mathbb Z}_{\ge 0} \alpha_b) \cap 
\ol{P}^+$. 
(See (\ref{eq:nlq3}) for the definition of 
$\gamma^{(a)}_\infty$.)
\end{theorem}
We shall denote the LHS of (\ref{th:completeness1}) by 
$\mbox{ch}W$.
Before entering the proof in section \ref{sec:proof}, 
several remarks are in order.

Suppose conversely that (\ref{th:completeness1}) holds for any
$\nu^{(a)}_i$'s. Then, for the function 
$\sum_\la N_\infty(W^{(a)}_j,\la,1)\mbox{ch}V(\la)$, one can derive
(A) and (B). The former is clear {}from the definition 
(\ref{eq:nlq1}) -- (\ref{eq:nlq2}) of 
$N_\infty(W^{(a)}_j,\la,1)$.
The latter follows {}from {\sc Theorem} \ref{th:rec} (\ref{th:rec5}) 
with $l \rightarrow \infty, q = 1$.

{\sc Theorem} \ref{th:completeness} asserts the uniqueness of 
$\{Q^{(a)}_j\}$ satisfying (A) -- (C), but does not 
assure the existence, i.e. it does not 
guarantee the functions 
$\{ \sum_\la N_\infty(W^{(a)}_j,\la,1)\mbox{ch}V(\la)\}$
fulfill\ (B) and (C). (They do (A).)
In this paper we have verified the existence  
in {\sc Theorem} \ref{th:domino}
for $X_n = A_n, B_n, C_n$ and $D_n$. 
Namely the choice $Q^{(a)}_j = \chi^{(a)}_j$ 
in (\ref{eq:dominoa})--(\ref{eq:dominobd}) 
indeed fulfills (A)--(C).
For the exceptional $X_n$, our {\sc Theorem} \ref{th:completeness}
also works once the conditions (A)--(C) are shown.

By combinatorial completeness of the string hypothesis
it is usually meant that

\begin{equation}\label{eq:comp}
\mbox{ch}W = \sum_\la M(W,\la,1)\mbox{ch}V(\la),
\end{equation}
where the sum ranges over the same domain 
as in (\ref{th:completeness1}).
Note the gap between (\ref{eq:comp}) and 
(\ref{th:completeness1}) has been left as 
{\sc Conjecture} \ref{con:MN} (\ref{con:MN1}) in general.
Our proof of {\sc Theorem} \ref{th:completeness} is a generalization 
of \cite{Ki2},
where not (\ref{eq:comp}) but (\ref{th:completeness1}) 
was shown similarly for $X_n = A_n$.
So far the above mentioned gap has been filled 
only by combinatorial means for 
$X_n = A_n, \nu^{(a)}_j = \nu_j\delta_{a,1}$ \cite{KR1} 
and $\nu^{(a)}_j$ general \cite{KSS}.
Although it is not ``physical" to allow negative vacancy 
numbers, the fermionic form $N_\infty(W,\la,q)$ is not less 
interesting than $M(W,\la,q)$ in view of the remarkable
symmetry in {\sc Remark} \ref{rem:weyl} and 
{\sc Conjecture} \ref{con:nmweyl}.
Of course there is another gap between the combinatorial 
completeness and the completeness in the literal sense.
To discuss it is beyond the scope of this paper.
See for example \cite{EKS}, \cite{JD}, 
\cite{LS}, \cite{TV} and the 
references therein.

String hypothesis is an origin of the fermionic formulae.
However it should be emphasized that our 
{\sc Theorem} \ref{th:completeness} stands totally 
independent of it; it is based purely on the 
$Q$-system (B) and the properties of the $Q$-functions (A) and (C).
In a sense the $Q$-system encodes the essential aspects 
of the string hypothesis as far as the combinatoral completeness is 
concerned.

\subsection{Proof of {\sc Theorem} \ref{th:completeness}}\label{sec:proof}

\subsubsection{Definition of $\psi_i$}\label{subsubdef}

For $1 \le i \le tl+1$, let $\mathcal{Z}_i$ denote the 
set of variables
\begin{equation}\label{eq:zdef}
\mathcal{Z}_i = \{ z^{(a)}_{j,i-1}\mid (a,j) \in H_l[i] \}
= \{z^{(a)}_{j,i-1} \mid \frac{t_a}{t}(i-1)<j \le t_al, 
1 \le a \le n \}.
\end{equation}
See (\ref{eq:cartan}) for the definition of $t$.
By definition $\mathcal{Z}_{tl+1} = \emptyset$.
Suppose the variables are related by
\begin{equation}\label{eq:zsystem}
z^{(a)}_{j,i} = z^{(a)}_{j,i-1}
\prod_{b=1}^n\left(1-z^{(b)}_{\frac{t_b}{t}i, i-1} 
\right)^{-C_{b a}(j-\frac{t_a}{t}i)}
\quad (a,j) \in H_l[i+1]
\end{equation}
for $1 \le i \le tl-1$.
We assume that $z^{(b)}_{\frac{t_b}{t}i, i-1} = 0$ 
unless $\frac{t_b}{t}i \in {\mathbb Z}$, for which
the power $C_{b a}(j-\frac{t_a}{t}i)$ is an integer.
Note that $z^{(a)}_{j,i} \in \mathcal{Z}_{i+1}$ in the LHS
while all the variables 
appearing in the RHS belong to $\mathcal{Z}_i$.
Thus we have the relation 
\begin{equation}\label{eq:fields}
{\mathbb C} = {\mathbb C}(\mathcal{Z}_{tl+1}) \subset
{\mathbb C}(\mathcal{Z}_{tl}) \subset \cdots
\subset {\mathbb C}(\mathcal{Z}_1)
\end{equation}
among the fields of rational functions in $\mathcal{Z}_i$.

For $1 \le i \le tl + 1$ we define the function
$\psi_i$ by $\psi_{tl+1} = 1$ and 
\begin{equation}\label{eq:psidef}
\psi_i = \prod_{(a,j) \in H_l[i]}
\left(1-z^{(a)}_{j, \frac{t}{t_a}j-1}
\right)^{-(\gamma^{(a)}_j - \mu_a + 1)}\quad 1 \le i \le tl,
\end{equation}
which involves $\gamma^{(a)}_j$ and $\mu_a$ 
specified in (\ref{eq:nlq3}) as integer parameters.
\begin{proposition}\label{pr:psirec}
\begin{align}
\psi_i &= \psi_{i+1} \prod_{a=1}^n
\left(1-z^{(a)}_{\frac{t_a}{t}i,i-1}
\right)^{-(\gamma^{(a)}_{\frac{t_a}{t}i}-\mu_a+1)}
\quad 1 \le i \le tl,\label{eq:psirec1} \\
\psi_i &\in {\mathbb C}(\mathcal{Z}_i) \quad 1 \le i \le tl+1.
\label{eq:psirec2}
\end{align}
\end{proposition}
For $\frac{t_a}{t}i \not\in \mathbb{Z}$, $\gamma^{(a)}_{\frac{t_a}{t}i}$  
is not needed since $z^{(a)}_{\frac{t_a}{t}i,i-1} = 0$.
\begin{proof}
(\ref{eq:psirec1}) is immediate by applying the decomposition 
(\ref{eq:Hdef2}) to (\ref{eq:psidef}).
We show (\ref{eq:psirec2}) by induction with respect to $i$.
It is trivially true for $i = tl+1$.
By virtue of (\ref{eq:fields}) and (\ref{eq:psirec1}) 
we only have to check
$z^{(a)}_{\frac{t_a}{t}i,i-1} \in \mathcal{Z}_i$
(\ref{eq:zdef}), which is obvious.
\end{proof}
\subsubsection{Series expansion of $\psi_i$}\label{subsubseries}
\begin{proposition}\label{pr:psiser}
As a formal power series in the variables 
${\mathcal Z}_i$, the function 
$\psi_i \,(1 \le i \le tl+1)$ has the 
expansion
\begin{equation}\label{eq:psiser}
\psi_i = \sum_{\{m\}} \prod_{(a,j) \in H_l[i]}
\left\{ \begin{array}{c} p^{(a)}_j + m^{(a)}_j \\
m^{(a)}_j \end{array} 
\right\}_1 \left(z^{(a)}_{j,i-1}\right)^{m^{(a)}_j}.
\end{equation}
Here the sum ranges over 
$\{m^{(a)}_j \in {\mathbb Z}_{\ge 0} \mid 
(a,j) \in H_l[i] \}$ and 
$p^{(a)}_j$ is given by (\ref{eq:nlq5}).
\end{proposition}
\begin{proof}
We again invoke induction with respect to $i$.
The case $i = tl+1$ is trivially valid.
Suppose (\ref{eq:psiser}) is valid for $\psi_{i+1}$.
Substituting it into (\ref{eq:psirec1}) we get 
\begin{equation*}
\psi_i = \prod_{a=1}^n
\left(1-z^{(a)}_{\frac{t_a}{t}i,i-1}
\right)^{-(\gamma^{(a)}_{\frac{t_a}{t}i}-\mu_a+1)}
\sum_{\{m\}} \prod_{(a,j) \in H_l[i+1]}
\left\{ \begin{array}{c} p^{(a)}_j + m^{(a)}_j \\
m^{(a)}_j \end{array} 
\right\}_1 \left(z^{(a)}_{j,i}\right)^{ m^{(a)}_j}.
\end{equation*}
By further substituting (\ref{eq:zsystem}) to express
$z^{(a)}_{j,i}$ in terms of $\mathcal{Z}_i$ variables,
the RHS becomes
\begin{equation*}
\sum_{\{m\}}
\prod_{a=1}^n
\left(1-z^{(a)}_{\frac{t_a}{t}i,i-1}
\right)^{-p^{(a)}_{\frac{t_a}{t}i}-1}
 \prod_{(a,j) \in H_l[i+1]}
\left\{ \begin{array}{c} p^{(a)}_j + m^{(a)}_j \\
m^{(a)}_j \end{array} 
\right\}_1 \left(z^{(a)}_{j,i-1}\right)^{m^{(a)}_j}
\end{equation*}
with the help of (\ref{eq:nlq5}).
Expanding the first factor by means of
\begin{equation}\label{eq:expand}
(1-z)^{-p-1} = \sum_{m \ge 0} \left\{\begin{array}{c}
p+m \\ m \end{array} \right\}_1 z^m\quad 
\mbox{ for any } p \in {\mathbb Z},
\end{equation}
and using (\ref{eq:Hdef2}), we obtain (\ref{eq:psiser}).
\end{proof}
\subsubsection{Specialization of $\psi_1$}\label{subsubspecial}
$\psi_1$ is a rational function of 
$\mathcal{Z}_1 = \{ z^{(a)}_{j,0} \mid (a,j) \in H_l \}$.
Now we begin specializing these variables to a 
certain combination of the $Q$-functions.
\begin{lemma}\label{lem:zsol}
Suppose $\{Q^{(a)}_j\}$ satisfies the $Q$-system (\ref{eq:qsys}) with 
$Q^{(a)}_0=1$.
Then the recursion relation (\ref{eq:zsystem}) has a solution
\begin{equation}\label{lem:zsol1}
z^{(a)}_{j,i} = \prod_{b=1}^n\left(
Q^{(b) j - \frac{t_a}{t_b}(i_b+1)}_{i_b}
Q^{(b) -j + \frac{t_a}{t_b}i_b}_{i_b+1}\right)^{C_{b a}},\quad
i_b = \left[\frac{t_b}{t}i\right]
\end{equation}
for $(a,j) \in H_l[i+1]$ and $0 \le i \le tl-1$.
\end{lemma}
\begin{proof}
For $\frac{t_a}{t}i \in {\mathbb Z}$ 
it is straightforward (though tedious) to check that
$z^{(a)}_{\frac{t_a}{t}i,i-1}$ given by (\ref{lem:zsol1})
satisfies 
\begin{equation}\label{lem:zsol2}
1 - z^{(a)}_{\frac{t_a}{t}i,i-1} = 
\frac{Q^{(a)}_{\frac{t_a}{t}i-1}Q^{(a)}_{\frac{t_a}{t}i+1}}
{Q^{(a) 2}_{\frac{t_a}{t}i}}
\end{equation}
by using the $Q$-system.
Upon substituting (\ref{lem:zsol1}) and (\ref{lem:zsol2}), both
sides of (\ref{eq:zsystem}) can be cast into products of 
$Q$-functions, which turn out to be equal.
\end{proof}

When $i=0$ (\ref{lem:zsol1}) reads
\begin{equation}\label{lem:zsol3}
z^{(a)}_{j,0} = \prod_{b=1}^n Q^{(b) -jC_{b a}}_1  \quad 
1 \le j \le tl
\end{equation}
for the $\mathcal{Z}_1$ variables.
As said in (B) we assume the $Q$-system and set
\begin{equation}\label{eq:bigpsidef}
\Psi = \psi_1\vert_{(\ref{lem:zsol3})} \in 
\mathbb{C}(Q^{(1)}_1, \ldots, Q^{(n)}_1),
\end{equation}
which depends on the integer parameters 
$\nu^{(a)}_i$ and $\la_a$ in (\ref{eq:nlq2}) -- (\ref{eq:nlq3})
as well.
Based on the results in sections \ref{subsubdef} and 
\ref{subsubseries}, it is easy to express $\Psi$ either as
a product or a power series.
\begin{proposition}\label{pr:bigpsi}
Let $\mbox{ch}W = \prod_{a=1}^n\prod_{j \ge 1}
\left(Q^{(a)}_j\right)^{\nu^{(a)}_j}$ be the LHS of 
(\ref{th:completeness1}). We have
\begin{align}
\Psi & = \mbox{ch}W \prod_{a=1}^n
Q^{(a) -\mu_a + 1}_1\left(
\frac{Q^{(a)}_{t_al}}{Q^{(a)}_{t_al+1}}\right)^{\la_a + 1},
\label{pr:bigpsi1}\\
\Psi & = \sum_{\{u\}} \left(
\sum_{\{m\}} \prod_{(a,j) \in H_l}
\left\{ \begin{array}{c} p^{(a)}_j + m^{(a)}_j \\
m^{(a)}_j \end{array} 
\right\}_1\right)
\prod_{a=1}^n Q^{(a) - \sum_{b=1}^n C_{a b} u_b}_1.
\label{pr:bigpsi2}
\end{align}
Here the sum $\sum_{\{u\}}$ extends over 
$\{ u_a \in {\mathbb Z}_{\ge 0} \mid 1 \le a \le n \}$
and the sum $\sum_{\{m\}}$ does over
$\{ m^{(a)}_j \in {\mathbb Z}_{\ge 0} \mid 
(a,j) \in H_l \}$ obeying the constraint
$\sum_{k=1}^{t_al} k m^{(a)}_k = u_a$ for 
$1 \le a \le n$.
\end{proposition}
\begin{proof}
To show (\ref{pr:bigpsi1}) apply (\ref{lem:zsol2})
to (\ref{eq:psidef}) with $i=1$
and use the identity
\begin{equation*}
-\gamma^{(a)}_{j-1} + 2\gamma^{(a)}_j - \gamma^{(a)}_{j+1} = 
\nu^{(a)}_j
\end{equation*}
for $1 \le j \le t_al$ ($\gamma^{(a)}_0 = 0$).
The result reads
\begin{equation*}
\Psi = \mbox{ch}W \prod_{a=1}^n Q^{(a) - \mu_a + 1}_1
\frac{Q^{(a) \gamma^{(a)}_{t_al + 1} - \mu_a + 1}_{t_al}}
{Q^{(a) \gamma^{(a)}_{t_al} - \mu_a + 1}_{t_al+1}}.
\end{equation*}
{}From (\ref{eq:nlq3}) it follows that 
$\gamma^{(a)}_{t_al+1} - \mu_a + 1 = 
\gamma^{(a)}_{t_al} - \mu_a + 1 = \la_a + 1$, hence
(\ref{pr:bigpsi1}).
To show (\ref{pr:bigpsi2}) substitute (\ref{lem:zsol3}) 
into (\ref{eq:psiser}) with $i=1$ and set 
$u_a = \sum_{k=1}^{t_al} km^{(a)}_k$.
\end{proof}
\subsubsection{Integral representation of $N_l(W,\la,1)$}
\label{subsubintegral}

Anticipating the change of variables {}from 
$\{Q^{(1)}_1, \ldots, Q^{(n)}_1 \}$ to 
$\{x_1, \ldots, x_n \}$, we first prepare
\begin{lemma}\label{lem:jacobian}
\begin{equation}\label{lem:jacobian1}
\mbox{det}_{1 \le a,b\le n}
\left(\frac{\partial Q^{(a)}_1}
{\partial x_b} \right) = 
\prod_{\alpha > 0}(1-e^{-\alpha}),
\end{equation}
where the RHS denotes the product over all  
positive roots $\alpha$ of $X_n$.
\end{lemma}
\begin{proof}
Under the action of the simple reflection
$w_a=w_{\alpha_a}$ ($a=1,\ldots,n$) in the Weyl group,
the variables $e^{-\alpha_a}$ and $x_a=e^{\bar \Lambda_a}$
transforms as
$$
w_a(e^{-\alpha_b})=e^{-\alpha_b+C_{ab} \alpha_a},
\quad
w_a(x_b)=x_b e^{-\delta_{ab} \alpha_b}.
$$ 
Moreover $w_a(\alpha_a)=-\alpha_a$ and
other positive roots are transformed with each other.
The RHS of the equation 
$$
\Delta=\prod_{\alpha >0} (1-e^{-\alpha}),
$$
has an expansion $\Delta=1+O(e^{-\alpha})$ and
has a transformation property such as, 
$$
w_a(\Delta)=(-e^{\alpha_a}) \Delta.
$$
On the
other hand,  the Jacobian $J$ in the LHS
$$
J=\frac{\partial(Q^{(1)}_1,\ldots,Q^{(n)}_1)}
{\partial(x_1,\ldots,x_n)}
$$
transforms in the same manner as $\Delta$, i.e. 
$w_a(J)=-e^{\alpha_a}J$, since
the Laurent polynomials $Q^{(a)}_1$
($a=1,\ldots, n$) are Weyl group invariants
and the volume form transforms as follows
\begin{eqnarray*}
&&w_a(dx_1 \wedge \cdots \wedge dx_n)=
dx_1 \wedge \cdots \wedge dx_{a-1} \wedge
d(w_a(x_a)) \wedge dx_{a+1} \wedge \cdots \wedge dx_n \cr
&&=dx_1 \wedge \cdots dx_{a-1} \wedge
(-\frac{1}{x_a^2} dx_a \prod_{b \neq a} x_b^{-C_{ba}})
\wedge dx_{a+1} \cdots \wedge dx_n
\cr
&&=-e^{-\alpha_a} dx_1 \wedge \cdots \wedge dx_n.
\end{eqnarray*}
Hence, the ratio $J/\Delta$ is a Weyl group invariant.
Since the polynomials $Q^{(a)}_1$ have the
expansion $Q^{(a)}_1=x_a(1+O(e^{-\alpha}))$ owing to (A),
the invariant $J/\Delta$ has an
expansion $J/\Delta=1+O(e^{-\alpha})$ and hence it must be
$1$.
\end{proof}

\begin{proposition}\label{pr:integral}
Suppose $W = \bigotimes_{(a,j) \in H_l} 
\bigl(W^{(a)}_j\bigr)^{\otimes \nu^{(a)}_j}$ and 
$\la = \la_1 \ol{\La}_1 + \cdots + 
\la_n \ol{\La}_n \in \ol{P}$ are chosen so that 
$\mu_a$ in (\ref{eq:nlq3}) has the property 
$\sum_{b=1}^n C^{-1}_{a b}\mu_b \in {\mathbb Z}_{\ge 0}$.
Then we have 
\begin{align}
N_l(W,\la,1) & = \oint \cdots \oint 
\frac{dQ^{(1)}_1 \wedge \cdots \wedge dQ^{(n)}_1}
{(2\pi i)^n} \mbox{ ch}W \prod_{a=1}^n
\left(\frac{Q^{(a)}_{t_al}}{Q^{(a)}_{t_al+1}}
\right)^{\la_a + 1}\label{pr:integral1}\\
& = \oint \cdots \oint 
\frac{dx_1 \wedge \cdots \wedge dx_n}
{(2\pi i)^nx_1 \cdots x_n} \, 
 \left(e^{\ol{\rho}}\prod_{\alpha > 0}(1-e^{-\alpha})\right)
\mbox{ ch}W \prod_{a=1}^n
\left(\frac{Q^{(a)}_{t_al}}{Q^{(a)}_{t_al+1}}
\right)^{\la_a + 1}. \label{pr:integral2} 
\end{align}
Here the integrand in (\ref{pr:integral1}) 
(resp. (\ref{pr:integral2}) ) is regarded as 
an element in ${\mathbb C}(Q^{(1)}_1,\ldots,Q^{(n)}_1)$
(resp. ${\mathbb C}(x_1,\ldots,x_n)$) and the 
integration contours encircle $\infty$ 
in the domain $\prod_{b=1} \vert Q^{(b)}_1 \vert^{C_{b a}} \gg 1$
(resp. $\vert e^{\alpha_a}\vert = 
\prod_{b=1} \vert x_b \vert^{C_{b a}} \gg 1$) for
$1 \le a \le n$.
\end{proposition}
\begin{proof}
By comparing (\ref{pr:bigpsi2}) and (\ref{eq:nlq1}) --
(\ref{eq:nlq4}) with $q=1$, we find that 
$N_l(W,\la,1) = \mbox{ coefficient of }
\prod_{a=1}^nQ^{(a) -\mu_a}_1$ in $\Psi$.
The integral (\ref{pr:integral1})
is picking the coefficient up as a multi-dimensional 
residue {}from $\Psi$ in (\ref{pr:bigpsi1}).
In view of (\ref{pr:bigpsi2}) the cycles of the integrals are 
to be taken in the domain 
$\prod_{b=1} \vert Q^{(b)}_1 \vert^{C_{b a}} \gg 1$.
It follows {}from this inequality that $\vert Q^{(a)}_1\vert \gg 1$
for all $1 \le a \le n$.
(\ref{pr:integral2}) can be derived {}from 
(\ref{pr:integral1}) by changing variables with the aid of
{\sc Lemma} \ref{lem:jacobian}.
When $\forall \vert e^{\alpha_a} \vert \gg 1$, 
$\forall \vert e^{\ol{\La}_a} \vert \gg 1$ also holds.
Thus the domain of 
$\{Q^{(a)}_1\}$-integrals matches that for 
$\{x_a\}$-integrals 
because $Q^{(a)}_1 = x_a(1 + O(e^{-\alpha_a}))$ owing to (A).
Note that such a domain indeed exists, for 
example, the vicinity of 
$\vert x_a \vert = e^{y(\ol{\La}_a \vert \ol{\rho})}$ with $y \gg 0$.
\end{proof}

\subsubsection{Fermionic form $N_\infty(W,\la,1)$ as 
branching coefficient}

Define the branching coefficient $[W:\omega] \in {\mathbb C}$ by
\begin{equation}\label{eq:branching}
\mbox{ch }W = \sum_{\omega}
[W:\omega] \mbox{ ch }V(\omega),
\end{equation}
where, due to (A), the sum is actually limited to 
$\omega \in (\sum_{a=1}^n \gamma^{(a)}_\infty \ol{\La}_a - 
\sum_{a=1}^n {\mathbb Z}_{\ge 0} \alpha_a) \cap \ol{P}^+$.
Now we are ready to finish the proof of 
{\sc Theorem} \ref{th:completeness}, namely, to 
establish
\begin{equation}\label{eq:final}
[W:\la] = N_\infty(W,\la,1) \quad \mbox{ for } \la \in \ol{P}^+.
\end{equation}
Consider the limit $l \rightarrow \infty$ in 
(\ref{pr:integral2}).
Thanks to the assumption (C), we may replace 
$\prod_{a=1}^n
\left(\frac{Q^{(a)}_{t_al}}{Q^{(a)}_{t_al+1}}
\right)^{\la_a + 1}$
by $\prod_{a=1}^nx_a^{-\la_a-1} = e^{-\la - \ol{\rho}}$, 
leading to
\begin{equation}\label{eq:integral3}
N_\infty(W,\la,1) = \oint \cdots \oint 
\frac{dx_1 \wedge \cdots \wedge dx_n}
{(2\pi i)^nx_1 \cdots x_n} \, 
\left( e^{\ol{\rho}}\prod_{\alpha > 0}(1-e^{-\alpha})\right)
 \ e^{-\la - \ol{\rho}}\mbox{ ch }W.
\end{equation}
Substitute (\ref{eq:branching}) into the RHS and express
$\mbox{ch}V(\omega)$ via the Weyl-Kac character formula.
The result reads
\begin{equation}\label{eq:nint}
N_\infty(W,\la,1) = \sum_{\omega \in \ol{P}^+}
[W:\omega] \oint \cdots \oint 
\frac{dx_1 \wedge \cdots \wedge dx_n}
{(2\pi i)^nx_1 \cdots x_n} \,
e^{-\la - \ol{\rho}}\sum_w \mbox{ det }w 
\ e^{w(\omega + \ol{\rho})},
\end{equation}
where $\sum_w$ denotes the sum over the Weyl group
of $X_n$.
Since $\la, \omega \in \ol{P}^+$, the integral equals
$\delta_{w, 1} \delta_{\omega, \la}$.
Thus we have (\ref{eq:final}), and thereby complete the proof of
{\sc Theorem} \ref{th:completeness}.

\begin{remark}\label{rem:weyl}
{}From (\ref{eq:nint}) it follows that
\begin{equation*}
N_\infty(W,w(\la+\ol{\rho})-\ol{\rho},1) = 
\det w \ N_\infty(W,\la,1)
\end{equation*}
for any $\la \in \ol{P}$ and any Weyl group element $w$.
In particular, $N_\infty(W,\la,1) = 0$ if there is a root 
$\alpha$ such that 
$\langle \alpha,\la+\ol{\rho} \rangle = 0$.
In fact our computer experiments indicate 
\begin{conjecture}\label{con:nmweyl}
\begin{equation*}
N_\infty(W,w(\la+\ol{\rho})-\ol{\rho},q) = 
\det w \ N_\infty(W,\la,q) 
\quad \mbox{ for any } \la \in \ol{P}.
\end{equation*}
\end{conjecture}
\end{remark}

\begin{remark}
The integral representations 
(\ref{pr:integral1}), (\ref{pr:integral2}) and 
(\ref{eq:integral3}) of the fermionic forms  are valid 
for arbitrary $\nu^{(a)}_j \in {\mathbb C}$ and 
$\la \in {\mathbb C}\ol{\La}_1 + \cdots +
{\mathbb C}\ol{\La}_n$ as long as 
$\mu = \sum_{a,j} j\nu^{(a)}_j \ol{\La}_a - \la$ is kept in 
${\mathbb Z} \alpha_1 + \cdots +
{\mathbb Z} \alpha_n$.
In such a case  $p^{(a)}_i \in {\mathbb C}$ 
in general, and (\ref{eq:integral3}) for example  should read as
\begin{equation*}
\begin{split}
&\sum_{\{m\}}\prod_{\begin{subarray}{c} 1 \le a \le n  \\
i \ge 1 \end{subarray}}
\left\{ \begin{array}{c} p^{(a)}_i +  m^{(a)}_i
 \\   m^{(a)}_i \end{array} \right\}_1 \\
&= \oint \cdots \oint 
\frac{dx_1 \wedge \cdots \wedge dx_n}
{(2\pi i)^nx_1 \cdots x_n} \, 
\left(\prod_{\alpha > 0}(1-e^{-\alpha})\right)
 \ e^{\mu}
\prod_{a=1}^n\prod_{j \ge 1}
\left(e^{-j \ol{\La}_a}Q^{(a)}_j\right)^{\nu^{(a)}_j},
\end{split}
\end{equation*}
where the sum in LHS is taken over all 
$m^{(a)}_j \in {\mathbb Z}_{\ge 0}$ such that 
$\sum_{a=1}^n \sum_{j \ge 1} j m^{(a)}_j \alpha_a = \mu$.
In particular, the above holds as $0 = 0$  when
$\mu \not \in {\mathbb Z}_{\ge 0} \alpha_1 + \cdots +
{\mathbb Z}_{\ge 0} \alpha_n$.
\end{remark}

%%%%%%%%%%%%%%%%%%%%%%%%%%%%%%%%%%%%%%%%%%%%%%%%%%%%%%%%

\appendix
\section{List of $M(W^{(r)}_s,\lambda,q^{-1})$}\label{app:hata}
Consider a formal linear combination
$$
{\mathcal W}^{(r)}_s =
\sum_{\lambda \in \overline{P}^{+}}
M(W^{(r)}_s,\lambda,q^{-1}) V(\lambda).
$$
In this section, we give a list of ${\mathcal W}^{(r)}_s$ for
$X_n=A_n,B_n,C_n,D_n$ and ${\mathcal W}^{(r)}_1$
( ${\mathcal W}^{(r)}_2$ or a conjecture for ${\mathcal W}^{(r)}_s$ 
in some cases) for $X_n=E_{6,7,8},F_4, G_2$.
Some conjectures here have arisen from the 
calculation using the efficient algorithm in \cite{Kl}.
When $q=1$ it reduces to (\ref{eq:dominoa}) -- (\ref{eq:dominobd})
for $A_n, B_n, C_n, D_n$ and those in 
\cite{Kl} for $E_{6,7,8}$.
\newline
\newline
$X_n=A_n:$ 
\[
{{\mathcal W}^{(r)}_s}=
V({s{\ol{\La} }_r}).
\]
\newline
$X_n=B_n:$
\[
{{\mathcal W}^{(r)}_s}=
\sum_{\lambda} q^{(\ol{\La}_n | s\ol{\La}_r - \lambda)}V(\lambda)\,,
\]
where the sum $\sum_{\lambda}$ is taken over
$\lambda \in \{k_{r_0}\ol{\La}_{r_0}+k_{r_0+2}\ol{\La}_{r_0+2}+\dots
+k_{r}\ol{\La}_r\in\ol{P}^+ \, |\, t_{r}(k_{r_0}+k_{r_0+2}+\dots
+k_{r-2})+k_{r}=s\}$ with
$\ol{\La}_0=0$ and $r_0 \equiv r \pmod{2},\;r_{0}=0$ or $1$.
In the computation of
$M(W^{(r)}_s,\,k_{r_0}\ol{\La}_{r_0}+k_{r_0+2}\ol{\La}_{r_0+2}+\dots
+k_{r}\ol{\La}_r,\,q^{-1})$,
the only choice of $\{m^{(a)}_j\}$ such that 
$\forall p^{(a)}_i \ge 0$ is the following:
\newline
if $r=2u$ is even, then
\[
\begin{cases}
m^{(2a-1)}_j=\sum_{b=1}^{\min(a-1,u)} \delta_{j,l_b}+
\sum_{b=1}^{\min(a,u)} \delta_{j,l_b} &
(a \ge 1,\, 2a-1 \le n-1,\,j \ge 1)\\
m^{(2a)}_j=2 \sum_{b=1}^{\min(a,u)} \delta_{j,l_b} &
(a \ge 1,\, 2a \le n-1,\,j \ge 1)\\
m^{(n)}_j=\sum_{b=1}^{u} \delta_{j,2l_b} &
(j \ge 1)
\end{cases},
\]
and if $r=2u+1$ is odd, then
\[
\begin{cases}
m^{(1)}_j=0 & (j \ge 1)\\
m^{(2a)}_j=\sum_{b=1}^{\min(a-1,u)} \delta_{j,l_b}+
\sum_{b=1}^{\min(a,u)} \delta_{j,l_b} &
(a \ge 1,\, 2a \le n-1,\,j \ge 1)\\
m^{(2a+1)}_j=2 \sum_{b=1}^{\min(a,u)} \delta_{j,l_b} &
(a \ge 1,\, 2a+1 \le n-1,\,j \ge 1)\\
m^{(n)}_j=\sum_{b=1}^{u} \delta_{j,2l_b} &
(j \ge 1)
\end{cases},
\]
where $l_b=\frac{1}{t_r}(s-k_r)-(k_{r-2}+k_{r-4}+\dots +k_{2b+r_0})$ for 
the both cases.
\newline
\newline
$X_n=C_n:$
\[
{{\mathcal W}^{(r)}_s}=
\begin{cases}
\sum_{\lambda} q^{(\ol{\La}_n | s\ol{\La}_r - \lambda)}V(\lambda) &
(1 \le r \le n-1) \\
V({s{\ol{\La} }_n}) & (r=n)
\end{cases},
\]
where the sum $\sum_{\lambda}$ is taken over
$\lambda \in \{k_{1}\ol{\La}_{1}+\dots
+k_{r}\ol{\La}_r\in\ol{P}^+ \, |\, k_{1}+\dots +k_{r}\le s,\,
k_a \equiv s \delta_{ar}\pmod{2} \}$.
In the computation of
$M(W^{(r)}_s,\,k_{1}\ol{\La}_{1}+\dots+k_{r}\ol{\La}_r,\,q^{-1})$,
the only choice of $\{m^{(a)}_j\}$ such that 
$\forall p^{(a)}_i \ge 0$ is the following:
\[
\begin{cases}
m^{(a)}_j=\sum_{b=1}^{\min(a,r)} \delta_{j,2l_b} &
(1 \le a \le n-1,\,j \ge 1)\\
m^{(n)}_j=\sum_{b=1}^{r} \delta_{j,l_b} &
(j \ge 1)
\end{cases},
\]
where $l_b=\frac{1}{2}(s-(k_b+k_{b+1}+\dots +k_{r}))$.
\newline
\newline
$X_n=D_n:$
\[
{{\mathcal W}^{(r)}_s}=
\begin{cases}
\sum_{\lambda} q^{(\ol{\La}_n | s\ol{\La}_r - \lambda)}V(\lambda) &
(1 \le r \le n-2) \\
V({s{\ol{\La} }_r}) & (r=n-1,\,n)
\end{cases},
\]
where the sum $\sum_{\lambda}$ is taken over
$\lambda \in \{k_{r_0}\ol{\La}_{r_0}+k_{r_0+2}\ol{\La}_{r_0+2}+\dots
+k_{r}\ol{\La}_r\in\ol{P}^+ \, |\, k_{r_0}+k_{r_0+2}+\dots
+k_{r}=s\}$ with
$\ol{\La}_0=0$ and $r_0 \equiv r \pmod{2},\;r_{0}=0$ or $1$.
In the computation of
$M(W^{(r)}_s,\,k_{r_0}\ol{\La}_{r_0}+k_{r_0+2}\ol{\La}_{r_0+2}+\dots
+k_{r}\ol{\La}_r,\,q^{-1})$,
the only choice of $\{m^{(a)}_j\}$ such that 
$\forall p^{(a)}_i \ge 0$ is the following:
\newline
if $r=2u$ is even, then
\[
\begin{cases}
m^{(2a-1)}_j=\sum_{b=1}^{\min(a-1,u)} \delta_{j,l_b}+
\sum_{b=1}^{\min(a,u)} \delta_{j,l_b} &
(a \ge 1,\, 2a-1 \le n-2,\,j \ge 1)\\
m^{(2a)}_j=2 \sum_{b=1}^{\min(a,u)} \delta_{j,l_b} &
(a \ge 1,\, 2a \le n-2,\,j \ge 1)\\
m^{(n-1)}_j=m^{(n)}_j=\sum_{b=1}^{u} \delta_{j,l_b} &
(j \ge 1)
\end{cases},
\]
and if $r=2u+1$ is odd, then
\[
\begin{cases}
m^{(1)}_j=0 & (j \ge 1)\\
m^{(2a)}_j=\sum_{b=1}^{\min(a-1,u)} \delta_{j,l_b}+
\sum_{b=1}^{\min(a,u)} \delta_{j,l_b} &
(a \ge 1,\, 2a \le n-2,\,j \ge 1)\\
m^{(2a+1)}_j=2 \sum_{b=1}^{\min(a,u)} \delta_{j,l_b} &
(a \ge 1,\, 2a+1 \le n-2,\,j \ge 1)\\
m^{(n-1)}_j=m^{(n)}_j=\sum_{b=1}^{u} \delta_{j,l_b} &
(j \ge 1)
\end{cases},
\]
where $l_b=s-(k_r+k_{r-2}+\dots +k_{2b+r_0})$ for 
the both cases.

In $C_n\,(1 \le r \le n-1)$ 
(resp. $B_n\,(1 \le r \le n-1),\,D_n\,(1 \le r \le n-2)$) case,
$(\ol{\La}_n | s\ol{\La}_r - \lambda)$ is the number of
$1 \times 2$ (resp. $2 \times 1$) dominoes to be removed {}from 
$r \times s$ rectangle Young diagram to obtain $\lambda$.
In $B_n,\, r=n$ case, we consider Young diagrams consisting of 
size $1 \times \frac{1}{2}$ elementary pieces.
Then the highest weights $\la = k_n \ol{\La}_n + k_{n-2}\ol{\La}_{n-2} 
+ \cdots \in \ol{P}^+$ occurring in $\mathcal{W}^{(n)}_s$ correspond to 
those diagrams obtained from $n \times \frac{s}{2}$ shape 
by removing $2 \times 1$ blocks (made of 4 elementary pieces).
The quantity $(\ol{\La}_n \vert s \ol{\La}_n - \la)$ is the 
number of the removed $2 \times 1$ blocks to obtain $\la$.
{\allowdisplaybreaks
\begin{align*}
\intertext{$X_n=E_6:$}
{{\mathcal W}^{(1)}_1}=&
V({{\ol{\La} }_1}),
\\
{{\mathcal W}^{(2)}_1}=&
V({{\ol{\La} }_2}) + q\,V({{\ol{\La} }_5}),
\\
{{\mathcal W}^{(3)}_1}=&
{q^3}\,V(0) + q\,V({{\ol{\La} }_1}+{{\ol{\La} }_5}) + V({{\ol{\La} }_3}) +
  \left( q + {q^2} \right) \,V({{\ol{\La} }_6}),
\\
{{\mathcal W}^{(4)}_1}=&
q\,V({{\ol{\La} }_1}) + V({{\ol{\La} }_4}),
\\
{{\mathcal W}^{(5)}_1}=&
V({{\ol{\La} }_5}),
\\
{{\mathcal W}^{(6)}_1}=&
q\,V(0) + V({{\ol{\La} }_6}).
\\
\intertext{In addition we have a conjecture for ${{\mathcal W}^{(r)}_s}$:}
{{\mathcal W}^{(1)}_s} =& V(s \ol{\La}_1 ),
\\
{{\mathcal W}^{(2)}_s} =& 
\sum_{k=0}^{s} q^{k} V((s-k)\ol{\La}_2 + k\ol{\La}_5),
\\
{{\mathcal W}^{(3)}_s} =& 
\sum_{
	\begin{subarray}{c}
		j_1 + 2 j_2 + j_3 + j_4 \le s\\
		j_1,j_2,j_3,j_4 \in \mathbb{Z}_{\ge 0}
	\end{subarray}
	}
\min \left( 1+j_3,\,1+s-j_1-2 j_2-j_3-j_4 \right)\,
q^{3s-2j_1-4j_2-3j_3-2j_4} \\
& \times \begin{bmatrix} j_4+1 \\ 1
\end{bmatrix}_q
V \left( j_1(\ol{\La}_1+\ol{\La}_5) +j_2(\ol{\La}_2+\ol{\La}_4)
+j_3 \ol{\La}_3 + j_4 \ol{\La}_{6} \right),
\\
{{\mathcal W}^{(4)}_s} =& 
\sum_{k=0}^{s} q^{k} V( k\ol{\La}_1 + (s-k)\ol{\La}_4 ),
\\
{{\mathcal W}^{(5)}_s} =& V(s \ol{\La}_5 ),
\\
{{\mathcal W}^{(6)}_s} =& 
\sum_{k=0}^{s} q^{s-k} V( k\ol{\La}_6 ).
\\
\intertext{The conjecture for $\mathcal{W}^{(3)}_s$ has been checked 
for  $1 \le s \le 7$.}
\intertext{$X_n=E_7:$}
{{\mathcal W}^{(1)}_1}=&
q\,V(0) + V({{\ol{\La} }_1}),
\\
{{\mathcal W}^{(2)}_1}=&
{q^3}\,V(0) + \left( q + {q^2} \right) \,V({{\ol{\La} }_1}) +
  V({{\ol{\La} }_2}) + q\,V({{\ol{\La} }_5}),
\\
%
\if0
\mathcal{W}^{(2)}_2 =&
{q^6}\,V(0) + \left( {q^4} + {q^5} \right) \,V({{\ol{\La}}_1}) + 
  \left( {q^2} + {q^3} + {q^4} \right) \,V(2\,{{\ol{\La}}_1}) + 
  2\,{q^3}\,V({{\ol{\La}}_2}) \\& +
  V(2\,{{\ol{\La}}_2}) + 
  \left( q + {q^2} \right) \,V({{\ol{\La}}_1} + {{\ol{\La}}_2}) + 
  {q^2}\,V({{\ol{\La}}_3}) + {q^4}\,V({{\ol{\La}}_5}) \\& + 
  {q^2}\,V(2\,{{\ol{\La}}_5}) + 
  \left( {q^2} + {q^3} \right) \,V({{\ol{\La}}_1} + {{\ol{\La}}_5}) + 
  q\,V({{\ol{\La}}_2} + {{\ol{\La}}_5}),
\\
\fi
%%
%%
{{\mathcal W}^{(3)}_1}=&
\left( {q^4} + {q^6} \right)
    \,V(0) + {q^2}\,V(2{{\ol{\La} }_1}) + {q^3}\,V(2{{\ol{\La} }_6}) +
  \left( 2\,{q^3} + {q^4} + {q^5} \right) \,V({{\ol{\La} }_1})  +
  q\,V({{\ol{\La} }_1}+{{\ol{\La} }_5}) \\
  & +
  \left( q + {q^2} + {q^3} \right) \,V({{\ol{\La} }_2}) + V({{\ol{\La} }_3}) +
  \left( 2\,{q^2} + {q^3} + {q^4} \right) \,V({{\ol{\La} }_5}) +
  \left( q + {q^2} \right) \,V({{\ol{\La} }_6}+{{\ol{\La} }_7}),
\\
\mathcal{W}^{(3)}_2 =&
\left( q^8+q^{10}+q^{12} \right) \,V(0) + 
\left( 2\,q^7+q^8+3\,q^9+q^{10}+q^{11} \right) \,V({{\ol{\La}}_1}) \\& + 
  \left( 4\,{q^6} + 2\,{q^7} + 4\,{q^8} + {q^9} + {q^{10}} \right) \, 
   V(2\,{{\ol{\La}}_1}) 
   + \left( 2\,{q^5} + {q^6} + {q^7} \right) \, 
   V(3\,{{\ol{\La}}_1}) \\ & + {q^4}\,V(4\,{{\ol{\La}}_1}) + 
  \left( q^5 + 2\,q^6 + 5\,q^7 + 3\,q^8 + 2\,q^9 \right) \,V({{\ol{\La}}_2}) \\& + 
  \left( {q^2} + {q^3} + 3\,{q^4} + {q^5} + {q^6}\
     \right) \,V(2\,{{\ol{\La}}_2}) + 
  \left( 2\,{q^4} + 5\,{q^5} + 5\,{q^6} + 3\,{q^7} + {q^8} \right) \, 
   V({{\ol{\La}}_1} + {{\ol{\La}}_2}) \\& + 
  \left( {q^3} + {q^4} + {q^5} \right) \, 
   V(2\,{{\ol{\La}}_1} + {{\ol{\La}}_2}) + 
   \left( 2\,q^4 + 3\,q^5 + 7\,q^6 + 2\,q^7 + q^8 \right) \,V({{\ol{\La}}_3}) \\& + 
  V(2\,{{\ol{\La}}_3}) + \left( 3\,{q^3} + 3\,{q^4} + 3\,{q^5} \right) \, 
   V({{\ol{\La}}_1} + {{\ol{\La}}_3}) + 
  {q^2}\,V(2\,{{\ol{\La}}_1} + {{\ol{\La}}_3}) \\& + 
  \left( q + {q^2} + {q^3} \right) \,V({{\ol{\La}}_2} + {{\ol{\La}}_3}) + 
  {q^3}\,V(2\,{{\ol{\La}}_4}) + 
  \left( 2\,q^6 + q^7 + 4\,q^8 + q^9 + q^{10} \right) \,V({{\ol{\La}}_5}) \\& + 
  \left( 3\,{q^4} + 2\,{q^5} + 4\,{q^6} + {q^7} + {q^8} \right) \, 
   V(2\,{{\ol{\La}}_5}) + \left( 5\,{q^5} + 6\,{q^6} + 7\,{q^7} + 2\,{q^8} + 
     {q^9} \right) \,V({{\ol{\La}}_1} + {{\ol{\La}}_5}) \\& + 
  \left( 4\,{q^4} + 2\,{q^5} + 2\,{q^6} \right) \, 
   V(2\,{{\ol{\La}}_1} + {{\ol{\La}}_5}) + 
  {q^3}\,V(3\,{{\ol{\La}}_1} + {{\ol{\La}}_5}) \\& + 
  \left( 2\,{q^3} + 5\,{q^4} + 7\,{q^5} + 3\,{q^6} + {q^7} \right) \, 
   V({{\ol{\La}}_2} + {{\ol{\La}}_5}) \\& + 
  \left( {q^2} + 2\,{q^3} + {q^4} \right) \, 
   V({{\ol{\La}}_1} + {{\ol{\La}}_2} + {{\ol{\La}}_5}) + 
  \left( 2\,{q^2} + 2\,{q^3} + 2\,{q^4} \right) \, 
   V({{\ol{\La}}_3} + {{\ol{\La}}_5}) \\& + 
  q\,V({{\ol{\La}}_1} + {{\ol{\La}}_3} + {{\ol{\La}}_5}) + 
  \left( 2\,{q^3} + {q^4} + {q^5} \right) \, 
   V({{\ol{\La}}_1} + 2\,{{\ol{\La}}_5}) + 
  {q^2}\,V(2\,{{\ol{\La}}_1} + 2\,{{\ol{\La}}_5}) \\& + 
  \left( {q^7} + {q^9} \right) \,V(2\,{{\ol{\La}}_6}) + 
  {q^6}\,V(4\,{{\ol{\La}}_6}) + 
  \left( {q^4} + 5\,{q^5} + 4\,{q^6} + 2\,{q^7} \right) \, 
   V({{\ol{\La}}_4} + {{\ol{\La}}_6}) \\& + 
  \left( 2\,{q^3} + 3\,{q^4} + {q^5} \right) \, 
   V({{\ol{\La}}_1} + {{\ol{\La}}_4} + {{\ol{\La}}_6}) + 
  {q^2}\,V({{\ol{\La}}_2} + {{\ol{\La}}_4} + {{\ol{\La}}_6}) \\& + 
  \left( 3\,{q^6} + {q^7} + {q^8} \right) \, 
   V({{\ol{\La}}_1} + 2\,{{\ol{\La}}_6}) + 
  {q^5}\,V(2\,{{\ol{\La}}_1} + 2\,{{\ol{\La}}_6}) \\& + 
  \left( 2\,{q^4} + 2\,{q^5} + {q^6} \right) \, 
   V({{\ol{\La}}_2} + 2\,{{\ol{\La}}_6}) + 
  2\,{q^3}\,V({{\ol{\La}}_3} + 2\,{{\ol{\La}}_6}) \\& + 
  \left( 2\,{q^5} + {q^6} + {q^7} \right) \, 
   V({{\ol{\La}}_5} + 2\,{{\ol{\La}}_6}) + 
  {q^4}\,V({{\ol{\La}}_1} + {{\ol{\La}}_5} + 2\,{{\ol{\La}}_6}) \\& + 
  \left( 2\,{q^5} + {q^6} + {q^7} \right) \,V(2\,{{\ol{\La}}_7}) + 
  \left( {q^3} + 4\,{q^4} + 3\,{q^5} + {q^6} \right) \, 
   V({{\ol{\La}}_4} + {{\ol{\La}}_7}) \\& + 
  \left( {q^2} + {q^3} \right) \, 
   V({{\ol{\La}}_1} + {{\ol{\La}}_4} + {{\ol{\La}}_7}) + 
  \left( {q^5} + 3\,{q^6} + 4\,{q^7} + 2\,{q^8} \right) \, 
   V({{\ol{\La}}_6} + {{\ol{\La}}_7}) \\& + 
  \left( 3\,{q^4} + 6\,{q^5} + 4\,{q^6} + {q^7} \right) \, 
   V({{\ol{\La}}_1} + {{\ol{\La}}_6} + {{\ol{\La}}_7}) + 
  \left( {q^3} + {q^4} \right) \, 
   V(2\,{{\ol{\La}}_1} + {{\ol{\La}}_6} + {{\ol{\La}}_7}) \\& + 
  \left( {q^2} + 3\,{q^3} + 3\,{q^4} + {q^5} \right) \, 
   V({{\ol{\La}}_2} + {{\ol{\La}}_6} + {{\ol{\La}}_7}) + 
  \left( q + {q^2} \right) \,V({{\ol{\La}}_3} + {{\ol{\La}}_6} + 
     {{\ol{\La}}_7}) \\& + 
  \left( 2\,{q^3} + 4\,{q^4} + 3\,{q^5} + {q^6}\
     \right) \,V({{\ol{\La}}_5} + {{\ol{\La}}_6} + {{\ol{\La}}_7}) \\& + 
  \left( {q^2} + {q^3} \right) \, 
   V({{\ol{\La}}_1} + {{\ol{\La}}_5} + {{\ol{\La}}_6} + {{\ol{\La}}_7}) + 
  \left( {q^4} + {q^5} \right) \,V(3\,{{\ol{\La}}_6} + {{\ol{\La}}_7}) \\& + 
  {q^4}\,V({{\ol{\La}}_1} + 2\,{{\ol{\La}}_7}) + 
  {q^3}\,V({{\ol{\La}}_5} + 2\,{{\ol{\La}}_7}) + 
  \left( {q^2} + {q^3} + {q^4} \right) \, 
   V(2\,{{\ol{\La}}_6} + 2\,{{\ol{\La}}_7}),
\\
{{\mathcal W}^{(4)}_1}=&
q\,V({{\ol{\La} }_1}+{{\ol{\La} }_6}) + V({{\ol{\La} }_4}) +
  \left( {q^2} + {q^3} \right) \,V({{\ol{\La} }_6}) +
  \left( q + {q^2} \right) \,V({{\ol{\La} }_7}),
\\
%
\if0
\mathcal{W}^{(4)}_2 =&
{q^4}\,V({{\ol{\La}}_2}) + 2\,{q^3}\,V({{\ol{\La}}_3}) + 
  V(2\,{{\ol{\La}}_4}) + {q^5}\,V({{\ol{\La}}_5}) + 
  \left( {q^3} + {q^4} \right) \,V({{\ol{\La}}_1} + {{\ol{\La}}_5}) \\& + 
  {q^2}\,V({{\ol{\La}}_2} + {{\ol{\La}}_5}) + 
  \left( {q^4} + {q^5} + {q^6} \right) \,V(2\,{{\ol{\La}}_6}) + 
  \left( {q^2} + {q^3} \right) \,V({{\ol{\La}}_4} + {{\ol{\La}}_6}) \\& + 
  q\,V({{\ol{\La}}_1} + {{\ol{\La}}_4} + {{\ol{\La}}_6}) + 
  \left( {q^3} + {q^4} \right) \,V({{\ol{\La}}_1} + 2\,{{\ol{\La}}_6}) + 
  {q^2}\,V(2\,{{\ol{\La}}_1} + 2\,{{\ol{\La}}_6}) \\& + 
  \left( {q^2} + {q^3} + {q^4} \right) \,V(2\,{{\ol{\La}}_7}) + 
  \left( q + {q^2} \right) \,V({{\ol{\La}}_4} + {{\ol{\La}}_7}) + 
  \left( {q^3} + 2\,{q^4} + {q^5} \right) \, 
   V({{\ol{\La}}_6} + {{\ol{\La}}_7}) \\& + 
  \left( {q^2} + {q^3} \right) \, 
   V({{\ol{\La}}_1} + {{\ol{\La}}_6} + {{\ol{\La}}_7}),
\\
\fi
%%
%%
{{\mathcal W}^{(5)}_1}=&
{q^2}\,V(0) + q\,V({{\ol{\La} }_1}) + V({{\ol{\La} }_5}),
\\
%
\if0
\mathcal{W}^{(5)}_2 =&
{q^4}\,V(0) + {q^3}\,V({{\ol{\La}}_1}) + {q^2}\,V(2\,{{\ol{\La}}_1}) + 
  {q^2}\,V({{\ol{\La}}_5}) + V(2\,{{\ol{\La}}_5}) + 
  q\,V({{\ol{\La}}_1} + {{\ol{\La}}_5}),
\\
\fi
%%
%%
{{\mathcal W}^{(6)}_1}=&
V({{\ol{\La} }_6}),
\\
{{\mathcal W}^{(7)}_1}=&
q\,V({{\ol{\La} }_6}) + V({{\ol{\La} }_7}).
\\
\intertext{In addition we have a conjecture for
${{\mathcal W}^{(r)}_s}\, (r = 1,2,4,5,6,7)$:}
{{\mathcal W}^{(1)}_s} =& \sum_{k=0}^{s} q^{s-k} V(k \ol{\La}_1 ),
\\
{{\mathcal W}^{(2)}_s} =& 
\sum_{
	\begin{subarray}{c}
		j_1 + j_2 + 2 j_3 + j_4 \le s\\
		j_1,j_2,j_3,j_4 \in \mathbb{Z}_{\ge 0}
	\end{subarray}
	}
\min \left( 1+j_2,\,1+s-j_1-j_2-2 j_3-j_4 \right)\,
q^{3s-2j_1-3j_2-4j_3-2j_4} \\
& \times \begin{bmatrix} j_1+1 \\ 1
\end{bmatrix}_q
V \left( j_1 \ol{\La}_1 +j_2 \ol{\La}_2
+j_3 \ol{\La}_3 + j_4 \ol{\La}_{5} \right),
\\
{{\mathcal W}^{(4)}_s} =& 
\sum_{\{j\}} (j_{\ol{\La}_3}+1) q^{c(\{j\})}
\\
& \times 
\begin{bmatrix} j_{\ol{\La}_6}+1 \\ 1 \end{bmatrix}_q
\begin{bmatrix} j_{\ol{\La}_7}+1 \\ 1 \end{bmatrix}_q
\begin{bmatrix} j_{\ol{\La}_1+\ol{\La}_5}+1 \\ 1 \end{bmatrix}_q
V \left( \sum_{\lambda \in T_1 \cup T_2 \cup T_3} 
j_{\lambda} \lambda \right),
\intertext{
where the sum $\sum_{\{j\}}$ is taken over
$\{j_\lambda \in \mathbb{Z}_{\ge 0} \, |\,
\lambda \in T_1 \cup T_2 \cup T_3 \}$ under the conditions
(1)
$\sum_{\lambda \in T_1} j_\lambda +
2 \sum_{\lambda \in T_2} j_\lambda +
3 \sum_{\lambda \in T_3} j_\lambda =s \,$,
(2) $j_{\ol{\La}_5} j_{\ol{\La}_1+\ol{\La}_6} = 0\,$,
(3) $j_{\ol{\La}_2+\ol{\La}_5} j_{\ol{\La}_1+\ol{\La}_4} = 0\,$.
Here $T_1=\{\, \ol{\La}_1+\ol{\La}_6,\, \ol{\La}_4,\,
\ol{\La}_6,\, \ol{\La}_7\,\},\,
T_2=\{\, \ol{\La}_2,\, \ol{\La}_3,\, \ol{\La}_5,\,
\ol{\La}_1+\ol{\La}_5,\, \ol{\La}_2+\ol{\La}_5\, \},\,
T_3=\{\, \ol{\La}_1+\ol{\La}_4,\, \ol{\La}_2+\ol{\La}_4\, \}$
and
$c(\{j\})=3s-2j_{\ol{\La}_1+\ol{\La}_6}-
3j_{\ol{\La}_4}-j_{\ol{\La}_6}-
2j_{\ol{\La}_7}-2j_{\ol{\La}_2}-3j_{\ol{\La}_3}-j_{\ol{\La}_5}
-3j_{\ol{\La}_1+\ol{\La}_5}-4j_{\ol{\La}_2+\ol{\La}_5}-
3j_{\ol{\La}_1+\ol{\La}_4}-4j_{\ol{\La}_2+\ol{\La}_4}$.
}
{{\mathcal W}^{(5)}_s} =&
\sum_{
	\begin{subarray}{c}
		j+k \le s\\
		j,k \in \mathbb{Z}_{\ge 0}
	\end{subarray}
	}
q^{2s-2k-j} V( j \ol{\La}_1 + k \ol{\La}_5 ),
\\
{{\mathcal W}^{(6)}_s} =& V( s \ol{\La}_6 ),
\\
{{\mathcal W}^{(7)}_s} =&
\sum_{k=0}^{s} q^{k} V \left( k \ol{\La}_6 + (s-k) \ol{\La}_7 \right).
\\
\intertext{The conjecture for $\mathcal{W}^{(2)}_s$ 
($\mathcal{W}^{(4)}_s$) has been checked 
for  $1 \le s \le 7$ ($1 \le s \le 8$).}
\intertext{$X_n=E_8:$}
{{\mathcal W}^{(1)}_1}=&
q\,V(0) + V({{\ol{\La} }_1}),
\\
{{\mathcal W}^{(2)}_1}=&
{q^3}\,V(0) + \left( q + {q^2} \right) \,V({{\ol{\La} }_1}) +
  V({{\ol{\La} }_2}) + q\,V({{\ol{\La} }_7}),
\\
{{\mathcal W}^{(3)}_1}=&
\left( {q^4} + {q^6} \right) \,V(0) + {q^2}\,V(2{{\ol{\La} }_1}) +
  \left( 2\,{q^3} + {q^4} + {q^5} \right) \,V({{\ol{\La} }_1}) +
  q\,V({{\ol{\La} }_1}+{{\ol{\La} }_7})\\
  & +
  \left( q + {q^2} + {q^3} \right) \,V({{\ol{\La} }_2}) + V({{\ol{\La} }_3}) +
  \left( 2\,{q^2} + {q^3} + {q^4} \right) \,V({{\ol{\La} }_7}) +
  \left( q + {q^2} \right) \,V({{\ol{\La} }_8}),
\\
{{\mathcal W}^{(4)}_1}=&
\left( {q^6} + {q^7} + {q^8} + {q^{10}} \right) \,V(0) +
  \left( {q^3} + 2\,{q^4} + {q^5} + {q^6} \right) \,V(2{{\ol{\La} }_1}) + {q^2}
      \,V(2{{\ol{\La} }_7})\\
  & +
  \left( {q^4} + 3\,{q^5} + 3\,{q^6} + 2\,{q^7} + {q^8} + {q^9} \right)
      \,V({{\ol{\La} }_1}) +
  \left( {q^2} + {q^3} \right) \,V({{\ol{\La} }_1}+{{\ol{\La} }_2}) \\
  & +
  \left( 2\,{q^2} + 3\,{q^3} + 2\,{q^4} + {q^5} \right)
      \,V({{\ol{\La} }_1}+{{\ol{\La} }_7}) +
  \left( q + {q^2} \right) \,V({{\ol{\La} }_1}+{{\ol{\La} }_8})\\
  & +
  \left( 3\,{q^3} + 3\,{q^4} + 3\,{q^5} + {q^6} + {q^7} \right)
      \,V({{\ol{\La} }_2}) + q\,V({{\ol{\La} }_2}+{{\ol{\La} }_7})\\
  & +
  \left( q + 2\,{q^2} + {q^3} + {q^4} \right) \,V({{\ol{\La} }_3}) +
  V({{\ol{\La} }_4}) + \left( q + {q^2} + {q^3} \right) \,V({{\ol{\La} }_6})\\
  & +
  \left( {q^3} + 4\,{q^4} + 3\,{q^5} + 2\,{q^6} + {q^7} + {q^8} \right)
      \,V({{\ol{\La} }_7}) +
  \left( 2\,{q^2} + 3\,{q^3} + 2\,{q^4} + {q^5} + {q^6} \right)
      \,V({{\ol{\La} }_8}),
\\
{{\mathcal W}^{(5)}_1}=&
\left( {q^7} + 3\,{q^9} + {q^{10}} +
        2\,{q^{11}} + {q^{12}} + {q^{13}} + {q^{15}} \right) \,V(0) \\
  & +
  \left( 5\,{q^5} + 4\,{q^6} + 6\,{q^7} + 3\,{q^8} +
          3\,{q^9} + {q^{10}} + {q^{11}} \right) \,V(2{{\ol{\La} }_1})\\
  & +
  \left( 2\,{q^3} + {q^4} + {q^5} \right) \,V(
      2{{\ol{\La} }_1}+{{\ol{\La} }_7}) + {q^3}\,V(2{{\ol{\La} }_2})\\
  & +
  \left( 3\,{q^3} + 2\,{q^4} + 3\,{q^5} + {q^6} + {q^7} \right) \,V(
      2{{\ol{\La} }_7}) + \left( {q^4} + {q^6} \right) \,V(3{{\ol{\La} }_1})\\
  & +
  \left( 4\,{q^6} + 5\,{q^7} + 8\,{q^8} + 5\,{q^9} + 5\,{q^{10}} +
          3\,{q^{11}} + 2\,{q^{12}} + {q^{13}} + {q^{14}} \right)
      \,V({{\ol{\La} }_1})\\
  & + {q^2}\,V({{\ol{\La} }_1}+2{{\ol{\La} }_7}) +
  \left( 2\,{q^3} + 5\,{q^4} + 5\,{q^5} + 3\,{q^6} + 2\,{q^7} + {q^8} \right)
      \,V({{\ol{\La} }_1}+{{\ol{\La} }_2})\\
  & +
  \left( 2\,{q^2} + {q^3} + {q^4} \right)
      \,V({{\ol{\La} }_1}+{{\ol{\La} }_3}) +
  \left( q + {q^2} + {q^3} \right) \,V({{\ol{\La} }_1}+{{\ol{\La} }_6})\\
  & +
  \left( 2\,{q^3} + 9\,{q^4} + 10\,{q^5} + 10\,{q^6} + 6\,{q^7} + 4\,{q^8} +
          2\,{q^9} + {q^{10}} \right) \,V({{\ol{\La} }_1}+{{\ol{\La} }_7})\\
  & +
  \left( 2\,{q^2} + 5\,{q^3} + 5\,{q^4} + 3\,{q^5} + 2\,{q^6} + {q^7} \right)
      \,V({{\ol{\La} }_1}+{{\ol{\La} }_8})\\
  & +
  \left( 3\,{q^4} + 6\,{q^5} + 11\,{q^6} + 8\,{q^7} + 7\,{q^8} + 4\,{q^9} +
          3\,{q^{10}} + {q^{11}} + {q^{12}} \right) \,V({{\ol{\La} }_2})\\
  & +
  \left( 3\,{q^2} + 4\,{q^3} + 4\,{q^4} + 2\,{q^5} + {q^6} \right)
      \,V({{\ol{\La} }_2}+{{\ol{\La} }_7}) +
  \left( q + {q^2} \right) \,V({{\ol{\La} }_2}+{{\ol{\La} }_8})\\
  & +
  \left( 5\,{q^3} + 5\,{q^4} + 7\,{q^5} + 4\,{q^6} +
          3\,{q^7} + {q^8} + {q^9} \right) \,V({{\ol{\La} }_3}) +
  q\,V({{\ol{\La} }_3}+{{\ol{\La} }_7})\\
  & +
  \left( q + 2\,{q^2} + 2\,{q^3} + {q^4} + {q^5} \right)
      \,V({{\ol{\La} }_4}) + V({{\ol{\La} }_5})\\
  & +
  \left( 2\,{q^2} + 4\,{q^3} + 6\,{q^4} + 4\,{q^5} +
          3\,{q^6} + {q^7} + {q^8} \right) \,V({{\ol{\La} }_6})\\
  & +
  \left( 5\,{q^5} + 6\,{q^6} + 9\,{q^7} + 6\,{q^8} + 6\,{q^9} + 3\,{q^{10}} +
          2\,{q^{11}} + {q^{12}} + {q^{13}} \right) \,V({{\ol{\La} }_7})\\
  & +
  \left( q + 2\,{q^2} + 2\,{q^3} + {q^4} \right)
      \,V({{\ol{\La} }_7}+{{\ol{\La} }_8})\\
  & +
  \left( {q^3} + 5\,{q^4} + 8\,{q^5} + 7\,{q^6} + 6\,{q^7} + 4\,{q^8} +
          2\,{q^9} + {q^{10}} + {q^{11}} \right) \,V({{\ol{\La} }_8}),
\\
{{\mathcal W}^{(6)}_1}=&
\left( {q^5} + {q^7} \right) \,V(0) + {q^3}\,V(2{{\ol{\La} }_1}) +
  \left( {q^3} + 2\,{q^4} + {q^5} + {q^6} \right) \,V({{\ol{\La} }_1}) +
  \left( q + {q^2} \right) \,V({{\ol{\La} }_1}+{{\ol{\La} }_7})
  \\& +
  \left( 2\,{q^2} + {q^3} + {q^4} \right) \,V({{\ol{\La} }_2}) +
  q\,V({{\ol{\La} }_3}) + V({{\ol{\La} }_6}) +
  \left( {q^2} + 2\,{q^3} + {q^4} + {q^5} \right) \,V({{\ol{\La} }_7})\\
  & +
  \left( q + {q^2} + {q^3} \right) \,V({{\ol{\La} }_8}),
\\
{{\mathcal W}^{(7)}_1}=&
{q^2}\,V(0) + q\,V({{\ol{\La} }_1}) + V({{\ol{\La} }_7}),
\\
{{\mathcal W}^{(8)}_1}=&
{q^4}\,V(0) + \left( {q^2} + {q^3} \right) \,V({{\ol{\La} }_1}) +
  q\,V({{\ol{\La} }_2}) + \left( q + {q^2} \right) \,V({{\ol{\La} }_7}) +
  V({{\ol{\La} }_8}).
\\
\intertext{In addition we have a conjecture for
${{\mathcal W}^{(r)}_s}\, (r=1,2,7)$:}
{{\mathcal W}^{(1)}_s} =&
\sum_{k=0}^{s} q^{s-k} V(k\ol{\La}_1),
\\
{{\mathcal W}^{(2)}_s} =& 
\sum_{
	\begin{subarray}{c}
		j_1 + j_2 + 2 j_3 + j_4 \le s\\
		j_1,j_2,j_3,j_4 \in \mathbb{Z}_{\ge 0}
	\end{subarray}
	}
\min \left( 1+j_2,\,1+s-j_1-j_2-2 j_3-j_4 \right)\,
q^{3s-2j_1-3j_2-4j_3-2j_4} \\
& \times \begin{bmatrix} j_1+1 \\ 1
\end{bmatrix}_q
V \left( j_1 \ol{\La}_1 +j_2 \ol{\La}_2
+j_3 \ol{\La}_3 + j_4 \ol{\La}_{7} \right),
\\
{{\mathcal W}^{(7)}_s} =&
\sum_{
	\begin{subarray}{c}
		j+k \le s\\
		j,k \in \mathbb{Z}_{\ge 0}
	\end{subarray}
	}
q^{2s-2k-j} V( j \ol{\La}_1 + k \ol{\La}_7 ).
\\
\intertext{The conjecture for $\mathcal{W}^{(2)}_s$ has been checked 
for  $1 \le s \le 7$.}
\intertext{$X_n=F_4:$}
{{\mathcal W}^{(1)}_1}=&
q\,V(0) + V({{\ol{\La} }_1}),
\\
{{\mathcal W}^{(2)}_1}=&
{q^3}\,V(0) + \left( q + {q^2} \right) \,V({{\ol{\La} }_1}) +
  V({{\ol{\La} }_2}) + q\,V({2{\ol{\La} }_4}),
\\
{{\mathcal W}^{(3)}_1}=&
V({{\ol{\La} }_3}) + q\,V({{\ol{\La} }_4}),
\\
\mathcal{W}^{(3)}_2 =&
\left( {q^4} + {q^6} \right) \,V(0) + 
  \left( 2\,{q^3} + {q^4} + {q^5} \right) \,V({{\ol{\La}}_1}) + 
  {q^2}\,V(2\,{{\ol{\La}}_1}) \\& + 
  \left( q + {q^2} + {q^3} \right) \,V({{\ol{\La}}_2}) + 
  {q^3}\,V({{\ol{\La}}_3}) + V(2\,{{\ol{\La}}_3}) + 
  \left( 2\,{q^2} + {q^3} + {q^4} \right) \,V(2\,{{\ol{\La}}_4}) \\& + 
  \left( q + {q^2} \right) \,V({{\ol{\La}}_3} + {{\ol{\La}}_4}) + 
  q\,V({{\ol{\La}}_1} + 2\,{{\ol{\La}}_4}),
\\
{{\mathcal W}^{(4)}_1}=&
V({{\ol{\La} }_4}).
\\
\intertext{In addition we have a conjecture for
${{\mathcal W}^{(r)}_s}\, (r = 1,2,4)$:}
{{\mathcal W}^{(1)}_s} =&
\sum_{k=0}^{s} q^{s-k} V(k\ol{\La}_1),
\\
{{\mathcal W}^{(2)}_s} =& 
\sum_{
	\begin{subarray}{c}
		j_1 + j_2 + 2 j_3 + j_4 \le s\\
		j_1,j_2,j_3,j_4 \in \mathbb{Z}_{\ge 0}
	\end{subarray}
	}
\min \left( 1+j_2,\,1+s-j_1-j_2-2 j_3-j_4 \right)\,
q^{3s-2j_1-3j_2-4j_3-2j_4} \\
& \times \begin{bmatrix} j_1+1 \\ 1
\end{bmatrix}_q
V \left( j_1 \ol{\La}_1 +j_2 \ol{\La}_2
+2 j_3 \ol{\La}_3 + 2 j_4 \ol{\La}_{4} \right),
\\
{{\mathcal W}^{(4)}_s} =&
\sum_{k=0}^{[s/2]} \sum_{j=0}^{k} q^{2k-j}
V \left( j \ol{\La}_1 + (s-2k)\ol{\La}_4 \right).
\\
\intertext{The conjecture for $\mathcal{W}^{(2)}_s$ has been checked 
for  $1 \le s \le 5$.}
\intertext{$X_n=G_2:$}
{{\mathcal W}^{(1)}_1}=&
q\,V(0) + V({{\ol{\La} }_1}),
\\
{{\mathcal W}^{(2)}_1}=&
V({{\ol{\La} }_2}).
\\
\intertext{In addition we have a conjecture for ${{\mathcal W}^{(r)}_s}$:}
{{\mathcal W}^{(1)}_s} =&
\sum_{k=0}^{s} q^{s-k} V(k\ol{\La}_1)
\\
{\mathcal W}^{(2)}_s = &
\sum_{k=0}^{[s/3]} \sum_{x=2k}^{s-k}
\left\{
\min\left(
	\left[ \frac{x-2k}{3} \right],\,
	\left[ \frac{s+k-2x}{3} \right] + \left[ \frac{x-2k}{3} \right]
	\right)
+1 \right\}\\
&\quad\quad\quad\quad\quad\quad \times q^{x-k}
\begin{bmatrix} k+1 \\ 1
\end{bmatrix}_q
V(k \ol{\La}_1+(s-x-k)\ol{\La}_2).
\\
\end{align*}}
The last conjecture has been checked for $1 \le s \le 16$.

\section{Example of one dimensional sums}\label{app:taka}
We present an example of the one dimensional configuration sums (1dsums)
and the fermionic forms for $C_2^{(1)}$.
First we consider the 1dsums.
Let us take a $U'_q(C_2^{(1)})$-crystal $B:=
B^{1,2} \otimes (B^{2,1})^{\otimes 3} 
\otimes (B^{1,1})^{\otimes 2}$.
As $U_q(C_2)$-crystals, $B^{1,2} \simeq B(2 \ol{\La}_1) \oplus B(0)$,
$B^{2,1} \simeq B(\ol{\La}_2) $ and $B^{1,1} \simeq B(\ol{\La}_1) $.
We use the parametrization of the $U_q(C_2)$-crystal $B(s \ol{\La}_r)$ 
given in \cite{N}. 
We denote the element in $B(0)$ by $\phi$, and
the elements in  $B(\ol{\La}_2) $ by $a,b,c,d,e$ for short, where
\vspace*{-4pt}
\begin{center}
$a=$
\begin{picture}(17,28)(0,12)
	\put( 0, 0){\line(0,1){28}}
	\put(14, 0){\line(0,1){28}}
	\put( 0, 0){\line(1,0){14}}
	\put( 0,14){\line(1,0){14}}
	\put( 0,28){\line(1,0){14}}
	\put( 7,21){\makebox(0,0)[c]{$1$}}
	\put( 7, 7){\makebox(0,0)[c]{$2$}}
\end{picture},\quad
$b=$
\begin{picture}(17,28)(0,12)
	\put( 0, 0){\line(0,1){28}}
	\put(14, 0){\line(0,1){28}}
	\put( 0, 0){\line(1,0){14}}
	\put( 0,14){\line(1,0){14}}
	\put( 0,28){\line(1,0){14}}
	\put( 7,21){\makebox(0,0)[c]{$1$}}
	\put( 7, 7){\makebox(0,0)[c]{$\ol{2}$}}
\end{picture},\quad
$c=$
\begin{picture}(17,28)(0,12)
	\put( 0, 0){\line(0,1){28}}
	\put(14, 0){\line(0,1){28}}
	\put( 0, 0){\line(1,0){14}}
	\put( 0,14){\line(1,0){14}}
	\put( 0,28){\line(1,0){14}}
	\put( 7,21){\makebox(0,0)[c]{$2$}}
	\put( 7, 7){\makebox(0,0)[c]{$\ol{2}$}}
\end{picture},\quad
$d=$
\begin{picture}(17,28)(0,12)
	\put( 0, 0){\line(0,1){28}}
	\put(14, 0){\line(0,1){28}}
	\put( 0, 0){\line(1,0){14}}
	\put( 0,14){\line(1,0){14}}
	\put( 0,28){\line(1,0){14}}
	\put( 7,21){\makebox(0,0)[c]{$2$}}
	\put( 7, 7){\makebox(0,0)[c]{$\ol{1}$}}
\end{picture},\quad
$e=$
\begin{picture}(17,28)(0,12)
	\put( 0, 0){\line(0,1){28}}
	\put(14, 0){\line(0,1){28}}
	\put( 0, 0){\line(1,0){14}}
	\put( 0,14){\line(1,0){14}}
	\put( 0,28){\line(1,0){14}}
	\put( 7,21){\makebox(0,0)[c]{$\ol{2}$}}
	\put( 7, 7){\makebox(0,0)[c]{$\ol{1}$}}
\end{picture}.
\end{center}
\vspace*{12pt}
The tableau in this appendix with the content  
$1^{x_1}2^{x_2}\bar{2}^{\bar{x}_2}\bar{1}^{\bar{x}_1}$ corresponds to 
$(x_1,x_2,\bar{x}_2,\bar{x}_1)$ in 
{\sc Example} \ref{ex:crystalstructure-C-row} or 
{\sc Example} \ref{ex:crystalstructure-C-col}.
In {\sc Table} \ref{tab1}, all the classically restricted paths in 
$B$ with zero classical weight, namely the
$\la=0$ paths, are listed.
The value of $\varepsilon_0$ means, for instance, 
$(\tilde{e}_0)^2 p \ne 0$ and $(\tilde{e}_0)^3 p = 0$ if $\varepsilon_0(p)=2$.
Here it stands for the restriction level.
In view of {\sc Conjecture} \ref{conj:X=M} we have to find suitable
$b_0$ from $B_0=B^{1,2}\ot B^{2,1}\ot B^{1,1}$. 
In this particular example, we further expect that we can take
$B_0=B^{1,2}$ and $b_0=\phi$ (note that $\vphi(\phi)=\La_0$).
We also have $b_0^\natural=\phi,b_1^\natural=11,
b_2^\natural=b_3^\natural=b_4^\natural=a,b_5^\natural=b_6^\natural=1$.
{}From these data, one can write down the expressions of 
the 1dsums $X_l(B,\la =0)$ over the 
level $l=1,2$ restricted paths, and $X(B,\la =0)$ over the
classically restricted paths,
\begin{eqnarray}\label{eq:Xresult}
X_1(B,0,q^{-1}) &=& q^{15}, \nonumber\\
X_2(B,0,q^{-1}) &=& q^{8} + 2 q^{9} + 2 q^{10} + 3 q^{11} + 2 q^{12} + q^{13} 
          + q^{15}, \nonumber\\
X (B,0,q^{-1})&=& q^{6} + 2 q^{7} + 2 q^{8} + 3 q^{9} + 2 q^{10} + 3 q^{11} 
          + 2 q^{12} + q^{13} + q^{15}. \nonumber\\
&&
\end{eqnarray}

We illustrate the calculation of the ``relative'' energy 
$E(p) = \sum_{0 \le i < j \le 6}H(b_i \ot b^{(i+1)}_j)
-\sum_{0 \le i < j \le 6}H(b_i^\natural \ot (b_j^\natural)^{(i+1)})$ 
for a path $p$,
i.e., we include the contribution from $c$ in {\sc Conjecture} \ref{conj:X=M}.
In this example $E(p)$ is reduced to 
\[
E(p)=\sum_{1 \le i < j \le 6}H(b_i \ot b^{(i+1)}_j)
+H(\phi \ot b_1)-H(\phi \ot 11).
\]
We show the computation of the relative energy of the path
$11 \otimes a \otimes d \otimes e \otimes \ol{2} \otimes \ol{1}$. 
The necessary data are shown at the end of this appendix.
It proceeds as 
\begin{eqnarray*}
11 \otimes (a) &\stackrel{0}{\mapsto}& 
(a) \otimes 11 ,
\nonumber \\
11 \otimes a \otimes (d)  &\stackrel{1}{\mapsto}& 
11 \otimes (a) \otimes d  \stackrel{0}{\mapsto}
(a) \otimes 11 \otimes d,
\nonumber \\
11 \otimes a \otimes d \otimes (e)  
&\stackrel{1}{\mapsto}& 
11 \otimes a \otimes (d) \otimes e   
\nonumber \\
&\stackrel{1}{\mapsto}& 
11 \otimes (a) \otimes d \otimes e   
\stackrel{0}{\mapsto}
(a) \otimes 11 \otimes d \otimes e  , 
\nonumber \\
11 \otimes a \otimes d \otimes e \otimes (\ol{2})  
&\stackrel{0}{\mapsto}& 
11 \otimes a \otimes d \otimes (\ol{2}) \otimes e  
\nonumber \\
&\stackrel{0}{\mapsto}&
11 \otimes a \otimes (\ol{1}) \otimes c \otimes e  
\nonumber \\
&\stackrel{1}{\mapsto}&
11 \otimes (1) \otimes d \otimes c \otimes e  
\nonumber \\
&\stackrel{0}{\mapsto}&
(1) \otimes 11 \otimes d \otimes c \otimes e  ,
\nonumber \\
11 \otimes a \otimes d \otimes e \otimes \ol{2}
\otimes   (\ol{1})
&\stackrel{1}{\mapsto}& 
11 \otimes a \otimes d \otimes e \otimes (\ol{2})
\otimes  \ol{1}
\nonumber \\
&\stackrel{0}{\mapsto}& \cdots \stackrel{0}{\mapsto} \cdots 
\stackrel{1}{\mapsto} \cdots \stackrel{0}{\mapsto}
(1) \otimes 11 \otimes d \otimes c \otimes e \otimes  \ol{1}.
\end{eqnarray*}
The parenthesized elements denote $b_j^{(i+1)}$.
They are to be moved  
to the left as $b_i \otimes b^{(i+1)}_j \mapsto b^{(i)}_j \otimes b'_i$ 
by the isomorphism of the crystals.
The values of the energy function $-H(b_i \otimes b_j^{(i+1)})$ are
shown above the arrows, which amount to 6.
We also show the computation of the relative energy
of the path $\phi\otimes a \otimes e 
\otimes a \otimes \ol{2} \otimes \ol{1}$;
\begin{eqnarray*}
% \ol{1}\ol{1} \otimes (\phi) &\stackrel{1}{\mapsto}& 
% (\ol{1}\ol{1}) \otimes \phi ,
% \nonumber \\
\phi \otimes (a) &\stackrel{1}{\mapsto}& 
(a) \otimes 2 \ol{2} ,
\nonumber \\
\phi \otimes a \otimes (e)  &\stackrel{2}{\mapsto}& 
\phi \otimes (a) \otimes e  \stackrel{1}{\mapsto}
(a) \otimes 2 \ol{2} \otimes e,
\nonumber \\
\phi \otimes a \otimes e \otimes (a)  
&\stackrel{0}{\mapsto}& 
\phi \otimes a \otimes (e) \otimes a   
\nonumber \\
&\stackrel{2}{\mapsto}& 
\phi \otimes (a) \otimes e \otimes a   
\stackrel{1}{\mapsto}
(a) \otimes 2 \ol{2} \otimes e \otimes a  , 
\nonumber \\
\phi \otimes a \otimes e \otimes a \otimes (\ol{2})  
&\stackrel{1}{\mapsto}& 
\phi \otimes a \otimes e \otimes (1) \otimes c  
\nonumber \\
&\stackrel{0}{\mapsto}&
\phi \otimes a \otimes (\ol{2}) \otimes c \otimes c  
\nonumber \\
&\stackrel{1}{\mapsto}&
\phi \otimes (1) \otimes c \otimes c \otimes c  
\nonumber \\
&\stackrel{1}{\mapsto}&
(1) \otimes 1 \ol{1} \otimes c \otimes c \otimes c  ,
\nonumber \\
\phi \otimes a \otimes e \otimes a \otimes \ol{2}
\otimes   (\ol{1})
&\stackrel{1}{\mapsto}& 
\phi \otimes a \otimes e \otimes a \otimes (\ol{2})
\otimes  \ol{1}
\nonumber \\
&\stackrel{1}{\mapsto}& \cdots \stackrel{0}{\mapsto} \cdots 
\stackrel{1}{\mapsto} \cdots \stackrel{1}{\mapsto}
(1) \otimes 1 \ol{1} \otimes c \otimes c \otimes c \otimes  \ol{1}.
\end{eqnarray*}
Here the energy amounts to 14.
One should add to it another 1, since 
$-(H(\phi\ot\phi)-H(\phi\ot 11))=1$ (see {\sc Table} \ref{tab5}).

\begin{table}[h]
	\caption{$\la=0$ paths on $B:=
B^{1,2} \otimes (B^{2,1})^{\otimes 3} 
\otimes (B^{1,1})^{\otimes 2}$}
	\label{tab1}
	\begin{center}
	\begin{tabular}{c|c|c} \hline
	{\it path} & $-E$ & $\varepsilon_0$ \\ \hline
 $11 \otimes a \otimes c \otimes e \otimes \ol{1} \otimes \ol{1}$ & 7 &\\
 $11 \otimes a \otimes e \otimes c \otimes \ol{1} \otimes \ol{1}$ & 8 &\\
 $11 \otimes a \otimes d \otimes e \otimes \ol{2} \otimes \ol{1}$ & 6 & 
 $3$ \\
 $11 \otimes a \otimes e \otimes d \otimes \ol{2} \otimes \ol{1}$ & 7 &\\
 $11 \otimes c \otimes a \otimes e \otimes \ol{1} \otimes \ol{1}$ & 9 &\\
	\hline
 $11 \otimes c \otimes c \otimes c \otimes \ol{1} \otimes \ol{1}$ &12 &\\
 $11 \otimes c \otimes c \otimes d \otimes \ol{2} \otimes \ol{1}$ &11 &\\
 $11 \otimes c \otimes d \otimes b \otimes \ol{1} \otimes \ol{1}$ & 9 &\\
 $11 \otimes c \otimes d \otimes e \otimes     1  \otimes \ol{1}$ &11 &\\
 $11 \otimes d \otimes a \otimes e \otimes \ol{2} \otimes \ol{1}$ &10 &\\
 $11 \otimes d \otimes b \otimes c \otimes \ol{1} \otimes \ol{1}$ & 8 &
 $2$ \\
 $11 \otimes d \otimes b \otimes d \otimes \ol{2} \otimes \ol{1}$ & 9 &\\
 $11 \otimes d \otimes e \otimes a \otimes \ol{2} \otimes \ol{1}$ &10 &\\
 $\phi\otimes a \otimes a \otimes e \otimes \ol{2} \otimes \ol{1}$ &11 &\\
 $\phi\otimes a \otimes b \otimes c \otimes \ol{1} \otimes \ol{1}$ &13 &\\
 $\phi\otimes a \otimes b \otimes d \otimes \ol{2} \otimes \ol{1}$ &12 &\\
 \hline
 $\phi\otimes a \otimes e \otimes a \otimes \ol{2} \otimes \ol{1}$ &15 &
 $1$ \\
	\hline
	\end{tabular}
	\end{center}
	\end{table}

Now we consider the fermionic forms.
Let $W := W^{(1)}_2 \otimes (W^{(2)}_1)^{\otimes 3} 
\otimes (W^{(1)}_1)^{\otimes 2}$.
Let us compute
the level $l=1,2$ restricted versions of the fermionic 
form $M_l(W,q^{-1})$ and $M(W,\la =0,q^{-1})$.
First recall that
\begin{equation*}
\left\{ (\alpha_a | \alpha_b) \right\} =
\left(
	\begin{array}{cc}
	1 & -1 \\
	-1 & 2
	\end{array}
\right), \quad
C^{-1} =
\left(
	\begin{array}{cc}
	1 & 1 \\
	\frac12 & 1
	\end{array}
\right), \quad
t_1 =2, \quad t_2 =1.
\end{equation*}
{}From these data
the constraint (\ref{eq:mlq5}) on the summation variables $m$'s
for $W$ reads
\begin{equation*}
\sum_{i=1}^{2l} i m^{(1)}_i=7 ,\quad \sum_{i=1}^{l} i m^{(2)}_i=5.
\end{equation*}
Among the 105 possible configurations of the $m$'s, only those 6 
listed in {\sc Table} \ref{tab2} have non-zero contributions to the 
fermionic forms.
Consequently,
\begin{eqnarray*}
M_1(W,q^{-1}) &=& q^{15}, \nonumber\\
M_2(W,q^{-1}) &=& q^{8} + 2 q^{9} + 2 q^{10} + 3 q^{11} + 2 q^{12} + q^{13} 
          + q^{15}, \nonumber\\
M(W,\la =0,q^{-1})&=& q^{6} + 2 q^{7} + 2 q^{8} + 3 q^{9} + 2 q^{10} + 3 q^{11} 
          + 2 q^{12} + q^{13} + q^{15} .\nonumber\\
&&
\end{eqnarray*}
They are exactly the same as (\ref{eq:Xresult}) calculated as the 1dsums.
\begin{table}[h]
	\caption{Relevant configurations to the fermionic forms}
	\label{tab2}
	\begin{center}
\begin{tabular}{l|l|c} \hline
 $m^{(1)},m^{(2)}$ &  $p^{(1)},p^{(2)}$ &contribution  \\ \hline
$(1,0,0,0,0,1), (0,1,1)$ & $(1,2,2,2,1,0), (2,0,0)$ & $q^6 + q^7$ \\
$(1,0,0,0,0,1), (2,0,1)$ & $(2,4,3,2,1,0), (0,0,0)$ & $q^7 + q^8 + q^9$ \\ \hline
$(0,0,1,1), (1,2)$ & $(2,2,0,0), (1,0)$ & $q^8 + q^9$ \\
$(1,1,0,1), (1,2)$ & $(0,0,0,0), (2,0)$ & $q^9 + q^{10} + q^{11}$ \\
$(1,1,0,1), (3,1)$ & $(1,2,1,0), (0,0)$ & $q^{10} + 2 q^{11} + 2 q^{12} + q^{13}$ \\ \hline
$(1,3), (5)$ & $(0,0), (0)$ & $q^{15}$ \\ \hline
\end{tabular}
\end{center}
\end{table}

Below we list the energy function and the
isomorphism (combinatorial R-matrix) of the 
$U'_q(C_2^{(1)})$ crystals
used in calculating the 1dsums.
The data on $B^{1,1} \otimes B^{1,1}$, $B^{2,1} \otimes B^{2,1}$ and
$B^{2,1} \otimes B^{1,1}$ are quoted from  \cite{Y}.
Here we added the cases including $B^{1,2}$.
For $B^{1,1} \otimes B^{1,1}$, $B^{2,1} \otimes B^{2,1}$ and 
$B^{1,2} \otimes B^{1,2}$ cases, the isomorphism is 
the trivial (identity) map.
The blanks in the tables signify that the energy is $0$.

\begin{table}[h]
	\caption{$-H(b_1 \otimes b_2)$ on $B^{1,1} \otimes B^{1,1}$}
	\label{tab3}
    \newcommand{\lw}[1]{\smash{\lower1.6ex\hbox{#1}}}
	\begin{center}
	\begin{tabular}{c|cccc}	\hline
	\lw{$b_1$} & \multicolumn{4}{c}{$b_2$} \\
% 	\cline{2-5}
	& $1$ & $2$ & $\ol{2}$ & $\ol{1}$ \\ \hline
	                  $1$  &   & $1$ &      $1$ &  $1$     \\
	                  $2$  &   &     &      $1$ &  $1$     \\
	              $\ol{2}$ &   &     &          &  $1$     \\
	              $\ol{1}$ &   &     &          &        \\ \hline
	\end{tabular}
	
	\end{center}
\end{table}

\begin{table}[h]
	\caption{$-H(b_1 \otimes b_2)$ on $B^{2,1} \otimes B^{2,1}$}
	\label{tab4}
    \newcommand{\lw}[1]{\smash{\lower1.6ex\hbox{#1}}}
	\begin{center}
	\begin{tabular}{c|ccccc}	\hline
	\lw{$b_1$} & \multicolumn{5}{c}{$b_2$} \\
	  & $a$ & $b$ & $c$ & $d$ & $e$ \\ \hline
	                $a$  &     & $1$ & $1$ & $1$ & $2$ \\
	                $b$  &     &     & $1$ & $1$ & $1$ \\
	                $c$  &     &     & $1$ & $1$ & $1$ \\
	                $d$  &     &     &     &     & $1$ \\
	                $e$  &     &     &     &     &     \\ \hline
	\end{tabular}
	
	\end{center}
\end{table}

\begin{table}[h]
	\caption{$-H(b_1 \otimes b_2)$ on $B^{1,2} \otimes B^{1,2}$}
	\label{tab5}
    \newcommand{\lw}[1]{\smash{\lower1.6ex\hbox{#1}}}
	\begin{center}
	\begin{tabular}{c|ccccccccccc}	\hline
	\lw{$b_1$} & \multicolumn{11}{c}{$b_2$} \\
 & $11$ & $12$ & $1\ol{2}$ & $1\ol{1}$ & $22$ & $2\ol{2}$ 
& $2\ol{1}$ & $\ol{2}\ol{2}$ & $\ol{2}\ol{1}$ & $\ol{1}\ol{1}$ & $\phi$ \\
\hline
$11$           &   & 1 & 1 & 2 & 2 & 2 & 2 & 2 & 2 & 2 & 1 \\
$12$           &   & 1 & 1 & 2 & 1 & 2 & 2 & 2 & 2 & 2 & 1 \\
$1\ol{2}$      &   & 1 & 1 & 2 & 1 & 1 & 2 & 1 & 2 & 2 & 1 \\
$1\ol{1}$      &   & 1 & 1 & 2 & 1 & 1 & 2 & 1 & 2 & 2 & 1 \\
$22$           &   &   & 1 & 1 &   & 2 & 1 & 2 & 2 & 2 & 1 \\
$2\ol{2}$      &   &   & 1 & 1 &   & 2 & 1 & 2 & 2 & 2 & 1 \\
$2\ol{1}$      &   &   & 1 & 1 &   & 1 & 1 & 1 & 1 & 1 & 1 \\
$\ol{2}\ol{2}$ &   &   &   & 1 &   &   & 1 &   & 1 & 2 & 1 \\
$\ol{2}\ol{1}$ &   &   &   & 1 &   &   & 1 &   & 1 & 1 & 1 \\
$\ol{1}\ol{1}$ &   &   &   &   &   &   &   &   &   &   & 1 \\
$\phi$         & 1 & 1 & 1 & 1 & 1 & 1 & 1 & 1 & 1 & 1 & 2 \\ \hline
	\end{tabular}
	
	\end{center}
\end{table}

\begin{table}[h]
	\caption{Isomorphism and 
	$-H(b_1 \otimes b_2)$ on $B^{2,1} \otimes B^{1,1}$}
	\label{tab6}
    \newcommand{\lw}[1]{\smash{\lower1.6ex\hbox{#1}}}
	\begin{center}
	\begin{tabular}{c|cccc||cccc}	\hline
	\lw{$b_1$} & \multicolumn{4}{c||}{$b_2$} & 
	\multicolumn{4}{c}{$b_2$} \\
	  & $1$ & $2$ & $\ol{2}$ & $\ol{1}$ 
	                     & $1$ & $2$ & $\ol{2}$ & $\ol{1}$ \\ \hline
$a$  & $1 \otimes a$ & $2 \otimes a$ & $1 \otimes c$ & $1 \otimes d$ 
     &               &               &             1 &             1\\
$b$  & $1 \otimes b$ & $\ol{2} \otimes a$ 
& $\ol{2} \otimes b$ & $1 \otimes e$ 
     &               &               &               &             1\\
$c$  & $2 \otimes b$ & $\ol{1} \otimes a$ 
& $\ol{1} \otimes b$ & $2 \otimes e$ 
     &               &               &               &             1\\
$d$  & $2 \otimes c$ & $2 \otimes d$ 
& $\ol{1} \otimes c$ & $\ol{1} \otimes d$ 
     &               &               &               &              \\
$e$  & $\ol{2} \otimes c$ & $\ol{2} \otimes d$ 
& $\ol{2} \otimes e$ & $\ol{1} \otimes e$ 
     &               &               &               &              \\ 
    \hline
	\end{tabular}
	
	\end{center}
\end{table}

\begin{table}[h]
	\caption{Isomorphism and 
	$-H(b_1 \otimes b_2)$ on $B^{1,2} \otimes B^{1,1}$}
	\label{tab7}
    \newcommand{\lw}[1]{\smash{\lower1.6ex\hbox{#1}}}
	\begin{center}
	\begin{tabular}{c|cccc||cccc}	\hline
	\lw{$b_1$} & \multicolumn{4}{c||}{$b_2$} & 
	\multicolumn{4}{c}{$b_2$} \\
	 & $1$ & $2$ & $\ol{2}$ & $\ol{1}$ 
	                     & $1$ & $2$ & $\ol{2}$ & $\ol{1}$ \\ \hline
$11$  & $1 \otimes 11$       & $1 \otimes 12$ 
      & $1 \otimes 1\ol{2}$  & $1 \otimes \phi$ 
     &               &             1 &             1 &             1\\
$12$  & $2 \otimes 11$       & $1 \otimes 22$ 
      & $2 \otimes 1\ol{2}$  & $2 \otimes \phi$ 
     &               &             1 &             1 &             1\\
$1\ol{2}$  & $\ol{2} \otimes 11$        & $1 \otimes 2\ol{2}$ 
           & $1 \otimes \ol{2}\ol{2}$   & $\ol{2} \otimes \phi$ 
     &               &             1 &             1 &             1\\
$1\ol{1}$  & $\ol{1} \otimes 11$        & $2 \otimes 2\ol{2}$ 
           & $2 \otimes \ol{2}\ol{2}$   & $\ol{1} \otimes \phi$ 
     &               &             1 &             1 &             1\\
$22$       & $2 \otimes 12$        & $2 \otimes 22$ 
           & $2 \otimes 1\ol{1}$   & $2 \otimes 2\ol{1}$ 
     &               &               &             1 &             1\\
$2\ol{2}$  & $\ol{2} \otimes 12$          & $\ol{2} \otimes 22$ 
           & $\ol{2} \otimes 1\ol{1}$     & $\ol{2} \otimes 2\ol{1}$ 
     &               &               &             1 &             1\\
$2\ol{1}$  & $\ol{1} \otimes 12$          & $\ol{1} \otimes 22$ 
           & $2 \otimes \ol{2}\ol{1}$     & $2 \otimes \ol{1}\ol{1}$ 
     &               &               &             1 &             1\\
$\ol{2}\ol{2}$  & $\ol{2} \otimes 1\ol{2}$     & $\ol{2} \otimes 2\ol{2}$ 
         & $\ol{2} \otimes \ol{2}\ol{2}$ & $\ol{2}  \otimes \ol{2}\ol{1}$ 
     &               &               &               &             1\\
$\ol{2}\ol{1}$  & $\ol{1} \otimes 1\ol{2}$     & $\ol{1} \otimes 2\ol{2}$ 
         & $\ol{1} \otimes \ol{2}\ol{2}$ & $\ol{2}  \otimes \ol{1}\ol{1}$ 
     &               &               &               &             1\\
$\ol{1}\ol{1}$  & $\ol{1} \otimes 1\ol{1}$     & $\ol{1} \otimes 2\ol{1}$ 
         & $\ol{1} \otimes \ol{2}\ol{1}$ & $\ol{1}  \otimes \ol{1}\ol{1}$ 
     &               &               &               &              \\
$\phi$  & $1 \otimes 1\ol{1}$     & $1 \otimes 2\ol{1}$ 
         & $1 \otimes \ol{2}\ol{1}$ & $1  \otimes \ol{1}\ol{1}$ 
     &            1  &            1  &            1  &           1  \\
    \hline
	\end{tabular}
	
	\end{center}
\end{table}

\begin{table}[h]
	\caption{Isomorphism and 
	$-H(b_1 \otimes b_2)$ on $B^{1,2} \otimes B^{2,1}$}
	\label{tab8}
	\begin{center}
    \newcommand{\lw}[1]{\smash{\lower1.6ex\hbox{#1}}}
	\begin{tabular}{c|ccccc||ccccc}	\hline
	\lw{$b_1$} & \multicolumn{5}{c||}{$b_2$} & 
	\multicolumn{5}{c}{$b_2$} \\
	 & $a$ & $b$ & $c$ & $d$ & $e$ 
	                     & $a$ & $b$ & $c$ & $d$ & $e$ \\ \hline
$11$  & $a \otimes 11$       & $b \otimes 11$ & $a \otimes 1\ol{2}$ 
      & $a \otimes \phi$  & $b \otimes \phi$ 
     &    &    &  1 &  1 &  1\\
$12$  & $a \otimes 12$       & $c \otimes 11$ & $a \otimes 1\ol{1}$
      & $a \otimes 2\ol{1}$  & $c \otimes \phi$ 
     &    &    &  1 &  1 &  1\\
$1\ol{2}$  & $b \otimes 12$ & $b \otimes 1\ol{2}$ & $b \otimes 1\ol{1}$ 
           & $b \otimes 2\ol{1}$   & $b \otimes \ol{2}\ol{1}$ 
     &    &    &  1 &  1 &  1\\
$1\ol{1}$  & $c \otimes 12$ & $c \otimes 1\ol{2}$ & $c \otimes 1\ol{1}$ 
           & $c \otimes 2\ol{1}$   & $c \otimes \ol{2}\ol{1}$  
      &    &    &  1 &  1 &  1\\
$22$       & $a \otimes 22$ & $d \otimes 11$ & $d \otimes 12$ 
           & $d \otimes 22$ & $d \otimes \phi$ 
      &    &    &    &    &  1\\
$2\ol{2}$  & $b \otimes 22$ & $e \otimes 11$ & $e \otimes 12$ 
           & $e \otimes 22$ & $e \otimes \phi$ 
      &    &    &    &    &  1\\
$2\ol{1}$  & $c \otimes 22$ & $d \otimes 1\ol{2}$ & $d \otimes 1\ol{1}$ 
           & $d \otimes 2\ol{1}$   & $c \otimes \ol{1}\ol{1}$  
      &    &    &    &    &  1\\
$\ol{2}\ol{2}$  
  & $b \otimes 2\ol{2}$ & $b \otimes \ol{2}\ol{2}$ & $e \otimes 1\ol{2}$ 
           & $e \otimes 2\ol{2}$   & $e \otimes \ol{2}\ol{2}$  
      &    &    &    &    &   \\
$\ol{2}\ol{1}$  
  & $c \otimes 2\ol{2}$ & $c \otimes \ol{2}\ol{2}$ & $e \otimes 1\ol{1}$ 
           & $e \otimes 2\ol{1}$   & $e \otimes \ol{2}\ol{1}$  
      &    &    &    &    &   \\
$\ol{1}\ol{1}$  
  & $d \otimes 2\ol{2}$        & $d \otimes \ol{2}\ol{2}$ 
  & $d \otimes \ol{2}\ol{1}$ 
  & $d \otimes \ol{1}\ol{1}$   & $e \otimes \ol{1}\ol{1}$  
      &    &    &    &    &   \\
$\phi$  
  & $a \otimes 2\ol{2}$        & $a \otimes \ol{2}\ol{2}$ 
  & $a \otimes \ol{2}\ol{1}$ 
  & $a \otimes \ol{1}\ol{1}$   & $b \otimes \ol{1}\ol{1}$  
     & 1  &  1 &  1 &  1  &  1 \\   \hline
	\end{tabular}
	
	\end{center}
\end{table}

%%%%%%%%%%%%%%%%%%%%%%%%%%%%%%%%%%%%%%%%%%%%%%
\clearpage
\section{Explicit form of $Q$-system}\label{app:qsys}
Let us write down the $Q$-system for 
non-simply laced $X_n$.
Below we assume that $Q^{(a)}_j = 1$ whenever 
$a \not \in \{1, \ldots, n\}$.

\noindent
$X_n=B_n:$
\begin{eqnarray*}
Q^{(a)^2}_j
  &=&Q_{j-1}^{(a)}Q_{j+1}^{(a)} + Q^{(a-1)}_j Q^{(a+1)}_j
\quad (1 \le a \le n-2 ), \\
Q^{(n-1)^2}_j &=& 
Q_{j-1}^{(n-1)} Q_{j+1}^{(n-1)} + Q^{(n-2)}_j Q^{(n)}_{2j}, \\
Q^{(n)^2}_{2j}
&=& Q_{2j-1}^{(n)} Q_{2j+1}^{(n)}
+ Q^{(n-1)^2}_j, \\
Q^{(n)^2}_{2j+1}
&=& Q_{2j}^{(n)}Q_{2j+2}^{(n)} 
+ Q^{(n-1)}_j Q^{(n-1)}_{j+1}.  \\
\end{eqnarray*}
$X_n=C_n:$
\begin{eqnarray*}
 Q^{(a)^2}_j &=& Q_{j-1}^{(a)} Q_{j+1}^{(a)}
 + Q^{(a-1)}_j Q^{(a+1)}_j 
\quad  (1 \le a \le n-2), \nonumber \\
Q^{(n-1)^2}_{2j}
      &=& Q_{2j-1}^{(n-1)} Q_{2j+1}^{(n-1)} 
+ Q^{(n-2)}_{2j}Q^{(n)^2}_j, \nonumber \\
Q^{(n-1)^2}_{2j+1}
&=& Q_{2j}^{(n-1)} Q_{2j+2}^{(n-1)}
+ Q^{(n-2)}_{2j+1} Q^{(n)}_j Q^{(n)}_{j+1}, \nonumber \\
Q^{(n)^2}_j
  &=& Q_{j-1}^{(n)} Q_{j+1}^{(n)}
+ Q^{(n-1)}_{2j}. \nonumber \\
\end{eqnarray*}
$X_n=F_4:$
\begin{eqnarray*}
 Q^{(1)^2}_j
  &=& Q_{j-1}^{(1)} Q_{j+1}^{(1)}
+ Q^{(2)}_{j}, \nonumber \\
Q^{(2)^2}_j
  &=& Q_{j-1}^{(2)} Q_{j+1}^{(2)}
+ Q^{(1)}_{j}Q^{(3)}_{2j}, \nonumber \\
Q^{(3)^2}_{2j}
      &=& Q_{2j-1}^{(3)} Q_{2j+1}^{(3)}
+ Q^{(2)^2}_j Q^{(4)}_{2j},\nonumber \\
Q^{(3)^2}_{2j+1}
      &=& Q_{2j}^{(3)} Q_{2j+2}^{(3)}
+ Q^{(2)}_j Q^{(2)}_{j+1} Q^{(4)}_{2j+1}, \nonumber \\
Q^{(4)^2}_j
&=& Q_{j-1}^{(4)} Q_{j+1}^{(4)}
+ Q^{(3)}_{j}. \nonumber \\
\end{eqnarray*}
$X_n=G_2:$
\begin{eqnarray*}
Q^{(1)^2}_j
  &=& Q_{j-1}^{(1)} Q_{j+1}^{(1)} + Q^{(2)}_{3j}, \nonumber \\
Q^{(2)^2}_{3j}
  &=& Q_{3j-1}^{(2)} Q_{3j+1}^{(2)}
+ Q^{(1)^3}_j, \nonumber \\
Q^{(2)^2}_{3j+1}
  &=& Q_{3j}^{(2)} Q_{3j+2}^{(2)}
+ Q^{(1)^2}_j Q^{(1)}_{j+1},\nonumber \\
Q^{(2)^2}_{3j+2}
  &=& Q_{3j+1}^{(2)} Q_{3j+3}^{(2)}
+ Q^{(1)}_{j} Q^{(1)^2}_{j+1}. \nonumber \\
\end{eqnarray*}

%%%%%%%%%%%%%%%%%%%%%%%%%%%%%%%%%%%%%%%%%%%%%

\end{document}